\documentclass[11pt,a4paper]{article}
\usepackage{amsmath}
\usepackage{latexsym}
\usepackage{theorem}
\usepackage[all]{xy}
\newtheorem{teorema}{Theorem}[section]
\newtheorem{definicion}[teorema]{Definition}
\include{amslatex}
\newtheorem{proposicion}[teorema]{Proposition}
\newtheorem{lema}[teorema]{Lemma}
\newtheorem{corolario}[teorema]{Corollary}
\textheight
21 cm \textwidth 14 cm \oddsidemargin 10 mm \topmargin  0 mm
{\theorembodyfont{\rmfamily}
\newtheorem{comentario}[teorema]{Remark}}
{\theorembodyfont{\rmfamily}
\newtheorem{ejemplo}[teorema]{Example}}

\numberwithin{equation}{section}

\include{amsmath}
\begin{document}

\title{\Large {\bf Geometry of generalized higher order fields and applications to classical linear electrodynamics}}

\maketitle
\author{

\begin{center}

Ricardo Gallego Torrom\'{e}\\
Departamento de Matem\'{a}tica\\
 Universidade de S\~{a}o Paulo \footnote{email:
rigato39@gmail.com; Currently at the Departamento de Matem\'atica, Universidade Federal de S\~{a}o Carlos, Brazil.}\\ [3pt]
\end{center}}
\begin{abstract}
{\small Motivated by obtaining a consistent mathematical description for the radiation reaction of point charged particles in linear classical electrodynamics,
a theory of {\it generalized higher order tensors and differential forms} is introduced. The generalization of some fundamental notions of the differential geometry and the theory of differential forms is presented. In particular,
 the cohomology and integration theories for generalized higher order forms are developed, including the Cartan calculus, a generalization of de Rham cohomology and a version of Thom's isomorphism theorem. We consider in detail a special type of generalized higher order tensors associated
  with {\it bounded maximal $n$-acceleration} and use it as a model of spacetime.
  A generalization of electrodynamic theory with higher order fields is introduced. Although the theory is non-local in the usual sense, it is free of some of the pathologies appearing in the standard linear classical electrodynamics. Indeed, we show that combining the generalized higher order fields with {\it maximal acceleration geometry}
the evolution of a point charged particle interacting with the generalized higher order fields can be described by solutions of an implicit
second order ordinary differential equation. In flat space such equation is Lorentz invariant, does not have pre-accelerated solutions of Dirac's type or run-away solutions, it is compatible with Newton's first law of dynamics and with the covariant Larmor's power radiation law. A generalization of the Maxwell-Lorentz theory is also introduced. The theory is linear in the field sector and it reduces to the standard Maxwell-Lorentz electrodynamics when the maximal acceleration is infinite. Finally, we discuss the assumptions of our framework in addition to some predictions of the theory. Apart from the many open questions already leave in this work, we indicate further research directions, including the full development of the cohomology theory of generalized forms, its relation with calibrated geometry and a theory of curvature of generalized metrics. From the physical side, we emphasize the extension to non-linear Yang-Mills theory and to gravity as well as the problem of quantizing  theories with generalized higher order fields.}
\end{abstract}
\tableofcontents{}
\newpage

\section{Introduction}
\subsection{Motivation}

The theory of classical electrodynamics of point charged particles suffers from severe theoretical problems
when one considers the coupled dynamics of charged particles with the total external and its own radiation electromagnetic field.
In such circumstances, the standard theory of classical electrodynamics, based on Maxwell's equations and the Lorentz-Dirac equation, makes un-physical predictions \cite{Jackson}.  A main theoretical
 difficulty is that although the derivation of the standard theory  is based on general principles, the
 equation has run-away (solutions with acceleration un-bounded, even if there is no external force) and pre-accelerated solutions (depending of the future sector of the world-line of the particle). Therefore, the problematic situation on which the classical electrodynamic theory lies permanently and the practical implications that a consistent theory of radiation reaction may be found in accelerator science and in modeling very dense, non-neutral, relativistic plasmas, demonstrates the relevance of having a consistent classical theory of charged particles and radiation reaction.

One can argue that such problems are of less relevance for fundamental physics, since from the advent of quantum field theory,
classical field theories describing fundamental interactions should be considered as the
{\it classical limit} of an appropriate quantum field theory. Following this line of thought, one could expect
that the unwanted problems of the classical electromagnetic theory (divergent Coulomb fields,
run away solutions and pre-accelerated solutions of the Lorentz-Dirac force equation
\cite{Dirac}) are cured in the framework of a convenient quantum theory. This point of
view is supported by the fact that some of the unwanted effects of the classical theory are believed to be
significant at scales where quantum effects become relevant.

Despite the fact that such problems make the classical theory unsatisfactory from a
theoretical viewpoint, the development of a quantum field theory able to overcome
them is yet to be realized. The situation is further complicated by the absence of
a well-defined regime of validity for the classical theory; without the security
of a \emph{valid} standard classical electromagnetic theory up to the energy scale
where quantum effects become relevant, a top-down approach to constructing a
quantum theory is inherently problematic. These comments specially apply to the
work of Moniz and Sharp \cite{MonizSharp}, where several hypothesis are not clearly justified
(in particular, the asymptotic condition of the states and some requirements on
analysity of the solutions). Thus, we consider that the conclusions drawn in
\cite{MonizSharp} are only of partial validity and under too restrictive conditions.

One proposed solution to the foundational problems of classical electromagnetism was detailed by
 Landau and Lifshitz \cite{LandauLifshitz} and cast in the framework of singular perturbation theory
(and developed thereafter) by H. Spohn and co-workers \cite{Spohn1, Spohn2}). Despite the success of  this second order
 differential equation in describing the dynamics of a point charged particle (instead of the third
order Dirac equation obtained from Dirac's theory \cite{Dirac}), the theory is still laden with theoretical difficulties. In particular,
the starting point in Landau-Lifshitz theory is the Lorentz-Dirac equation, which is substituted
by a more convenient yet {\it equivalent} equation. This procedure is seemingly ad-hoc
(one starts from a differential equation which is physically unacceptable) and does not
appeal to {\it fundamental principles}. Furthermore, it is unclear whether pre-accelerated
solutions can be obtained from the Landau-Lifshitz equation. Although these deficiencies
 may be overcome by future theoretical developments, one cannot be sure a priori whether
this will be achieved within the framework of standard classical electromagnetic theory.

Another approach with the of
flavor the standard classical theory was introduced more than an century by Larmor and reviewed later
by W. Bonnor \cite{Bonnor} (see also \cite{Herrera1990a}). The main idea was that {\it the observable rest mass $m$} of a  point charged particle can vary with time, providing the origin of the energy-momentum radiated and that in effective terms eliminates the Schott's term in the Lorentz-Dirac equation. The main difficulty with Bonnor's theory relies on the physical motivation and the interpretation of a variable mass for a point charged particle, that usually is interpreted as a {\it fundamental particle}. Although Bonnor's proposal is consistent with experiment, the idea of varying mass for elementary particles is disfavored by the complications in the interpretation of fundamental particle, for example as labeled by a irreducible representation of the Poincar\'e's group.

 T. S. Mo and C. H. Papas introduced the idea that the external electromagnetic force on a point charged particle could be dissipative \cite{MoPapas}. Thus, a modification of the Lorentz force in such a way that a term proportional to the acceleration was introduced. However, apart from the required non-local modification of the energy-momentum tensor, the theory is lacking of a clear foundation from fundamental principles. In Mo-Papas' theory it is unclear the origin of the radiation term as well as the particular form of the  generalized Lorentz force. Nevertheless, the proposed differential equation does not have pre-acceleration, run away solutions and is compatible with radiation damping.

A very different approach to the foundations of classical electrodynamics was developed by R. P. Feynman and J. A. Wheeler \cite{FeynmanWheeler}, based on previous work of K. Schwarzschild, H. Tetrode and V. Fokker.
For instance, in the Wheeler-Feynman theory, the fields are still living on the manifold $M$ and the particle-field interaction is formulated in a time-symmetric way. As a result, the radiation damping term is the same than in the standard Abraham-Lorentz theory. Such damping term it is the problematic one. Thus, Wheeler-Feynman's appears more as a justification for the theory of Abraham-Lorentz-Dirac than a solution to the problem.

Without being exhaustive in the review of the different ideas and theories aiming to solve the
problems of classical electrodynamics, it is clear that the current frameworks
 are not entirely satisfactory. The experimental difficulties in testing
fundamental postulates of theories of classical electrodynamics, particularly in relation
to the problems described before, partially explain the lack of a solution to the
foundational problems of the theory. It is therefore reasonable to
investigate new perspectives and investigate the possible new experimental consequences.

In this work we present a theory of generalized higher order jet bundle tensors and differential forms  and we apply
it to classical electrodynamics, in order to find a model of electrodynamics of point charged particles free from the problems
of run away solutions and pre-accelerated solutions of Dirac's type\footnote{The point particle will be considered structureless, therefore disregarding spin, etc...}. Motivated by experimental testability and economy of postulates, the
formalism for our theory is based on a {\it phenomenological description} of the
electromagnetic field $F$, as determined by its effect on the motion of point charged particles.
Point charged particles are characterized by certain smooth curves on the spacetime $M$
that are interpreted as the world-line curves of point charged particles interacting with the electromagnetic
field that we want to measure. It is the freedom that one has in the description of the electromagnetic field that in combination with the bound on $n$-acceleration sufficient to find a consistent field-particle dynamics. Each of the both new elements alone does not provide a consistent theory.

The theory proposed does not aim to solve the problem
of the singularities of the electromagnetic fields, yet we will use the more
phenomenological technique of renormalization of mass. In this way, although the problems
 of the infinities in classical electrodynamics is a difficult one, we show that it is possible to understand other problems. Indeed,
we believe that
 the problem of infinities in electrodynamics is of different nature than the run-away and pre-accelerated solutions and probably requires a
 better understanding of the ultra-violet limits of classical field theories.
  Possible avenues to solve the divergent problems could be higher order modifications of the non-linear electrodynamics of Born and Infeld  \cite{BornInfeld}, Bopp-Podolsky theory \cite{Bopp, Podolsky} or action at a distance theories on the way of Wheeler-Feynman theory. However, we do not consider that topic here and we show that indeed one can have an effective, consistent description of classical electrodynamics, if one consider maximal acceleration geometry structure, generalized higher order fields and mass renormalization.

\subsection{Criticism of the notions of external electromagnetic field and point test particle}

In the description of point charged {\it test particle} dynamics influenced by the
action of an external field, an strategy that one can adopt is the
following. One starts with the hypothesis that the {\it test particle}  only `see'
the external field that we wish to measure (called the {\it external field}). The standard description of the trajectory of
the test particle is then obtained as a solution to the Lorentz force equation, determined
by the external field. Such a na\"{i}ve approach, based on experimental determination of
the evolution of the test particle, would however yield an inaccurate description of the external field.
This is because an accelerating charged particle radiates electromagnetic waves, and this
radiation of energy-momentum would have a twofold effect on the motion of the test particle.
Firstly, the \emph{total} electromagnetic field would now be a linear combination of the
external and radiation fields, as opposed to the sole external field. Secondly, the radiation
of energy by the test particle would change its energy-momentum four-vector, thereby
producing a dynamical effect on the motion of the test particle. Since the motion of
the test particle is determined by the total field (rather than the sole external field),
this suggests that the above description of the motion of a test particle is not accurate enough to describe
the behavior of point charged in external electromagnetic fields.

Let us focus our attention initially on the related issue of how an electromagnetic field
is defined in terms of physical measurements on the motion of point charged particles
(note we have replaced the denomination of `test particle' with `point charged particle').
The discussion of the previous paragraph indicates that the motion of a point charged particle
is affected not only by the external field, but also by the radiation of the particle
emitted as it accelerates. Therefore, the notion of defining an {\it external electromagnetic field}
 based on the measurement of the motion of test particles is not as clear-cut as we would like,
 since it cannot be characterized by measurements on the motion of charged particles without a
 complete description of the self-field. Both the point charged particle and the total field form a
 coupled dynamical system. In such systems it is difficult to abstract a notion of an external
 electromagnetic field that is consistent with a phenomenological characterization.
  Further, we observe that it is natural to interpret the standard theory
 of classical electrodynamics as an effective theory. With this interpretation in mind, the
 notion of an external field must appear along with an associated notion of a test particle.
 In addition, a covariant notion of weak-strong fields should be possible
 \footnote{The pair (particle, external field) is a convenient representation
 if the external field has a much stronger effect than the radiative-reaction
 in the motion of the test particle. A notion for {\it strong} or {\it week}
 field can be introduced in a covariant way using the norm-operator associated
 with the Riemannian metric determined by an observer \cite{HawEllis}.}.

\subsection{A theory of generalized higher order electromagnetic fields}

In order to describe the classical electrodynamics of point charged particles interacting with
electromagnetic fields, a modification of the notion of electromagnetic field can
be useful. We propose that the description of an electromagnetic field should depend
not only on the macroscopic source, but also on the state of motion of the charged
particle which is used in each particular measurement setting. The motion of the
particle is not prescribed a priori, but we assumed that such evolution exists and that it is regular in
the domain of definition, assumptions that in the classical domain are reasonable.
We then follow {\it fundamental physical and geometric principles} to obtain a mathematical formalism for
 the electromagnetic field, the equation of motion of the point charged particle and
 the equation of evolution for the fields. The guiding principles when pursuing this strategy are as follows:

\begin{enumerate}

\item A minimal higher order extension of the notion of the field. We think that a formalism capable
of accommodating the dynamical system of point charged particles interacting with electromagnetic
fields can be constructed if the electromagnetic and other fields are described as
sections of certain sub-bundles of higher order jet bundles over the spacetime manifold $M$.
In particular, the generalization of the notion of a field introduces degrees of freedom in
such a way that one can obtain consistent dynamics of particles and fields. This procedure
is performed in a minimal way, introducing the minimal number of new degrees of freedom,
attempting to be as conservative as possible in the process of the generalization. In this case,
the generalized electromagnetic field lives in the third jet $J^3_0(M)$.

\item {\it Geometry of maximal acceleration}. In the Theory of General Relativity, the geometry of the
spacetime is dynamical. On the other hand, if the modified electromagnetic field enters
in Einstein's equations, then in an analogous way the substitute for the metric field needs
 to be formulated in the same framework of higher order jet bundles. Therefore, the
  substitute of the spacetime metric tensor will be a tensor defined on higher order
   bundles over the spacetime manifold $M$. One way to achieve this is by using {\it geometries
   of maximal acceleration}. This provides a minimum extension such that the principles
   third to five described bellow are fulfilled. Furthermore, a description of geometries with maximal covariant acceleration
   requires the use of generalized metrics, with coefficients living on the jet bundle $J^2_0(M)$.
   This can be seen as a justification of the formalism of generalized higher order fields.

\item We assume that the theory is an effective theory, in the sense that it depends on a small
parameter in such a way that in the limit when such parameter goes to zero, one obtains the standard
classical Maxwell-Lorentz theory. This parameter is related to the inverse square of the maximal
covariant acceleration described in the point before. However, we do not discuss a particular
mechanism producing the maximal acceleration and only general considerations are presented.

\item Conservation of the energy-momentum of the system. This implies that the loss of energy-momentum
of the particle must be compensated by its emitted radiation, following a covariant Larmor formula \cite{Jackson, Rohrlich}. That the covariant Larmor's law is still valid is currently an assumption. However, it is consistent with the fact that Maxwell's equation is a good effective description of the electromagnetic fields in the regime that accelerations are small compared with the maximal acceleration.

\item The equation of motion of the point charged particles must be a second order ordinary differential equation.
This is a requirement that we follow, in order to avoid the problems that plague
the standard formulation of classical mechanics. This principle is very restrictive, as we will see.

\end{enumerate}
We construct a theory where all the physical fields (electromagnetic excitation, electromagnetic field,
density currents and generalized metrics in the case of electrodynamics) have coefficients living on higher order jet bundles over the
spacetime manifold $M$ \footnote{An alternative description for generalized higher order fields  is as sections of the bundle of forms
$\Lambda^p(J^k_0(M))$ that are horizontal, in the sense that will be indicated later. Fortunately both descriptions are equivalent.}. A
maximal covariant $n$-acceleration is introduced, providing a book-keeping parameter for the theory.
We show that in the limit of large maximal $n$-acceleration, there is a convenient definition
of external field and test particle which recovers the standard notions.

As the result of our analysis of the extension of electromagnetic fields to higher order jet bundles
on the spacetime manifold $M$ and with point particles moving with bounded speed and acceleration,
we derive an effective theory living on $M$ for point charged particles and electromagnetic
fields such that it does not contain some of the pathologies of the standard classical theory.
In particular, the electromagnetic fields are described by a set of equations analogous to Maxwell's equations, whereas the
point charged particle is described by the covariant equation $(\ref{covariantequationofmotion})$
(or in its normal coordinate version form, equation. (\ref{equationofmotion})), which is of
second order and does not suffer most of the problems found with the Lorentz-Dirac equation.
The effective full dynamics of fields and point charged particles is described in a general covariant way by equations
$(\ref{covariantequationofmotion})$ and the generalized Maxwell equations $(\ref{equationforF}),
\,(\ref{equationfor*F})$ together with the charge conservation law $(\ref{equationfor J})$.

Our considerations will involve localized point charged particles with fixed charge $e$ and mass $m$.
The electromagnetic media will be the {\it vacuum}, and therefore the permittivity
and permeability are constant fields living on the spacetime manifold $M$. The study of other
more convenient media and charge distributions will be postponed to subsequent investigations.

\subsection{Structure of the work}

In this paper we explore a classical electrodynamic theory with electromagnetic fields to be
sections of the bundle $\Lambda^2(M, \mathcal{F}(J^k_0( M)))$. These are differential forms
over $M$ that, when applied to sections of $\Gamma TM$ gives a  section of
$\mathcal{F}(J^k_0(M))$. In {\it section 2} such differential forms are introduced.
Particular attention is payed in the case of generalized metric, since it will be of
relevance to the construction of {\it geometries of maximal acceleration}. The fundamental notion of non-linear
 connection that we need to define natural objects in the relevant bundles is discussed.
 We describe the fundamental geometric and cohomological notions of generalized forms
 required in later developments. Then we consider the Cartan's calculus of differential forms, the rudiments (including a version of Stokes' theorem) and the corresponding de Rham cohomology (including a version of Poincar\'e's lemma and Thom's isomorphism theorem). Although the section is primarily of mathematical character, we have restricted ourself to developed the notions and results that are strictly necessary for the application to generalized higher order electrodynamics.

 In {\it section 3} the notion and fundamental properties of maximal
 acceleration geometry are briefly presented. Maximal covariant acceleration introduces a natural
 perturbation parameter (the inverse of the square of the norm of the maximal covariant acceleration),
 which is fundamental for our treatment and will eventually forbid run-away solutions in a natural way. We did not pursue a more fundamental explanation for maximal acceleration, restricting our considerations to a kinematic formulation of the new geometry.

  The generalized higher order fields are introduced
 in {\it section 4} from the cohomological perspective discussed in {\it section 2}.
 Also, we discuss the behavior of the singularities of the fields. The combination of
 the cohomology theory of {\it section 2} with the analytical structure
 that we assume for the electromagnetic fields provides an unambiguous way of generalizing
 the notion of electromagnetic field in {\it section 4}.

 In {\it section 5} the Lorentz-Dirac equation is
 discussed in the standard framework of classical electrodynamics, in the new setting
 of generalized higher order fields and in the contest of geometries of maximal acceleration. In particular, we use a simple way to obtain the Lorentz-Dirac
 equation, in a similar way as it was discussed by Rohrlich \cite{Rohrlich}, but for each of the possible different frameworks (standard electrodynamics in Minkowski spacetime, standard electromagnetic fields but with a maximal acceleration geometry and generalized electromagnetic fields in Minkowski spacetime).

 In {\it section 6}
 the second order differential equation \eqref{covariantequationofmotion}  describing the motion of a charged particle is obtained and its basic properties discussed. The argument is based strongly on the notions presented in {\it sections 3} and {\it section 4}
 and the method explained in {\it section 5}, applied in this case to generalized electromagnetic fields in a maximal acceleration geometry. We show that for generalized higher order fields, maximal acceleration
 and conservation of energy momentum are compatible. The resulting differential equation for a point charged particle is general covariant, free of run-away and pre-accelerated
 solutions of Dirac's type.

 In {\it section 7} a  set of differential equations for
 an effective electromagnetic theory on $M$ are obtained and the skeleton of a generalized gauge theory described. In particular, the explicit form of the generalized electromagnetic field and current is obtained, and prove that follow a generalized Maxwell's equations. Gauge invariance, the notion of vacuum and boundary conditions are also briefly discussed.

 In {\it section 8}, the results of the paper are briefly
 discussed and some relations with other approaches to electrodynamics highlighted. Finally, some directions of further investigation are indicated like the  development of the cohomology theory of generalized forms, its relation with calibrated geometry, the need of  a theory of curvature of generalized metrics, the extension to non-linear Yang-Mills theory and to gravity as well as the quantization of the theories with generalized higher order fields.

\subsection{Conventions}

$M$ is an $n$-dimensional, Hausdorff, second countable, smooth manifold (in applications to physics $n=4$,
although this restriction is not necessary for most of the mathematical results that we present). That $M$ hold such technical properties  $M$ assures the existence and uniqueness of natural geometric operators, in particular the existence of exterior algebra of differential forms and partitions of unity.  An arbitrary point in the manifold
$M$ is denoted by Latin letters, usually by $x$ or $p$. However, in the case of $x$, the same symbol will
be used to designate its coordinates in a local coordinate chart on $M$ or for curves on $M$. An arbitrary local coordinate
chart on $M$ is denoted by $(U,x)$, where $U$ is an open neighborhood. A curve on $M$ will be a smooth
map $x: I\to M$, being $I$ an interval of $R$.  The distinction between the three different
 meanings of $x$ will be clear from the context, although parameterized curves will be often denoted by
  $x(\sigma)$ or by $x(\tau)$, indicating the dependence in the given parameter. Greek letters are used for
   spacetime indices and run from $1$ to $n$. General smooth differential forms will also be
    denoted by Greek letters.
 However, for the electromagnetic fields, currents and electromagnetic potential $1$-form  we
     follow the traditional notation and denoted by  $F,\, G,\,A,\,J$ the Faraday,
     excitation, potential and current density forms. When considering fiber bundles
      over $M$, indices of the fiber over $x$ in the jet bundle are denoted by capital Latin letters.
       For instance for fibers over $x$ of the jet bundle $J^k_0(M)$ they run in $\{1,...,nk\}$.
        Einstein's summation convention is used, if anything else is not stated. In some occasions
         we will use multi-index notation, in order to simplify algebraic expressions. The manifold $M$ is
          equipped with a pseudo-Riemannian metric $\eta$ with signature $(-1,1,..,1)$. In more concrete applications, the
            spacetime manifold will be four dimensional and will have signature
             $(-1,1,1,1)$. In this case, when a $3+1$ decomposition is possible, three
             dimensional objects are indicated by arrowed vectors like $\vec{V}$, $\vec{A},\,\vec{B}$, etc...
\section{Differential geometry of jet bundles and
cohomology of generalized differential forms}

In this {\it section} we introduce the notions generalized tensors and generalized differential forms that we will use later.
Firstly, a short introduction to $k$-jets bundles $J^k_0(M)$ is necessary, since
the geometric objects that we will consider are constructed on such bundles.
Then a relevant class of connections defined on each $k$-jet bundle are explained.
Once we have chosen a connection, the notion of generalized tensors
and forms is given: they will be tensor and forms along maps on $M$ with values
on the algebra of functions $\mathcal{F}(J^k_0(M))$. The justification to introduce such
sophisticated objects is based on physical considerations and will be fully explained
in later sections. Notice that it allows a notion of electromagnetic field which is
capable to accommodate the changes produce by radiating point charged particles.

We will show that the fundamental notions of differential geometry can be transplanted to the
 framework of  generalized tensors and forms. For instance, a notion of generalized metric is introduced as generalized tensor,
as well as associated notions as distance, isometries, Hodge operator and isometry. The fundamental notions of the  causal theory for generalized metrics with Lorentzian signature are also discussed. The foundations of the
 Cartan's calculus for generalized forms is developed. With a notion of generalized metric and the
Cartan's calculus for generalized forms on hand, the fundamental results of the cohomology of generalized differential forms are proved.
In particular, integration along the fiber implies a $k$-jet bundle version of {\it Thom's isomorphism theorem}.
 We finish this {\it section}  introducing the fundamental notions of a
 theory of integration of generalized forms. We provide the fundamental properties of the integral, including a version of Stokes' theorem.

\subsection{Jet bundles and natural lifts}

Jet theory is a natural framework for the study of the geometry of ordinary and partial
differential equations (see for instance for a
introduction to jet theory the references \cite{BucataruConstantinescuDahl, deLeonRodrigues, KolarMichorSlovak, Michor1980, Saunders}).
Given a smooth curve $x:I\to { M}$, the set of derivatives
$(x(0),\,\frac{dx}{d\sigma},\,\frac{dx}{d \sigma}|_{0},\,...,\,\frac{d^k x}{d\sigma^k}|_{0})$ determines a point in the space of
jets of curves $J^k_0(p)$ in a neighborhood of the point $p=\,x(0)\in \,M$,
\begin{align*}
J^k_0(p):=\left\{(x(0),\,\frac{dx}{d \sigma}|_{0},\,...,\,\frac{d^k x}{d\sigma^k}|_{0}),
\,\,\forall\,\mathcal{C}^{k}\, x:I\to M,\, x(0)=p\in M,\,0\in I\right\}.
\end{align*}
 The jet bundle $J^k_0(M)$ over $M$ is then the disjoint union
\begin{align*}
J^k_0(M):=\bigsqcup_{x\in M}\,J^k_0(x).
\end{align*}
The projection map is
\begin{align*}
^k\pi  :J^k_0(M)\to M,\quad
 (x(0),\,\frac{dx^\mu}{d\sigma}\big|_0,\,\frac{d^2 x^\mu}{d\sigma^2}\big|_0,..., \frac{d^k x^\mu}{d\sigma^k}\big|_0)\mapsto x(0).
\end{align*}
Since it is smooth, the projection $^k\pi  :J^k_0(M)\to M$ determines the fiber bundle $(J^k_0(M),M,\,^k\pi)$.
The fiber over $p\in M$  $J^k_0(p):=\,^k\pi^{-1}(p)$ is a manifold of dimension $nk$.
Local coordinates on the fiber $J^k_0(p)$ are induced from the local coordinate system $(U,x)$
 over $M$ such that a natural system of local coordinates associated with the $k$-jet at $x(0)$ of the map $x:I\to M$ is
$(x^\mu(\sigma),\frac{dx^\mu}{d\sigma},\,\frac{d^2 x^\mu}{d\sigma^2},..., \frac{d^k x^\mu}{d\sigma^k})|_{\sigma_0=0}$.
A natural system of local coordinates of a point in $J^k_0(M)$ is denoted by  $(x^1,...,x^n,y^1,...,y^{nk})$.
The linear map $^k\pi:J^k_0(M)\to M$ is differentiable and the differential
 of the projection $^k\pi_u$ at $u\in\,^k\pi^{-1}(x)$ is the linear map
\begin{align*}
(\,^k\pi_*)_u:T_uJ^k_0(M)\to T_xM.
\end{align*}
We will denote by $\,^k\pi_*$ the projection $(\,^k\pi_*):TJ^k_0(M)\to TM$ such that at each $u$ it is the linear map $(\,^k\pi_*)_u$.

Given a curve $x:I\to M$, a $k$-lift is the curve $^kx:I\to J^k_0(M)$ such that the following diagram commutes,
\begin{align*}
\xymatrix{ &
{ J^k_0(M)} \ar[d]^{^k\pi}\\
{ I} \ar[ur]^{^kx}  \ar[r]^{x} & { M}.}
\end{align*}
In order to fix a $k$-lift one needs to specified each point $(x^1(\sigma),...,x^n(\sigma),y^1,...,y^{nk}(\sigma))$ of the $k$-lift:
each point on the lift $^kx:I\to J^k_0(M)$ has local coordinates given by
$(x^\mu(\sigma),\,\frac{dx^\mu(\sigma)}{d\sigma},...,\frac{d^k x^\mu(\sigma)}{d\sigma^k})$,
 where $\sigma$ is the parameter of the curves $x:I\to M$ and $^kx:I\to J^k_0(M)$.

There are also the notions of lift of  tangent vectors and smooth functions:
\begin{itemize}
\item Let $X\in T_xM$ be a tangent vector at $x\in M$. A lift $^kl_u(X)$ at $u\in\,^k\pi^{-1}(x)$
is a tangent vector at $u$ such that
\begin{align*}
\,^k\pi_* (\,^kl_u(X))=X.
\end{align*}

\item Let us denote by $\mathcal{F}(J^k_0( M))$ the algebra
 of real smooth functions over $J^k_0(M)$. Then there is defined the lift of a function $f\in\mathcal{F}(M)$
  to $\mathcal{F}(J^k_0(M))$ by
  \begin{align*}
 ^k\pi^*(f)(u)=f(\,x),\,\forall u\in\,^k\pi^{-1}(x).
  \end{align*}

\end{itemize}

The kernel of $(\,^k\pi_*)_x$ at $x\in\,M$ is the vector space
\begin{align*}
 ^k\mathcal{V}_x:=(\,^k\pi_*)^{-1}(0_x).
 \end{align*}
 where $0_x$ is the zero vector in $T_xM$.
 Then vertical bundle over $M$ is the manifold
 \begin{align*}
 ^k\mathcal{V}:=\,\bigsqcup_{x\in\,M}\, ^k\mathcal{V}_x
 \end{align*}
with the induced projection $^k\tilde{\pi}_V:\,^k\mathcal{V}\to M. $
The vertical bundle over $J^k_0(M)$ is determined by the surjection
\begin{align*}
^k\pi_V:\,^k\mathcal{V}\to J^k_0(M),\quad (\,^k\pi_*)^{-1}(0_x)\ni \,\xi^v \,\mapsto u\in\,(\pi^k)^{-1}(x).
\end{align*}
$^k\mathcal{V}$ is a real vector bundle over $J^k_0(M)$, since it is the kernel of $\,^k\pi_*$.
The composition of $^k\pi_V$ with $^k\pi$ determines also a real vector bundle  over $M$,
\begin{align*}
^k\pi\circ\,^k\pi_V:\,^k\mathcal{V}\to M.
\end{align*}
One can introduce the notation $^k\tilde{\pi}_V $
and the vertical bundle over $J^k_0(M)$ is determined by the surjection
\begin{align*}
^k\pi_V:\,^k\mathcal{V}\to J^k_0(M),\quad (\,^k\pi_*)^{-1}(0_x)\ni \,\xi^v \,\mapsto u\in\,(\pi^k)^{-1}(x).
\end{align*}
$^k\mathcal{V}$ is a real vector bundle over $J^k_0(M)$, since it is the kernel of $\,^k\pi_*$.
The composition of $^k\pi_V$ with $^k\pi$ determines also a real vector bundle  over $M$,
\begin{align*}
^k\pi\circ\,^k\pi_V:\,^k\mathcal{V}\to M.
\end{align*}

\subsection{Connections on $TJ^k_0(M)$ and a covariant derivative on $TM$}

There are several general notions of connections in differential geometry.
Let us recall the notion of general connection that will be particularly useful for us (see for instance \cite{KolarMichorSlovak}).
On an arbitrary fiber bundle $\pi_{\mathcal{E}}:\mathcal{E}\to M$,  the vertical
bundle is $\mathcal{VE}:=\,ker( (\pi_{\mathcal{E}})_*)$.
\begin{definicion}
A connection on the bundle $\pi_{\mathcal{E}}:\mathcal{E}\to M$ is a vector valued $1$-form
$\Phi_{\mathcal{E}}:T\mathcal{E}\to \mathcal{VE}$
 such that
 \begin{enumerate}
\item $\Phi^2_{\mathcal{E}}=\Phi_{\mathcal{E}}$,

\item $Im \,\Phi=\mathcal{VE}$.
\end{enumerate}
\label{generalconnection}
\end{definicion}
Therefore, a connection is a projection operator on $T\mathcal{E}$.
 For the map
 \begin{align*}
 h_{\mathcal{E}}=\,(Id\,-\Phi_{\mathcal{E}}):T\mathcal{E}\to T\mathcal{E}
 \end{align*}
  one has that
\begin{align*}
h^2_{\mathcal{E}}=\,h_{\mathcal{E}}.
\end{align*}
$ker(\Phi)$ is a sub-bundle $\mathcal{HE}$ of $\mathcal{E}$ called the {\it  horizontal bundle}.
 Also, it is direct that the following decomposition holds,
\begin{align*}
T\mathcal{E}=\,\Phi_{\mathcal{E}}(T\mathcal{E})
\oplus h_{\mathcal{E}}(T\mathcal{E}):=\,\mathcal{VE}\oplus\mathcal{HE}.
\end{align*}
Thus
for each $u\in \mathcal{E}$ there is a unique decomposition
$T_u\mathcal{E}=\,\mathcal{H}_u\mathcal{E}\oplus\mathcal{V}_u\mathcal{E}.$
Given the connection $\mathcal{HE}$, the horizontal lift of a vector field $X\in\,\Gamma\, TM$ is a horizontal
 field $^hX\in\Gamma\, \mathcal{HE}$ such that $^k\pi_* (\,^hX)=X$.
\begin{comentario}
It is equivalent to give a connection $\Phi_\mathcal{E}$ or give a distribution $\mathcal{HE}$ satisfying \ref{generalconnection}. The distribution $\mathcal{HE}$ is not necessarily integrable and the integrability obstruction is given by the curvature of the connection.
\end{comentario}
\begin{definicion}
A connection on $^k\pi:J^k_0(M)\to M$ is a global splitting
\begin{align*}
TJ^k_0(M)=\,\mathcal{V}J^k_0(M)\oplus{\mathcal{H}}J^k_0(M),
\end{align*}
and  with $\mathcal{V}_uJ^k_0(M)=\,^k\mathcal{V}_u=\,ker(\,^k\pi_*)_u$ for each $u\in\,J^k_0(M)$.
\label{definicionnonlinearconnection}
\end{definicion}
We will denote by $^k\hat{\mathcal{H}}=\mathcal{H}J^k_0(M)$ to the horizontal distributions.

Let us consider
 the Levi-Civita connection $D$ of the semi-Riemannian metric $\eta$ on $M$.
Our objective in this sub-section is to introduce a {\it natural connection}  $\mathcal{H}_D\hookrightarrow TTM$ on $TM$. Related with this connection, it will be shown later the existence of connections on the jet bundles $J^k_0(M)$.

In order to introduce the connection $\mathcal{H}_D$ on $TM$, let us consider the commutative diagram
\begin{align}
\xymatrix{TTM\ar[d]_{\pi_1} \ar[r]^{\pi_2} &
 T M \ar[d]^{\pi}\\
TM \ar[r]^{\pi} &  M}
\label{diagramacommutativo0}
\end{align}
and let $\mathcal{V}=\,ker (\pi_1)\,\subset TTM$ be the vertical sub-bundle.
The vertical lift of a vector field is canonically defined \cite{Crampin} by the relations
\begin{align*}
^vl:TM\to \mathcal{V},\quad \hat{X}\mapsto\, \,(^v\hat{X})_{u=(x,y)}:=\frac{d(x,y+t\hat{X})}{dt}\Big|_{t=0}.
\end{align*}
The geodesic spray $\mathcal{S}$ is the section $\mathcal{S}\in\,\Gamma TTM$ such that $\pi_* \circ \mathcal{S}=\,\pi_1\circ S$ and the geodesic curves of $\mathcal{S}$ is the geodesic flow of $\eta$.
Then the  connection $\mathcal{H}_D$ is determined by the distribution of $TTM$ determined by the conditions \cite{Crampin, Crampin00},
\begin{enumerate}
\item It is torsion free,
\begin{align}
[\,^hX,\,^vY]-\,[\,^hY,\,^vX]-\,^v[X,Y]=0,
\label{torsionfreecondition}
\end{align}
where $^hX$ is the horizontal lift of the connection $D$.
\item The geodesic spray
\begin{align}
  \mathcal{S}=\,y^i\frac{\partial}{\partial x^i}+\,\gamma^i\,_{jk}(x)\,y^jy^k
\end{align}
 of the Riemannian connection ${D}$ determines the horizontal lift by the formula
\begin{align}
^hX_u:=\frac{1}{2}([\,^vX,\,\mathcal{S}]_u\,+(\,^cX)_u),\quad X\in\,T_x M,\,u\in\,\pi^{-1}(x).
\label{Connectionassociatedtoaspray}
\end{align}
\end{enumerate}
One has the following
\begin{proposicion}
The set $\{(\,^hX)_u,\,X\in\,T_xM,\,u\in\pi^{-1}(x),\,x\in\,M\}$ determines a vector space of sections
that defines a connection $\mathcal{H}_D$ on $TM$.
\end{proposicion}
 Since $dim(\mathcal{V})=n$, one has that $dim(\mathcal{H}_D)=n$.
 The connection $\mathcal{H}_D$ has associated the horizontal projection operator $h_D:TTM\to\mathcal{H}_D$ and
$\mathcal{H}_D=Im(h_D)$.
The horizontal lift of an arbitrary vector field $X\in\,\Gamma TM$ is
denoted by $^hX\in\,\Gamma TTM$. Given $\hat{X}\in \,\Gamma\,TTM$, its horizontal
component is denoted by $^h\hat{X}$.
\begin{ejemplo}
Let $(M,\eta)$ be a Lorentzian manifold.
Given an holonomic local frame $\{\frac{\partial}{\partial x^{\mu}},\,\mu=1,...,n\}$
the connection coefficients of the Levi-Civita connection $D$ are given by the expression
\begin{align*}
\big(\,D_{\partial_{\nu}}\,\partial_{\rho}\big)^{\mu}:= {\gamma}^{\mu }\,_{\nu \rho}=\frac{1}{2}\eta^{\mu s}(\frac{\partial \eta_{s\nu}}{\partial
x^{\rho}}-\frac{\partial \eta_{\rho \nu}}{\partial
x^{s}}+\frac{\partial \eta_{s\rho}}{\partial x^{\nu}}),\quad \mu ,\nu
,\rho , s=1,...,n.
\end{align*}
\subsection*{A connection on $TTM$ associated with the connection $D$ on $TM$}

Given the Levi-Civita connection $D$ of the metric $\eta$,
one can introduce
the homomorphism
\begin{align*}
N:& TTM \to TTM,\quad  \frac{\partial}{\partial x^{\mu}}\mapsto N^{\nu}\,_{\mu}\,\frac{\partial}{\partial y^{\nu}},\quad \frac{\partial}{\partial y^{\mu}}\mapsto \frac{\partial}{\partial y^{\mu}},\quad \mu,\nu=1,...,n.
\end{align*}
such that the projector $h_D$ is given by
\begin{align}
h_D:& TTM \to TTM,\quad \hat{X}\mapsto\, (Id-\,N)(\hat{X}).
\label{homomorphismhD}
\end{align}
Then $h_D$ is determined by the components of $N$,
\begin{align}
{N^{\mu}\,_{\nu}}(u)=\gamma^{\mu}\,_{\nu \rho}(x)\,y^{\rho},\quad \mu , \nu ,\rho =1,...,n.
\label{nonlinearconnection}
\end{align}
\begin{proposicion}
The homomorphism \eqref{homomorphismhD} determines  a connection in the sense of the Definition \ref{generalconnection} on the bundle $TM$.
\end{proposicion}

An adapted basis for $\mathcal{H}_D$ is given by the globally defined distribution
\begin{align}
\left\{ \frac{{\delta}}{{\delta} x^{1}}\big|_u ,...,\frac{{\delta}}{{\delta} x^{n}} \big|_u, 
\frac{\partial}{\partial y^{1}} \big|_u,...,\frac{\partial}{\partial
y^{n}} \big|_u\right\},\quad
 \frac{{\delta}}{{\delta} x^{\nu}}\big|_u =\frac{\partial}{\partial
x^{\nu}}\big|_u -N^{\mu}\,_{\nu}(u)\frac{\partial}{\partial y^{\mu}}\big|_u
\end{align}
for $\mu,\nu=1,...,n.$ Then the set of local sections
\begin{align}
\{ \frac{{\delta}}{{\delta} x^{1}}|_u 
,...,\frac{{\delta}}{{\delta} x^{n}}|_u,\, u\in {TU} \}
\label{generatorofthehorizontaldistributiononTM}
\end{align}
generates the local
horizontal distribution $\mathcal{H}_D|_U $, while the set of local sections
\begin{align}
\{ \frac{\partial}{\partial y^{1}}|_u ,..., 
\frac{\partial}{\partial y^{n}}|_u, u\in { TU} \}
\label{generatorlocalverticaldistributiononTM}
\end{align}
generates the local vertical distribution 
$\mathcal{V}|_U$.
\end{ejemplo}

\subsection*{The covariant derivative of D}

Given a connection on ${TM}$, there is associated a covariant derivative on $M$, defined in local coordinates by the expression
\begin{align}
(D _{\dot{x}}\,Z)^{\sigma}:=\dot{Z}^{\sigma} + N^{\sigma}\, _{\bf \rho}((x,\dot{x})) Z^{\rho},\quad u=(x,\dot{x})\in\, TM,\,Z\in\,\Gamma TM.
\label{nonlinearcovariantderivative}
\end{align}
The derivatives are performed respect to $\sigma$ (in this case, an affine parameter).
The geodesic equation for the linear connection $D$ is given as the curve on $M$ solution of the differential equation
\begin{align}
D_{\dot{x}}\,\dot{x}=0,
\label{geodesicequationforD}
\end{align}
Equation (\ref{geodesicequationforD}) corresponds to the geodesic equation of the pseudo-Riemannian metric $\eta$.

\subsection{Induced connections on $TJ^k_0(M)$ from connections on $TM$}

Given a connection $\mathcal{H}$ on $TM$, we can show the existence of a connection $^k\hat{\mathcal{H}}$
on $J^k_0(M)$ related with $\mathcal{H}$. The construction is similar to the
{\it induced connection} of the references \cite{KN, KolarMichorSlovak}. The connections $^k\hat{\mathcal{H}}$
will be used in the construction of geometric objects and some relevant natural operators on $J^k_0(M)$.

Let us consider the commutative diagram,
\begin{align}
\xymatrix{J^k_0(M)\ar[d]_{^k\pi} \ar[r]^{^kpr} &
 TM \ar[d]^{\pi}\\
M \ar[r]^{id} &  M.}
\label{diagramacommutativo1}
\end{align}
The commutative diagram (\ref{diagramacommutativo1}) induces another commutative diagram
for each $^kx\in \,J^k_0(M)$; if $(x,y)\in\,TM$ and $x\in\, M$, with $pr(\,^kx)=(x,y)$ and $\pi(x,y)=x$, then one has that
\begin{align}
\xymatrix{T_{\,^kx}J^k_0(M)\ar[d]_{\,^k\pi_*} \ar[r]^{^kpr_*} &
 T_{(x,y)}TM \ar[d]^{\pi_*}\\
T_xM \ar[r]^{id} &  T_xM}
\label{diagramacommutativo2}
\end{align}
is also commutative.
\begin{lema}
 Let us consider a distribution $^k\hat{\mathcal{H}}\subset\,TJ^k_0(M)$ such that
 $pr_*(\,^k\hat{\mathcal{H}})\subset\, \mathcal{H}$ and that the diagram (\ref{diagramacommutativo2})
  commutes. Then it holds that $^k\hat{\mathcal{H}}\cap\,^k\mathcal{V}=\,0$.
\label{lemasobrehatH1}
\end{lema}
{\bf Proof}. If the vector field $\hat{X}_0\in\,^k\hat{\mathcal{H}}$ and is vertical,
then $(\,^k\pi_*)(\hat{X}_0)=0.$ The commutative diagram (\ref{diagramacommutativo2}) implies
that $^kpr_*(\hat{X}_0)\in\,\mathcal{H}_D$ and $^k\pi_* (\,^kpr_*(\hat{X}_0))=0$. Thus the relation $^kpr_*(\hat{X}_0)=0$ holds, since the only
 horizontal vector $U$ such that $^kpr_*(U)_=0$ is the null vector.\hfill$\Box$
\begin{lema}
Let us consider $^k\hat{\mathcal{H}}$  such that $^kpr_*(\,^k\hat{\mathcal{H}})\subset\, \mathcal{H}$ and that
 the diagrams (\ref{diagramacommutativo2}) commute. Then  $dim(\,^k\hat{\mathcal{H}})\leq \,dim(M)=n$.
\end{lema}
{\bf Proof}. Assume that it has dimension $n+1$. Then there is a basis $\{\hat{X}_1,...,\hat{X}_n,\hat{Z}\}$
of $^k\hat{\mathcal{H}}$ such that $\hat{X}_i=\,X_i$ for each $i=1,...,n$. Therefore,
\begin{align*}
 ^kpr_*(\hat{Z})=a^iX_i=\,a^i\,^kpr_*(\hat{X}_i)=\,^kpr_*(a^i\,\hat{X}_i),
 \end{align*}
 from what follows that $^kpr_*(\hat{Z}-\,a^i\,\hat{X}_i)=0$.
Since diagrams (\ref{diagramacommutativo2}) are commutative, it follows that
$\hat{Z}-\,a^i\,\hat{X}_i$ is vertical. Therefore, it must be zero (by {\it Lemma} \ref{lemasobrehatH1}).\hfill$\Box$
\begin{proposicion}
Let us consider $^k\hat{\mathcal{H}}_m$ of maximal dimension such that $^kpr_*(\,^k\hat{\mathcal{H}})\subset\, \mathcal{H}$. Then it has dimension $dim(\,^k\hat{\mathcal{H}}_m)n$.
\end{proposicion}
{\bf Proof}. We have that $dim(TJ^k_0(M))=n(k+1)$, $dim(ker(\,^kpr_*))=\,nk$. On the other
 hand, $^kpr_*$ is surjective. Therefore, $dim(Im (pr_*))=n$. Assume that there is a $X_0\in \,\mathcal{H}_m$
 such that it is not the image of some $\hat{X}_0\in \,^k\hat{\mathcal{H}}$. This is
  in contradiction with the following facts (just counting dimensions):
\begin{enumerate}
\item  $\,^k\hat{\mathcal{H}}_m$ is of maximal dimension,
\item  $pr_*(\,^k\hat{\mathcal{H}}_m)\subset\,\mathcal{H}$,
\item The projection $^kpr_*$ is surjective.
\end{enumerate}
Therefore, it follows that $dim(\,^k\hat{\mathcal{H}}_m)=n$.
\hfill$\Box$
\begin{teorema}
A distribution $^k\hat{\mathcal{H}}_m$ of maximal dimension such that $pr_*(\,^k\hat{\mathcal{H}})\subset\, \mathcal{H}$
and that the diagrams (\ref{diagramacommutativo2}) commutes is a connection on $J^k_0(M)$.
\label{teoremasobreconexiones}
\end{teorema}
{\bf Proof}. Let $^k\hat{\mathcal{H}}_m$ a maximal distribution of dimension $n$. As we showed
in {\it Lemma} \ref{lemasobrehatH1} $^k\hat{\mathcal{H}}_m\,\cap\,^k\mathcal{V}=0$. This
shows that $TJ^k_0(M)=\,^k\hat{\mathcal{H}}_m\,\oplus \,^k\mathcal{V}$.\hfill$\Box$
\begin{definicion}
Given $\mathcal{H}$ a connection on $TM$,  a connection $^k\hat{\mathcal{H}}$  on $J^k_0(M)$ as
in {\it Theorem} \ref{teoremasobreconexiones} is an {\it induced connection} on $J^k_0(M)$.
\end{definicion}
We considered an induced connection $^k\hat{\mathcal{H}}$ on $J^k_0(M)$ induced
by the connection $\mathcal{H}_D$, determined from the Levi-Civita connection $D$. Let us note that
{\it Theorem} \ref{teoremasobreconexiones} does not fix the connection $^k\hat{\mathcal{H}}$,
but only shows the existence of such connections. Indeed, because of the difference between
the dimensionality between $T_{\,^kx}J^k_0(M)$ and $T_{(x,y)}TM$, it is not possible to fix
the connection in such a way and a bunch of connections exists that full-fill {\it Theorem} \ref{teoremasobreconexiones}.
 However, one can characterize such connections as follows,
\begin{proposicion}
Let us consider two connections $^k\hat{\mathcal{H}}_1,\,^k\hat{\mathcal{H}}_2$ constructed
as in {\it Theorem} \ref{teoremasobreconexiones}. If the
corresponding horizontal lifts  $^{h_1}X\in\,^k\hat{\mathcal{H}}_1|_{\,^kx}$ and
 $^{h_2}X\in\,^k\hat{\mathcal{H}}_2|_{\,^kx}$  of $X\in T_xM$ are such that $pr_*(\,^{h_1}X\,-\,^{h_2}X)=0$, then the vector field $^{h_1}X\,-\,^{h_2}X\in \,TJ^k_0(M)$ is vertical.
\label{connexionesconnection independentes}
\end{proposicion}
{\bf Proof}. It is a direct consequence of the fact that the diagrams (\ref{diagramacommutativo2}) commute.\hfill$\Box$

The following notion captures when a property is independent of the connection:
\begin{definicion}{\it Equivariant relations.}
\begin{itemize}

\item Two connections $^k\hat{\mathcal{H}}_1$ and $^k\hat{\mathcal{H}}_2$ are said to be equivalent  iff $\,^{h_1}X\,-\,^{h_2}X\,\in\,^k\mathcal{V},\,\forall X\in \,TM$.

\item A property is connection independent if only depends on the induced connection as in
in {\it Theorem} \ref{teoremasobreconexiones} and if it is true for any set of connection independent connections.
\end{itemize}
\label{definiciondeconexionesquivariantes}
\end{definicion}
If we denote by $^k{\bf H}$ the set of connections as in {\it Theorem} \ref{teoremasobreconexiones}, the following result holds:
\begin{proposicion}
The relation of being connection independent for induced connections \ref{definiciondeconexionesquivariantes} is an equivalence relation in $^k{\bf H}$.
\end{proposicion}

Each induced connection $^k\hat{\mathcal{H}}$ has associated a projection map $\hat{h}_k$,
\begin{align*}
\hat{h}_k:TJ^k_0(M)\to\,^k\hat{\mathcal{H}}
\end{align*}
with the property that $(\hat{h}_k)^2=\,\hat{h}_k$.
\begin{comentario}
We have the following remarks:
\begin{itemize}
\item One way to define a connection on $J^k_0(M)$ is by using a semi-spray of order $k$
plus additional conditions on the variation equations \cite{BucataruConstantinescuDahl}. However, for the problem that we
are interested on this paper (the dynamics of point charged particles with the electromagnetic field),
we do not have  any natural semi-spray of order $k$ which is free of pathological solutions.
We have two natural semi-sprays defined in the problem. The first one corresponds
 to the Lorentz-Dirac equation. This is a spray of third
 order. However, the Lorentz-Dirac spray is not physically acceptable: it is well known that it contains
 un-physical solutions. The natural spray corresponds to the geodesic equation of
 the pseudo-Riemannian metric $\eta$. This is a spray of second order $k=2$.
 It determines geodesics describing the dynamics of free point charged particles without interacting with an electromagnetic field, thus not suitable to our problem (description of the motion of a point charged particle).

 \item Kaluza-Klein  dimensional reduction theories does not provides a natural spray to define a convenient connection. The reduction
  reproduces the Lorentz force equation, that brings to the game several connections
   (see for instance \cite{Miron2004, Miron2006, Ricardo09}). However, it does not take into account radiation reaction effects.

 \item The equation governing the dynamics of a point charged particle is not necessarily a (semi)-spray. It will turn out that is an implicit differential equation of second order and does not have the form of a spray.

\end{itemize}
\end{comentario}

Given a connection $^k\hat{H}$ on $J^k_0(M)$, the local connection coefficients of the endomorphism $N$ in the natural basis
\begin{align}
\left\{\frac{\partial}{\partial x^1},...,\frac{\partial}{\partial x^n},\frac{\partial}{\partial y^1},...,\frac{\partial}{\partial y^{nk}}\right\}.
\end{align}
are $\{N^A\,_{\mu},\,\mu=1,...,\,A=1,...,nk\,\}$.
 The vertical distribution is generated locally by the local frame,
\begin{align*}
\left\{ \frac{\partial}{\partial y^{1}}|_x ,...,\frac{\partial}{\partial y^{nk}}|_x, x\in { U} \subset \,M\right\}.
\end{align*}
Let us adopt a connection as in {\it Theorem} \ref{teoremasobreconexiones}.
Then there is a splitting of any local frame of $TJ^k_0(M)$ such that the horizontal distribution is locally generated by
\begin{align}
\left\{ \frac{{\delta}}{{\delta} x^{1}}|_x ,...,\frac{{\delta}}{{\delta} x^{n}} |_x, 
\frac{\partial}{\partial y^{1}} |_x,...,\frac{\partial}{\partial
y^{nk}} |_x\right\},\quad
 \frac{{\delta}}{{\delta} x^{\nu}}|_x =\frac{\partial}{\partial
x^{\nu}}|_x -N^{A}\,_{\nu}\frac{\partial}{\partial y^{A}}|_x
\end{align}
with $\nu=1,...,n,\, A=1,...,nk.$ Therefore, the functions $kn^2$ functions $N^A\,_{\mu}(x,y)$
 determines the connection in the given coordinate system.
From {\it Theorem} \ref{teoremasobreconexiones}, $n^2$ of those functions are fixed.
 Therefore, there are still $(k-1)n^2$ to be fixed by imposing additional conditions.

A way to impose conditions on the connection is the following,
\begin{definicion}
Let us consider a connection $^k\mathcal{H}$ of $J^k_0(M)$.
An horizontal $p$-form  is such that $\omega(...,V,...)=0$ for any vertical vector field $V\in\,\hat{\mathcal{V}}$.
\label{definicionhorizontallform}
\end{definicion}
 A local frame of vertical forms is determined by exterior product of vertical $1$-forms
\begin{align}
\delta y^{A}=\,dy^{A}\,+N^{A}\,_{\nu}\,dx^{\nu},\quad \nu=1,...,n,\quad A=1,...,nk,
\label{verticalforms}
\end{align}
where the dual co-frame of $1$-forms $\{dx^1,...,dx^n,dy^1,...,dy^{nk}\}$ defined by the relations
\begin{align*}
dx^i(\frac{\partial }{\partial x^j})=\delta^i_j,\quad dx^i(\frac{\partial }{\partial x^A})=0,
\end{align*}
\begin{align*}
 dx^A(\frac{\partial }{\partial x^j})=0,\quad dy^A(\frac{\partial}{\partial y^B})=\delta^A_B,\quad\,i,j=1,...,n,\,\, A,B=1,...,nk.
\end{align*}
\begin{proposicion}
Every horizontal  $p$-form is connection independent.
\label{Propositionnotionhorizontalisconnection independent}
\end{proposicion}
{\bf Proof}. Given two equivalent connections  $^k\mathcal{H}_1$ and  $ ^k\mathcal{H}_2$, if a vector $Z$ is horizontal respect to the first, it must be re-written in the form $^hZ_1=\,^h{Z}_2+ V(\,^hZ_1)$, with $\,^h{Z}\in\,^k\mathcal{H}_2$ and $V(\,^hZ_1)$ vertical. Thus
\begin{align*}
\omega(...,\,^hZ_1,...)& =\,\omega(...,\,^hZ_2+V(Z),...)\\
& = \omega(...,\,^hZ_2,...)+\omega(...,\,V(Z),...)\\
& = \omega(...,\,^hZ_2,...).
\end{align*}\hfill$\Box$

The exterior product of horizontal forms is also horizontal and connection independent, forming the algebra of horizontal differential forms.
\begin{ejemplo}
A tensor field of type $(p,q)$ $T$ is a section of the bundle $T^{(p,q)}J^k_0(M)$. Given a connection $^k\hat{\mathcal{H}}$ on $J^k_0(M)$,
 one can define the notions of
horizontal and vertical tensors of any order.

A tensor field $T$ of type $(2,0)$ is a section of the bundle $T^{(2,0)}J^k_0(M)$.
 Locally, the sections of  $T^{(2,0)}J^k_0(M)$ are spanned (with coefficients on $\mathcal{F}(J^k_0(M))$
  by the tensor product of elements of the frame $\{e_1(\,^kx), ..., e_{n(k+1)}(\,^kx)\}$
\begin{align*}
T^{(2,0)}J^k_0(M)=\,span \left\{e_i(\,^kx)\otimes\,e_j(\,^kx),\,i,j=1,...,n(k+1)\right\}.
\end{align*}
If there is defined a connection $^k\hat{\mathcal{H}}$ as in {\it Theorem} \ref{teoremasobreconexiones},
 the tensor bundle of $(0,2)$ horizontal tensors is locally spanned as
\begin{align*}
T^{(2,0)}_hJ^k_0(M)=\,span \left\{\frac{\delta}{\delta x^{\mu}}\Big|_{\,^kx}\otimes \,\frac{\delta}{\delta x^{\nu}}\Big|_{\,^kx},\,\mu,\nu=1,...,n\right\}.
\end{align*}
Similarly, the tensor bundle of $(2,0)$ vertical tensors is locally spanned as
\begin{align*}
T^{(2,0)}_hJ^k_0(M)=\,span \left\{\frac{\partial}{\partial y^{A}}\Big|_{\,^kx}\otimes
 \,\frac{\partial}{\partial y^{B}}\Big|_{\,^kx},\,\, A,B=1,...,nk\right\}.
\end{align*}
The bundle over $J^k_0(M)$ of $hv$ tensors of type $(1,1)$ is locally spanned
\begin{align*}
T^{(1,1)}_hJ^k_0(M)=\,span \left\{\frac{\delta}{\delta x^{\mu}}\Big|_{\,^kx}\otimes
 \,\frac{\partial}{\partial y^{A}}\Big|_{\,^kx},\,\, \mu=1,...,n,\,A=1,...,nk\right\}.
\end{align*}
One can consider an horizontal $(1,1)$ tensor, that generically will have the following expression in local coordinates
\begin{align*}
T(\,^kx)=\,T^i\,_j(\,^kx) \frac{\delta}{\delta x^i}\Big|_{\,^kx}\otimes \delta x^j\Big|_{\,^kx},\,T^i\,_j(\,^kx)\in\,\mathcal{F}(J^k_0(M)).
 \end{align*}
Horizontal $p$-forms can be spanned locally in a similar way.
For instance,
a $2$-form horizontal can be expressed in local coordinates as
\begin{align*}
\omega(\,^kx)=\,\omega_{ij}(\,^kx)dx^i|_{\,^kx}\wedge dx^j|_{\,^kx}.
\end{align*}
Note that horizontal forms does not depend on the specific connection $^k\hat{\mathcal{H}}$ that we can choose.

The tensor product of horizontal tensors is an horizontal tensor. In a similar way,
the exterior product of horizontal forms is an horizontal form. We will show below  that there is a notion of {\it horizontal exterior derivative}, for which definition we need a connection on $J^k_0(M)$. We will also show that the notion of horizontal exterior derivative is connection independent.
\end{ejemplo}

\subsection{Generalized forms and tensors fields}

Let $J^k_0(M)$ be the $k$-jet bundle over $M$ and $\eta$ a pseudo-Riemannian metric on $M$.
Let us consider a connection $^k\hat{\mathcal{H}}$ on $J^k_0(M)$ as in {\it Theorem} \ref{teoremasobreconexiones}.
\begin{definicion}
Let $pr:T^{(p,q)}J^k_0(M)\to J^k_0(M)$ be a tensor bundle over  $J^k_0(M)$.
A $k$-tensor along the curve $^kx:I\to M$ (in short, a tensor along the curve $^kx:I\to M$) is a map
$\hat{S}:I\to T^{(p,q)}J^k_0(M)$ such that for the lift $^kx:I\to M$ to $J^k_0(M)$, the following diagram commutes:
\begin{align*}
\xymatrix{ &
{T^{(p,q)} J^k_0(M)} \ar[d]^{pr}\\
{ I} \ar[ur]^{\hat{S}}  \ar[r]^{^kx} & { J^k_0(M)}.}
\end{align*}
A tensor  $\hat{S}$ along the curve $^kx:I\to M$ is horizontal if when acting on any
 arbitrary vertical vector or vertical $1$-form the result is zero.
 A tensor $\hat{S}$ is vertical if when acting on any horizontal vector the result is zero.
\label{definicionofhorizontaltensor}
\end{definicion}
One defines in a similar way differential forms along the curve $^kx:I\to M$ as a map $\hat{\omega}:M\to \Lambda^p J^k_0(M)$
 such that the following diagram commutes,
\begin{align*}
\xymatrix{ &
{\Lambda^p J^k_0(M)} \ar[d]^{pr}\\
{ I} \ar[ur]^{\hat{\omega}}  \ar[r]^{^kx} & { J^k_0(M)}.}
\end{align*}
Let us define the following projections:
\begin{align*}
\hat{h}_k:\,T^{(p,q)}(J^k_0(M))&\to \,T^{(p,q)}_h(J^k_0(M)),\quad T\mapsto \hat{h}_k(T)
\end{align*}
such that on horizontal vectors and forms $T=\,\hat{h}_k(T)$, but if $\hat{h}_k(T)$ acting over any vertical vector
 or differential $1$-form is zero.  Then the horizontal tensor bundle is
\begin{align*}
\rho^{(p,q)}:T^{(p,q)}_h(J^k_0(M))\to J^k_0(M).
\end{align*}
In a similar way, one can define the horizontal forms bundle $\Lambda^p_h(J^k_0(M))$ and the corresponding projection
$pr:\Lambda^k(J^k_0(M))\to J^k_0(M)$.
We have the following direct generalization of {\it Proposition} \ref{Propositionnotionhorizontalisconnection independent},
\begin{proposicion}
Vertical tensors $T\in T^{(p,0)}J^k_0(M)$  are connection independent.
\label{equivarianceofhorizontaltensors}
\end{proposicion}
{\bf Proof}.  For vertical tensors $T\in T^{(0,q)}J^k_0(M)$ acting on horizontal $1$-forms $\theta=\theta_\mu dx^\mu$ we have that
\begin{align*}
T(...,\theta,...)=T(...,\theta_\mu\, dx^\mu,...)=0
\end{align*}
and is clear that the result does not depend on the connection.\hfill$\Box$

The tensorial product of connection independent tensors is connection independent, thus forming the sub-algebra of connection independent tensors of the tensor algebra $(\sum_p T^{(p,0)}_h J^k_0(M),\otimes )$.
 However, similar calculations show that the notion of vertical form and horizontal contravariant tensor is not connection independent.
\subsection*{Weak definition of generalized tensors and forms}
\begin{definicion}
Let $(M,\eta)$ be a pseudo-Riemannian manifold.
A $k$-generalized $(p,q)$-tensor $\hat{S}$ is a map such that:
\begin{enumerate}
\item To each causal curve
 $x:I\to M$ associates a horizontal tensor $\hat{S}:I\to T^{(p,q)}_hJ^k_0(x(I))$
  along $^kx:I\to M$ as in definition \ref{definicionofhorizontaltensor}.

 \item If two lifted curves $^kx_1:I\to J^k_0(M)$ and
$^kx_2:I\to J^k_0(M)$ intersect at the point $^kx(s_0)\in\,^kx_1\subset\,^kx_2$,
 then the value at the intersection of the generalized tensor coincides,
 \begin{align}
 \hat{S}(\,^kx_1(s_0))=\,\hat{S}(\,^kx_2(s_0)).
 \label{microlocalitycondition}
 \end{align}
 \end{enumerate}
A $k$-generalized $p$-form is a map such that to each timelike curve
 $x:I\to M$ associates a horizontal differential form $\hat{\omega}:I\to \Lambda^p_h J^k_0(x(I))$
  along $x:I\to M$.
\label{generalizedfieldalongacurve}
\end{definicion}
Therefore, given a curve $x:I\to M$ and the lift of the curve (with a prescribed initial point over the initial point $x(\sigma_0)$)
$^kx:I\to J^k_0(M)$ a $k$-generalized $(p,q)$-tensor associates a unique section of $T^{(p,q)}_hJ^k_0(M)$
along the map $^kx:I\to J^k_0(M)$.
\begin{comentario}
We give some remarks on  {\it definition} \ref{generalizedfieldalongacurve}:
 \begin{itemize}

\item The second condition in {\it definition} \ref{generalizedfieldalongacurve} is related with {\it locality character} of physical fields and generalization of the
 idea o local field to higher order jet bundles. The physical interpretation is that the fields are determined by the trajectories of the probe particles.

\item We restrict to consider horizontal fields only. These fields have un-ambiguous interpretations in terms macroscopic
 flux across spatial surfaces, which is how energy-momentum variations are computed.

\end{itemize}
We will restrict our attention to causal curves $x:I\to M$ in {\it definition} \ref{generalizedfieldalongacurve} because its relevance for physical applications.
\end{comentario}

\subsection*{Strong definition of generalized tensors and forms}

Definition \ref{generalizedfieldalongacurve} provides a
notion of tensors along lifted curves to $J^k_0(M)$. However, this notion is not
enough to define some natural operations on generalized tensors and generalized differential forms (in particular to consider an exterior derivative operator and a generalized Lie derivative).
 In order to be able to define such operations, one needs to define them on open neighborhoods
 of $J^k_0(M)$ and on {\it variations}. One can achieve a more convenient definition of generalized field by considering
  fields over lifts to $J^k_0(M)$ of {\it tubes} in $M$. A way to construct convenient tubes is the following.
 Let $\Sigma^{n-1}\hookrightarrow M$ be a sub-manifold of $M$ and
 \begin{align*}
 \Delta:I\times\,\Sigma^{n-1} \to M,
 \end{align*}
 be a congruence of curves such that the lift
 \begin{align}
 ^k\Delta:=\,^k\pi^{-1}(\Delta(I\times\,\Sigma^{n-1}))
\label{liftoftube}
\end{align}
 is a tube in $J^k_0(M)$ with not intersections between the lifted curves in the sense that
  \begin{align*}
  \Delta(\cdot,\delta_1)\cap\,\Delta(\cdot,\delta_2)\,=\emptyset.
  \end{align*}
  Note that two curves can intersect at different times.
\begin{definicion}
Let $(M,\eta)$ be a pseudo-Riemannian manifold.
 A $k$-generalized $(p,q)$-tensor $\hat{S}$ is a map such that:
 \begin{enumerate}

 \item To each tube
$\Delta:I\times \,\Sigma^{n-1}\to M$ associates a
horizontal tensor along the lifted tube $^k\Delta$, which is a map
$\hat{S}:I\to T^{(p,q)}_hJ^k_0(\Delta)$.

\item If two tubes intersect, the value of the tensor $\hat{S}$ coincides in the intersection,
\begin{align}
 \hat{S}(\,^k\Delta_1)=\,\hat{S}(\,^k\Delta_2),\,\textrm{ for any  }\,^kx\in\,^k\Delta_1\,\cap\,^k\Delta_2.
\label{microlocality in tubes}
\end{align}

\end{enumerate}
Similarly, a $k$-generalized $p$-form is a map such that to each tube
$\Delta:I\times \,\Sigma^{n-1}\to M$ associates
a horizontal differential form along the lifted $^k\Delta$, which is a map
$\hat{S}:I\to \Lambda^p J^k_0(\Delta)$.
\label{generalizedfield}
\end{definicion}
 It is clear that the set of generalized tensors form a real
vector space, that we denote by $T^{(p,q)}_h(J^k_0 (M))$ (respectively, $\Lambda^p_h(J^k_0 (M)))$ for differential forms.
\begin{ejemplo}
If $k=0$, $\hat{S}$ is  a standard tensor or form defined in a tube $\hat{S}\in \,\Gamma \,T^{(p,q)}\Delta$
along a curve $x:I\to M$.
\end{ejemplo}
In the case that the manifold $(M,\eta)$ is Lorentzian, one can specify to work with tubes  composed by causal curves.
Also It is also interesting the possibility to work with congruences of curves, solutions of a given
differential equation on $M$. Then one needs to check that given
an arbitrary point $^kx\in\,J^k_0(M)$, it can be surrounded by a tube composed by lifting to $J^k_0(M)$
local solutions of the differential equation.
\begin{comentario}
The notion of tube in $J^k_0(M)$ that we use assumes that there are non-intersections in the curves composing the tube, at fixed time.
 This is natural from the point of view of locality, that becomes even a more restrictive notion than in usual geometry.
In this sense, one can speak about {\it local properties up to $k$-order}.
\end{comentario}

\subsection{Tensors and forms with values on $\mathcal{F}(J^k_0(M))$}
\begin{definicion}
A generalized tensor $T$ of type $(p,q)$ with values on $\mathcal{F}(J^k_0(M))$ is a smooth section of the bundle of $\mathcal{F}(M)$-linear homomorphisms
\begin{align*}
T^{(p,q)}(M,\mathcal{F}(J^k_0(M))) :=\,Hom(T^*M\times...^p...\times T^*M\times TM\times ...^q...\times TM,\,\mathcal{F}(J^k_0(M))).
\end{align*}
A $p$-form $\omega$ with values on $\mathcal{F}(J^k_0(M))$ is a smooth section of the bundle of $\mathcal{F}(M)$-linear completely alternate homomorphisms
\begin{align*}
\Lambda^p(M,\mathcal{F}(J^k_0(M))) :=\,Alt(TM\times ...^p...\times TM,\,\mathcal{F}(J^k_0(M))).
\end{align*}
The space of $0$-forms is $\Gamma\,\Lambda^0(M,\mathcal{F}(J^k_0(M))):=\mathcal{F}(J^k_0(M))$.
\label{definiciontensoerFJMvaluados}
\end{definicion}
Given a connection $^k\hat{\mathcal{H}}$ on $J^k_0(M)$,
there is an alternative way to describe generalized tensors (for both weak and strong definitions of generalized tensors and forms).
The connection $^k\hat{\mathcal{H}}$ determines the horizontal lifts
\begin{align*}
 T_xM\,\ni X\,\mapsto\,^hX\in\,T_{\,^kx}J^k_0(M)
 \end{align*}
and
\begin{align*}
T^*_xM\,\ni \alpha\,\mapsto\,^h\alpha \in\,T^*_{\,^kx}J^k_0(M).
\end{align*}
Then the tensor $\hat{S}\in  \Gamma\,T^{(p,q)}(M,\mathcal{F}(J^k_0(M)))$ defines a $\mathcal{F}(J^k_0(M))$-multilineal form
\begin{align*}
&\bar{S}:T^*M\times...^p...\times T^*M\times TM\times ...^q...TM  \to \mathcal{F}(J^k_0(M))
\end{align*}
by the rule
\begin{align}
& (\alpha_1,...,\alpha_p,\,X_1,...,X_q)\mapsto \hat{S}(\,^h\alpha_1,...,\,^h\alpha_p,\,^hX_1,...,
\,^hX_k,...,\,^hX_q).
\label{deficiniconhatSbarS}
\end{align}
This rule determines the homomorphism
\begin{align}
\phi^{-1}:T^{(p,q)}_h(J^k_0 (M))\to T^{(p,q)}(M,\mathcal{F}(J^k_0 (M))),\quad \hat{S}\,\mapsto \bar{S}.
\label{homomorphismtensorspaces}
\end{align}
There is a similar construction for generalized forms,
\begin{align}
\phi^{-1}:\Lambda^p_h(J^k_0 (M))\to \Lambda^p(M,\mathcal{F}(J^k_0 (M))),\quad \hat{\omega}\,\mapsto \bar{\omega}
\label{homomorphismgeneralizedforms}
\end{align}
such that it is an homomorphism of graded algebras,
\begin{align*}
\phi^{-1}(\alpha_1\wedge\,\alpha_2)=\,\phi^{-1}(\alpha_1)\,\wedge \phi^{-1}(\alpha_2),\,\alpha_i\in \,\Lambda^{p_i}_h(J^k_0(M)).
\end{align*}
The relation between definition \ref{definiciontensoerFJMvaluados} and this new alternative definition are given by the following two results,
 \begin{proposicion}
The definition of horizontal forms $\mathcal{F}(J^k_0(M))$-valued is connection independent.
\end{proposicion}
{\bf Proof}. Given two equivalent connections as the specified in {\it Theorem} \ref{teoremasobreconexiones},
 for any vector field $X\in \,\Gamma TM$, $^{h_1}X-\,^{h_2}X$ is vertical (similarly for the lifting of differential forms). Therefore,
 since the tensor $\hat{T}$ is horizontal,
 \begin{align*}
\hat\omega(\,^{h_1}X_1,...,\,^{h_1}X_q)
 =\,\hat\omega(\,^{h_2}X_1,...,\,^{h_2}X_q).
 \end{align*}\hfill$\Box$

Given a connection $^k\mathcal{H}$, it is possible to define the horizontal component of forms by using the projection
\begin{align}
\hat{h}_k: \Lambda^p (J^k_0(M))\to \Lambda^p (J^k_0(M))
\label{horizontalprojectionforforms}
\end{align}
that associates to each form $\hat{\omega}\in \Gamma \Lambda^p (J^k_0(M))$ its horizontal component. Thus we have the following {\it Proposition},
\begin{proposicion} Fixed a connection $^k\mathcal{H}$,
the homomorphisms (\ref{homomorphismtensorspaces}) and (\ref{homomorphismgeneralizedforms})
are $\mathcal{F}(J^k_0(M))$-isomorphisms.
\label{proposiciondelisomorphismode tensores}
\end{proposicion}
{\bf Proof}.
Let us fix our attention in the first homomorphism $(\ref{homomorphismtensorspaces})$.
Since $\hat{S}\in \,T^{(p,q)}_h(J^k_0 (M))$, it is horizontal and the fact that the horizontal
lift is an injective homomorphism, we have that (\ref{homomorphismtensorspaces}) is injective.
To prove that it is also surjective, one can consider an arbitrary
$\bar{S}\in\, \Gamma T^{(p,q)}(M,\mathcal{F}(J^k_0 (M)))$ defined at each point $^kx$ in the following way.
 For the horizontal $(\,^h\alpha_i,..., \,^h\alpha_p)$ $1$-forms
and horizontal $(X_1,...,X_p)$ vector fields on $M$ one has that,
\begin{align}
\hat{S}(\,^h\alpha_1,...,\,^h\alpha_p,\,^hX_1,...,\,^hX_q)(\,^kx):=
\bar{S}(\,^k\pi_*(\,^h\alpha_1),...,\,^k\pi_*(\,^h\alpha_p),\,^k\,\pi_*(\,^hX_1),...,\,^k\pi_*(\,^hX_q))(x);
\label{definitionofhatSfrombarS}
\end{align}
for any other vectors and forms on $J^k_0(M)$ one has
$\hat{S}(\beta_1,...,\beta_p,\,Y_1,...,Y_q)=0,$ {if any vector $Y$ or $1$-form $\beta$ is vertical}.
Then $\hat{S}$ is horizontal and such that $\phi(\hat{S})=\bar{S}$. Therefore, $\phi$ is also surjective.
Note that the computations are done pointwise.
The localization property follows in a similar way as in {\it Proposition 3.1} in \cite{KN}.
\hfill$\Box$

Because of the equivalence \ref{proposiciondelisomorphismode tensores}, we will use the two definitions of generalized field, as sometimes one could be more convenient than the other.
Any section of $Hom(T^*M\times...^p...\times T^*M\times TM\times ...^q...\times TM,\,\mathcal{F}(J^k_0(M)))$ is defined point-wise as
\begin{align*}
(\alpha_1,...,\alpha_p,X_1,...,X_q)\mapsto\,\hat{T}(\,^h\alpha_1,...,\,^h\alpha_p,\,^hX_1,...,\,^hX_q)\in \,\mathcal{F}(J^k_0(M),
\end{align*}
that is linear in any of its arguments, for any $x\in \,M$ and $\alpha_i\in\,T^*_xM$, $X_i\in\,T_xM$ and with $\hat{T}\in\,T^{(p,q)}J^k_0(M)$.
Similarly, any section $\bar{\omega}$ of $Hom(TM\times ,...^p,\times ,...,TM,\mathcal{F}(J^k_0(M)))$ is defined point-wise as
\begin{align*}
(X_1,...,X_p)\mapsto\,\hat{\omega}(\,^hX_1,...,\,^hX_p)\in \,\mathcal{F}(J^k_0(M),
\end{align*}
for a unique
such that it is alternate for any $x\in \,M$ and $X_i\in\,T_xM$ and with $\hat{\omega}\in \Lambda^p_hJ^k_0(M)$.

The isomorphism $\phi$ defined by (\ref{homomorphismgeneralizedforms}) induces the injection (on the image of $\phi$)
\begin{align}
^k\zeta:=\,\phi^{-1}:\Gamma\,\Lambda^p(M, \mathcal{F}(J^k_0(M)))& \to \Gamma\Lambda^p_h(J^k_0(M))
\label{definiciondezeta}
\end{align}
and is defined by the formula \eqref{deficiniconhatSbarS}.
In local coordinates  equation \eqref{definiciondezeta} is such that
\begin{align*}
^k\zeta:& \Gamma\Lambda^p(M, \mathcal{F} (J^k_0(M))) \to \Gamma\Lambda^p_h(J^k_0(M))\\
&\theta_I(\,^kx)\,dx^I\mapsto\theta_I(\,^kx)\,\phi^{-1}(dx^I),
\end{align*}
where $I$ is a multi-index. $^k\zeta$ is well defined, independent of coordinates and is an algebra homomorphism,
\begin{align*}
^k\zeta(\alpha_1\wedge\alpha_2)=\,^k\zeta(\alpha_1)\wedge\,^k\zeta(\alpha_2).
\end{align*}
Note that $^k\zeta$ is an isomorphism on the image. Therefore,  one has that
\begin{align*}
& ^k\zeta^{-1}(\hat{\omega}_1\wedge \hat{\omega}_2)=\,^k\zeta^{-1}(\hat{\omega}_1)\wedge\,^k\zeta^{-1}(\hat{\omega}_2),\quad \hat{\omega}_1,\,\hat{\omega}_2\in
\zeta(\Gamma\Lambda^*(M,\mathcal{F}(J^k_0(M))))\subset \Gamma\Lambda^p_h(J^k_0(M)).
\end{align*}
\begin{proposicion}
Let $^k\mathcal{H}_1$ and $^k\mathcal{H}_2$ be two connections on $J^k_0(M)$. Then
\begin{align*}
^k\zeta_{1}=\, ^k\zeta_{2}.
\end{align*}
\label{zetaesconnection independente}
\end{proposicion}
{\bf Proof}. Writing the homomorphisms $^k\zeta_{1}$ and $^k\zeta_{2}$ in local coordinate, it is clear that it does not depend on the connection:
in local coordinates and using multi-index notation for a form $\omega=\,\theta_I\,d\hat{x}^I$, $^k\zeta_i, i=1,2$ are given by expressions of the form
\begin{align*}
^k\zeta_1(\theta_I(\,^kx)d\bar{x}^I)=\,\theta_I(\,^kx)d\hat{x}^I=\,^k\zeta_2(\theta_I(\,^kx)d\bar{x}^I).
\end{align*}\hfill$\Box$
\begin{corolario}
The projection $\hat{h}_k: \Lambda^p (J^k_0(M))\to \Lambda^p_h (J^k_0(M))$ does not depend on the connection $^k\mathcal{H}$.
\label{corolariosobreinvarianceofkh}
\end{corolario}
Thus, for two different connections $^k\mathcal{H}_1$ and $^k\mathcal{H}_2$ we have the same projection operator, $\hat{h}^1_k(\hat{\omega})=\,\hat{h}^2_k(\hat{\omega})$.
\begin{comentario} We have the following remarks:
\begin{itemize}
\item The generalized higher order fields considered are smooth
sections of the generalized forms
\begin{align*}
\bar{\alpha}\in \Gamma\,\Lambda^p(M, \mathcal{F}(J^k_0 (M)))
\end{align*}
 or sections of a given generalized tensor bundle
 \begin{align*}
 \bar{T}\in \,\Gamma T^{(p,q)}(M,\mathcal{F}(J^k_0 (M))),
\end{align*}
for some specific natural number $k$, indicating the $k$-jet order where these fields take values.

\item When applied to classical electrodynamics, there is the same number of horizontal degrees of freedom that in the standard electromagnetic theory.

\item Generalized fields can be understood as weak generalized tensor or strong generalized tensors as was described before.
When one needs to speak about exterior derivatives or other differential operators acting on generalized higher order fields, they will be understood in a strong way.

\item The manifold $M$ is of relevance, since it is a requirement to define tensor along curves
or along tubes that contain an specific $x\in \,M$.
\end{itemize}
\end{comentario}

The {\it generalized tensor algebra} is
\begin{align}
T^*(M,\mathcal{F}(J^k_0 (M))):=\sum_{p,q}\oplus T^{(p,q)}(M,\mathcal{F}(J^k_0 (M))),
\label{tensoralgebra}
\end{align}
where the product is induced from the tensor product on
\begin{align*}
T^*(J^k_0(M))=\,\sum_{p,q}\oplus T^{(p,q)}(J^k_0 (M)).
\end{align*}

In a similar way, one defines the {\it generalized exterior algebra},
\begin{align}
\Lambda (M, \mathcal{F}(J^k_0 (M))):=\sum^n_{p=0}\oplus \Lambda^p(M, \mathcal{F}(J^k_0 (M))),
\label{tensorforms}
\end{align}
where the exterior product is induced from exterior product of the standard exterior algebra over $M$,
\begin{align*}
\Lambda (J^k_0 (M)))=\sum^n_{p=0}\oplus \Lambda^p((J^k_0 (M)).
\end{align*}
 Therefore, the product in the exterior algebra is defined as
\begin{align*}
\alpha(\,^kx)\wedge \beta(\,^kx)=(\alpha_I(\,^kx)e^I(x))\wedge\,(\beta_J(\,^kx)e^J(x)):=\,\alpha_I(\,^kx)\,\beta_J(\,^kx)\,e^I(x)\wedge e^J(x).
\end{align*}
Given a local frame $\{e_I(x),\,I=1,...,dim(\Lambda^p(M))\}$ for $\Lambda^p(M)$,
the homomorphism \eqref{homomorphismgeneralizedforms} implies that
a local frame for $\Lambda^p(M, \mathcal{F}(J^k_0 (M)))$ is obtained as the linear closure of$\{e_I(x),\,I=1,...,dim(\Lambda^p(M))\}$ with coefficients in $\mathcal{F}(J^k_0(M))$.

Given a $p$-form $\alpha\in \Lambda^pM$, there is a form $\varphi({\alpha})\in \Lambda^p(M, \mathcal{F}(J^k_0(M)))$
 such that for $(X_1,\,...,X_p) \subset \Gamma TM$ is defined by
 \begin{align}
  \varphi({\alpha})(X_1,\,...,\,X_p)=\alpha(X_1,...,X_p).
  \end{align}

  Therefore, one can establish the homomorphism of vector spaces
\begin{align}
\varphi:\Lambda^p M\to \Lambda^p(M,\mathcal{F}(J^k_0(M))),\quad  \varphi({\alpha})_u(X_1,...,X_p)=\alpha_x(X_1,...,X_p).
\label{embeddingforms}
\end{align}
 The value $\varphi(\alpha)_u(X_1,...,X_p)$ is constant along the fiber $\forall u\in\pi^{-1}(x)$.

\subsection{Generalized metric structures}

A relevant type of generalized tensor that we consider in some detail
are {\it generalized metrics}. For the purposes of this work
 a generalized metric $\bar{g}$ will be a weak generalized higher order tensor.
\begin{definicion}
A generalized metric $\hat{g}$ is a section $\hat{g}\in\,\Gamma\,T^{(0,2)}_h(J^k_0(M))$, a week generalized tensor such that:
\begin{enumerate}
\item It is smooth: given a curve $x:I\to\,M$ and a lift $^kx:I\to\,J^k_0(M)$,
 for any two smooth vector fields $\hat{X}_1$, $\hat{X}_2$ along $^kx:I\to J^k_0(M)$,
  the function $\hat{g}(^kx)(\hat{X}_1,\hat{X}_2):I\to \,R$ is smooth, except if $\hat{g}(^kx)(\hat{X}_1,\hat{X}_2)=0$.

\item It is homogeneous of degree zero in the following sense: if the lift $^kx:I\to M$ has local coordinates $(x^{\mu}(s),\dot{x}^{\mu}(s),\ddot{x}^{\mu}(s)...)$, then
\begin{align*}
\hat{g}(x^{\mu}(s),\lambda_1\dot{x}^{\mu}(s),\lambda^2_1\ddot{x}^{\mu}(s)+&\lambda_2\dot{x}^{\mu}(s),
\lambda_3\dot{x}+3\lambda_2\lambda_1\ddot{x}+\lambda^3_1\dddot{x},...)(X,X)=\\
&\hat{g}(x^{\mu}(s),\dot{x}^{\mu}(s),\ddot{x}^{\mu}(s),...\,x^{(k)}(s))(X,X)
\end{align*}
for all $\hat{X}$ along $x:I\to M$ and $\lambda_i>0,\,i=1,...,k$.

\item It is a symmetric form in the sense that
$\hat{g}(\,^kx)(\hat{X}_1,\hat{X}_2)=\,\hat{g}(\,^kx)(\hat{X}_2,\hat{X}_1)$ for any smooth pair
of vector fields $\hat{X}_1,\hat{X}_2$ along $x:I\to\,M$.

\item It is bilinear in the sense that
\begin{align}
\hat{g}(\,^kx)(\hat{X}_1\,+\,f(\,^kx)\hat{X}_2,\hat{X}_3)=\,\hat{g}(\,^kx)(\hat{X}_1 ,\hat{X}_3)\,+\,f(\,^kx)\hat{g}(\,^kx)(\hat{X}_1 ,\hat{X}_3),
\label{bilinerityproperty}
\end{align}

with $\hat{X}_1,\,\hat{X}_2,\,\hat{X}_3$ vectors fields along the curve $^kx:I\to \,J^k_0(M)$ and $f\in\,\mathcal{F}(J^k_0(M)).$

\item It is {\it weak non-degenerate} in the following sense: if the condition
\begin{align}
\hat{g}(\,^kx)(\hat{X},\hat{Z})=0,
\end{align}
holds for any horizontal smooth vector field $\hat{Z}$ along any curve $^kx:I\to \,J^k_0(M)$,
then the horizontal field $\hat{X}=0$.
\end{enumerate}
\label{generalizedmetric}
\end{definicion}
\begin{comentario}
The homogeneity condition implies that $\hat{g}(\,^kx)(\hat{X},\hat{Z})$ is invariant under the action of the chain
 rule $\big(\frac{d}{d\sigma}\big)^k\,=\big(\frac{d \rho}{d\sigma}\frac{d}{d\rho}\big)^k$.
 Therefore, it is invariant under positive parameterization of the curve $x:I\to M$.
\end{comentario}

Given a generalized metric, there is defined a
generalized bilinear form $g_u$ at the point $u\in\,^k\pi^{-1}(x)$,
\begin{align*}
\hat{g}_u: &T_xM \times T_xM\to R,\quad (X_1,X_2)\mapsto\,g_u(\,^hX_1,\,^hX_2)
\end{align*}
 such that it is positive homogeneous in the sense that $g_u( a\,y,a\, y)=\,a^2\,g(y,y)$ for any $a\in\, R^+$ and symmetric. Note that $\hat{g}$ does not depend on the particular non-linear connection.
Then following {\it Proposition} \ref{proposiciondelisomorphismode tensores}, it is possible to interpret a generalized
metric $\mathcal{F}(J^k_0(M))$ as the bilinear form
\begin{align}
\bar{g} :TM\times \,TM\to R, \quad (X_1,X_2)\mapsto \,\hat{g}_u(\,^hX_1,\,^hX_2).
\label{metricageneralizadaonTM}
\end{align}
\begin{proposicion}
The generalized metric $\hat{g}\in\,\Gamma\,T^{(0,2)}_h J^k_0(M)$ has associated a unique section
$\bar{g}\in \,\Gamma\, T^{(0,2)}(M,\mathcal{F}(J^k_0(M)))$ such that conditions $1$ to $5$ in {\it definition} \ref{{generalizedmetric}} hold.
\begin{enumerate}
\item $\bar{g}$ is smooth in the sense that for all $X_1,X_2$ smooth vector fields along the curve $x:I\to M$,
 the function $\bar{g}(X_1,X_2)$ is smooth except when it takes the zero value.

\item It is homogeneous of degree zero: if $^kx:I\to M$ has local coordinates $(x^{\mu}(s),\dot{x}^{\mu}(s),\ddot{x}^{\mu}(s)...)$, then
\begin{align*}
\bar{g}(x^{\mu}(s),\lambda_1\dot{x}^{\mu}(s),\lambda^2_1\ddot{x}^{\mu}(s)+&\lambda_2\dot{x}^{\mu}(s),
\lambda_3\dot{x}+3\lambda_2\lambda_1\ddot{x}+\lambda^3_1\dddot{x},...)(X,X)=\\
&\hat{g}(x^{\mu}(s),\dot{x}^{\mu}(s),\ddot{x}^{\mu}(s),...\,x^{(k)}(s))(X,X)
\end{align*}
for all ${X}\in T_xM$ and $\lambda_i>0,\,i=1,...,k$.
\item $\bar{g}$ is symmetric, in the sense that
\begin{align*}
\bar{g}(X_1,X_2)=\bar{g}(X_2,X_1)
\end{align*}
for all $X_1, ,X_2$ smooth vector fields along the curve $x:I\to  M$.

\item It is bilinear,
\begin{align*}
\bar{g}(X_1+\,f(\,^kx)X_2,X_3)=\bar{g}(X_1,X_3)\,+f(\,^kx)\bar{g}(X_2,X_3),
\end{align*}
for all $X_1,X_2,X_3$ arbitrary smooth vector fields along $x:I\to M$ and $f\in \,\mathcal{F}(M)$.

\item It is non-degenerate, in the sense that if $\bar{g}(X,Z)=0$ for all $Z$ smooth along $x:I\to M$, then $X=0$.

\end{enumerate}
\label{equialentedefiniciondelggeneralizada}
\end{proposicion}
{\bf Proof}. By {\it Proposition} \ref{proposiciondelisomorphismode tensores}, it is enough to show that the above
properties follow one by one from the corresponding properties in definition \ref{generalizedmetric}. Let us check the non-degeneracy condition.
The other properties follow using similar calculations. If $\hat{g}$ is non-degenerate,
then $\bar{g}(X,Z)=\,\hat{g}(h(X),h(Z))=0$ for all $Z\in\,\Gamma\,TM$ implies that,
since $h:T_xM\to \mathcal{H}J^k_0(x)$ is an isomorphism of real vector spaces, $X(x)=0$ at any $x\in\,M$. \hfill$\Box$

From the bilinear property it follows that locally the generalized metric can be written as
\begin{align*}
\bar{g}=\,g_{\mu\nu}\,d_4 x^\mu\otimes d_4x^\nu.
\end{align*}
The forms $\{d_4 x^\mu,\,\mu=1,...,n\}$ act linearly on vectors and define an adequate basis.
\begin{ejemplo} We can consider the following examples of generalized metrics:
\begin{itemize}
\item A relevant example that will appear in later {\it sections} is a generalized metric of the form
\begin{align}
\,^\lambda\hat{g}=\lambda(\,^kx)\, \eta(x),
\label{definitionofglambda}
\end{align}
where $\lambda$ is a positive, homogeneous of degree zero element of  $\mathcal{F}(J^k_0(M))$ for some fixed positive integer $k$, $\eta$ is the Minkowski metric and $^kx:I\to J^k_0(M)$ the lift of the curve $x:I\to M$.

\item Generalized metrics are provided by Finsler structures  \cite{Asanov, Bao, BaoChernShen}.
In this case the vertical hessian $g_{ij}(u)$ of the Finsler function $F$ is smooth on $TM\setminus\{0\}$. This example is of positive signature.
\end{itemize}
\label{ejemplodemetricageneralizada}
\end{ejemplo}
Of particular interest for the considerations in this paper are the generalized metrics of the form (\ref{definitionofglambda}),
\begin{align*}
g(u)=\lambda(u)\,\eta(x),\,\, x=\,^k\pi(u)\,\subset J^k_0(M).
\end{align*}
where $\eta$ is a Lorentzian metric on $M$.
The exact form of the conformal factor $\lambda$ will be elucidated later in the context of geometries of maximal covariant acceleration.
Therefore, the bilinear form $g$ is not a (pseudo)-Riemannian metric, since  $\lambda$ is not a function on $M$ but a function on $J^2_0(M)$.
In particular, we assume that $\lambda\in \mathcal{F}(J^2_0(M))$ is invariant under local transformations of coordinates on $M$
\footnote{In order to simplify the treatment and the calculations, we are considering flat spacetime in the sense that $\eta$ is flat.
Therefore, $\eta$ is Minkowski and $M$ is a flat domain of ${\bf R}^n$ of dimension $n$.
In physical applications the action of the Lorentz group as transformation group is the standard one.
For $n=4$, the factor $\lambda$ will be proved that is Lorentz invariant.}.
\begin{definicion}
Fixed a curve $x:I \to M$, the signature of the generalized metric $\bar{g}$ at $u=\,^kx(s)$ is the signature
 of the symmetric bilinear metric $\bar{g}_u:T_x M\times T_x M \to R$.
\end{definicion}
Thus the signature of a generalized metric will depend on the curve $^kx:I\to M$ where it is evaluated. This implies a natural classification
 of curves $\mathcal{F}(I,J^k_0(M))$
 with the corresponding curves have the same signature: fixed a curve, $x:I\to M$, we denote by the
 {\it signature sector associated with} $x:I\to M$ to the set of all the generalized metrics with the same along
 corresponding curves $z:I'\to M$ with the same signature $\bar{g}$ on $x:I\to M$.
\begin{ejemplo}
If one considers the generalized metric $\,^\lambda\hat{g}$ \eqref{definitionofglambda} in the first case in the example
 \ref{ejemplodemetricageneralizada}: if $\lambda$ is positive, the metric \eqref{definitionofglambda} has the same signature than the metric $\eta$.
\end{ejemplo}

Given a curve $x:I\to M$ and a vector field $W$ along the curve, if the signature of $g$ is $(-1,1,...,1)$
 one can define a bilinear, positive definite,
symmetric form
\begin{align*}
{g}_+\in\,\Gamma T^{(0,2)}(M,\mathcal{F}(J^k_0(M)))
\end{align*}
 such that along the curve
\begin{align}
{g}_+(X,X)(u)=\,g(X,X)(u)\,-2\,\frac{g^2(X,W)}{g(W,W)}(u),\quad g(W,W)<0,
\label{positivebilinearform}
\end{align}
with $u=\,^kx(s)\in J^k_0(M),\quad X\in \Gamma TM.$
For a generalized metric of the form $g=\lambda \eta$, the bilinear, positive form \ref{positivebilinearform} is
\begin{align*}
{g}_+(X,X)(u)=\lambda (u)\,{\eta}_+(X,X)(x),\quad g(W,W)<0,\quad u\in J^k_0(M),\quad X\in \Gamma TM,
\end{align*}
with
\begin{align*}
{\eta}_+(X,X)(x)=\,\eta(X,X)(x)\,-2\frac{\eta^2(X,W)}{\eta(W,W)}(x).
\end{align*}
Both generalized Riemannian metrics ${\eta}_+$ and ${g}_+$ can be used to provide a definition
of length along the curve $x:I\to M$.

\subsection{Group of isometries of a generalized metric}

Let ${\psi}:J^k_0(M)\to J^k_0(M)$ be a diffeomorphism of $M$  and denote the tangent map by
${\psi}_*:TJ^k_0(M)\to TJ^k_0(M)$. Then one can define the following homomorphism of the tangent bundle
$TJ^k_o(M)$,
\begin{align*}
\hat{\psi}_x: T_x& M\to T_{\,^k\pi(\psi(u))}(M),\quad X\mapsto\, ^k\pi_*\circ (\psi_*)|_u\circ \,^{h}X,
\end{align*}
with $ ^k\pi(u)=x$ and $ ^k\pi_*$, $\psi_*$ the differential maps.
\begin{definicion}
An isometry of a generalized metric $\bar{g}$ is a diffeomorphism $\psi:J^k_0(M)\to J^k_0(M)$
 such that on any arbitrary curve $x:I\to M$ it holds that
\begin{align}
\hat{g}_{{\,^kx}}(X_1,X_2)=\,\hat{g}_{\psi(\,^kx)}(\hat{\psi}_x({X}_1),\hat{\psi}_x({X}_2)),
\label{isometry}
\end{align}
for any two vector fields along $x:I\to M$.
\label{definicionofisometry}
\end{definicion}
The following {\it Proposition} follows easily from the definition \ref{definicionofisometry},
\begin{proposicion}
The group of isometries of a generalized metric does not depend on the non-linear connection $^k\mathcal{H}$ used in the definition of the corresponding diffeomorphisms.
\end{proposicion}
{\bf Proof}. Since a generalized metric can be locally written as $\bar{g}=\,\bar{g}_{\mu\nu}(\,^kx)\,dx^\mu \otimes dx^\nu$,
it is clear that the isometry group does not depend on the connection.\hfill $\Box$
\begin{comentario}
The set of isometries $Iso(\bar{g})$ of the generalized metric $\bar{g}$ is a subgroup of the group of diffeomorphism $Diff(J^k_0(M))$.
However, it is not so clear if it is a Lie group.
\end{comentario}
\begin{ejemplo}
We can consider the{ isometry groups of the following generalized metrics:}
\begin{itemize}
\item Isometries of the metric of type $\,^\lambda\bar{g}(\,^kx)=\,\lambda(\,^kx)\eta(x)$. Let $Iso(\,^\lambda\bar{g})$ be the group of
 isometries of $g$, and let $iso(\lambda)$ the group of diffeomorphisms
    \begin{align*}
    iso({\lambda}):=\left\{\phi:J^k_0(M)\to J^k_0(M),\,\,s.t.\quad\,\lambda\circ\phi=\,\lambda\,\right\}.
    \end{align*}
    The isometry group of $\eta$ is $Iso(\eta)$. Then one can write
    \begin{align*}
    Iso(\,^\lambda\bar{g})=\,(Iso(\eta)\,\cap iso(\lambda)) \times\, \mathcal{D}
     \end{align*}
     where $\mathcal{D}$ is the multiplicative group of {\it generalized dilatations},
\begin{align*}
\mathcal{D}:=\big\{f:J^k_0(M)\to J^k_0(M),\,\,s.t.\quad \lambda \to f^{-1}\,\lambda,\,\eta\to f\,\eta \big\}.
\end{align*}
It is clear that $iso(\lambda)\subset Iso(\eta)$ and
      that it is a subgroup of $Iso(\eta)$.

In the particular case  when $\lambda$ is invariant under $Iso(\eta)$, one has that
\begin{align*}
Iso(\eta)\cap iso(\lambda)=Iso(\eta).
\end{align*}
 Thus in this case
the isometry group is
\begin{align*}
Iso(\,^\lambda\bar{g})=\, Iso(\eta)\times \mathcal{D}.
\end{align*}
In general,the group of
generalized dilatations $\mathcal{D}$ is not necessarily a Lie group of transformations of $M$. However, if $\mathcal{D}$ is a Lie group,
the full group $ Iso(\eta)\times \mathcal{D}$ is a Lie group too.

\item Isometries of a Finsler metric. In this case, the group of isometries is a Lie group \cite{DengHou}.
In particular, the group of transformations that preserve the Finsler function $F$ (see for instance \cite{Bao, BaoChernShen}
for standard notation in Finsler geometry) coincides with the group of transformations that leave the {\it Finslerian distance} invariant.
They also prove that the group of isometries is a differentiable manifold and a Lie group.
\end{itemize}
\end{ejemplo}

\subsection{Generalized g-dual isomorphisms}

Let us consider the pull-back bundle
\begin{align*}
\xymatrix{^k\pi^* TM \ar[d]_{\rho_1} \ar[r]^{\rho_2} &
 TM \ar[d]^{\pi}\\
J^k_0(M) \ar[r]^{^k\pi} &  M}
\end{align*}
Then a generalized metric  ${g}(u)$ (without ''hat'' from now on) can be understood as a  metric on the fiber $\rho^{-1}_1(u)\subset \, ^k\pi^*{ TM}$.
Each fiber is the pull-back is a vector space generated by the pull-back of the local frame on $M$
\begin{align*}
\{^k\pi^*_u(\tilde{e}_1(x)),...,\,^k\pi^*_u(\tilde{e}_n(x)),\,\tilde{e}_i(x)\,\in\,\Gamma T_U M.
\end{align*}
In such local frame
for the pull-back fibers, the fiber metric has components
 \begin{align*}
g(u)_{\mu\nu}:=\,g|_u(^k\pi^*_u(\tilde{e}_\mu(x)),\,^k\pi^*_u(\tilde{e}_\nu(x))).
\end{align*}

One can diagonalize the metric component matrix  $g_{\mu\nu}(u=\,^kx(s))$ along $x:I\to M$ using the Gram-Schmidt's procedure at each point.
 The difference with the Gram-Schmidt for pseudo-Riemannian
manifolds is that the matrix of the transformation lives on $J^k_0(M)$. In such orthonormal frame
\begin{align*}
&\big(g_u(\,^k\pi^*e_\mu(u),\,^k\pi^*e_\nu(u))\big)=\,diag\big(-1,...,1\big)
\end{align*}
 on $ \rho^{-1}_1(u=\,^kx(s))\subset\,^k\pi^* TM$ holds.
Therefore, the projection $\rho_2$ provides a local frame
\begin{align*}
\{\rho_2(\,^k\pi^*e_1(u)),...,\rho_2(\,^k\pi^*e_n(u))\}
\end{align*}
along the curve $x:I\to \,M$ such that the generalized metric $g(x(s))$ is diagonal,
 \begin{align}
g_u(e_\mu,e_\nu):=\,\big(g(u)\big)_{\mu\nu}=\delta_{\mu\nu}.
\end{align}
The inverse of the matrix $g(u)_{\mu\nu}$ has components $g(u)^{\mu\nu}$ and determines an element of
$\Gamma T^{(2,0)}(M,\mathcal{F}(J^k_0(M)))$ which is non-degenerate and symmetric.
The raising and lowering indices operations are defined using $g^{-1}$ and $g$ (if anything else is specified).
$g^{-1}$ denotes the inverse matrix components of $g$. Thus, if we fix a frame on $M$, one
obtains $g^{-1}=\,\big(g^{-1}\big)^{\mu\nu}\,e_\mu\otimes e_\nu$ such that $\big(g^{-1}\big)^{\mu\rho}\,g_{\rho\nu}=\,\delta^\mu_\nu$.

The generalized metric $g$ determines the isomorphism,
\begin{align}
g:\Gamma T^{(1,0)} (M,\mathcal{F}(M,J^k_0(M)))\to \Gamma T^{(0,1)}(M,\mathcal{F}(M,J^k_0(M))),\quad X\,\mapsto \,\tilde{X}:=g(X,\cdot).
\end{align}
In similar way, the generalized tensor $g^{-1}$ determines the following canonical isomorphism,
\begin{align}
g^{-1}:\Gamma \Lambda^1 (M,\mathcal{F}(M,J^k_0(M)))\to\Gamma T^{(1,0)}(M,\mathcal{F}(M,J^k_0(M))),\quad \omega\,\mapsto g^{-1}(\omega,\cdot).
\end{align}
and a lowering index operation
\begin{align}
\kappa:\,\Gamma\Lambda^2(M,\mathcal{F}(J^k_0(M)))\to \,\Gamma T^{(1,1)}(M,\mathcal{F}(J^k_0(M))),\quad \kappa(\alpha)(\,\omega,X):=
\,\alpha(g^{-1}(\omega,\cdot),X),
\label{definiciondekappa}
\end{align}
with $\quad X\in\,\Gamma TM,\,
\alpha\,\in \Gamma \Lambda^2(M,\mathcal{F}(J^k_0(M))),\,\omega\,\in \Gamma \Lambda^1(M,\mathcal{F}(J^k_0(M)))$.

\subsection{A generalized star operator}

In order to introduce a star operator $\star_g$ associated with the generalized metric $g$, we will use local coordinate
expressions (see for instance  \cite{GoeckelerSchucker}).
 The generalized metric $g$ determines a star  operator $\star_g$ acting on the algebra
 \begin{align}
 \Lambda^p(M,\mathcal{F}(J^k_0(M))):=\,\sum^n_{p=0}\oplus \Lambda^p(M, \mathcal{F}(J^k_0 (M)))
 \end{align}
  in the following way. Let $\{e^\mu(u),\,\mu=1,...,n\}$ be a local, orthonormal frame respect of $g$ defined as before.
   The Levi-Civita symbol is denoted by $\epsilon_{\mu_1... \mu_n}$. Then the $\star_g$ operator of the algebra
   $\Lambda^p(M, \mathcal{F}(J^k_0 (M)))$ is the $\mathcal{F}(J^k_0(M))$-multilineal map determined
   by the image on the elements $e^{\mu_1}\wedge\cdot\cdot\cdot \wedge e^{\mu_p}$,
\begin{align*}
\star_g :\Gamma\Lambda^p (M,\mathcal{F}(J^k_0(M)))\to \Gamma\Lambda^{n-p} (M,\mathcal{F}(J^k_0(M)))
\end{align*}
\begin{align}
& (e^{\mu_1}(\,^kx)\wedge\cdot\cdot\cdot \wedge e^{\mu_p}(\,^kx))\mapsto \epsilon_{\nu_1...\nu_n}\,g^{\mu_1 \nu_1}\,
\cdot\cdot\cdot g^{\mu_p \nu_p}e^{\nu_{p+1}}(\,^kx)\wedge\cdot\cdot\cdot \wedge e^{\nu_n}(\,^kx)
\label{hodgestaroperator}
\end{align}
and then it is extended to an arbitrary generalized form
\begin{align*}
\alpha=\alpha_{\mu_1....\mu_p}(\,^kx)\,e^{\mu_1}(\,^kx)\wedge\cdot\cdot\cdot
\wedge e^{\mu_p}(\,^kx)\,\in\Gamma \Lambda^p(M, \mathcal{F}(J^k_0 (M)))
\end{align*}
 by the multilineal property \cite{GoeckelerSchucker},
\begin{align*}
\star\,\alpha(\,^kx)=\,\alpha_{\mu_1....\mu_p}(\,^kx)\epsilon_{\nu_1...\nu_n}\,g^{\mu_1 \nu_1}\,
\cdot\cdot\cdot g^{\mu_p \nu_p}e^{\nu_{p+1}}(\,^kx)\wedge\cdot\cdot\cdot \wedge e^{\nu_n}(\,^kx).
\end{align*}
By direct computation in a orthogonal frame, one can prove the following generalization of the Hodge star operator as in standard Riemannian geometry,
\begin{proposicion}
\begin{align}
\star \, \star \alpha\,= (-1)^{p(n-1)+s(g)}\alpha,\quad \alpha\in \, \Gamma \Lambda^p(M,\mathcal{F}(J^k_0(M))),
\label{hodgestarsquare}
\end{align}
where $s(g)$ is the signature of the generalized metric $g$.
\end{proposicion}

If $M$ is a four dimensional manifold, the $\star_g$ operator on $2$-forms is invariant under conformal transformations of $g$.
Therefore, when acting on an ordinary differential $2$-form on ${ M}$, the operators $\star_g$ and $\star_{\eta}$
determined by $g=\lambda\,\eta$ and $\eta$ respectively coincide: if $\{\tilde{e}^{\nu},\,\nu=1,...,4\}$ is a local orthonormal dual frame for $\eta$ and $\{{e}^{\nu},\,\nu=1,...,4\}$ is a local orthonormal dual frame for $g$, the relation between both is the conformal relation
\begin{align*}
\tilde{e}^{\nu}=\,\lambda^{-1}\,{e}^{\nu},\,\nu=1,...,4.
\end{align*}
Thus, one has that
\begin{align*}
\star_g\,\alpha(\,^kx) & =\,\alpha_{\mu_1\mu_2}(\,^kx)\epsilon_{\nu_1...\nu_4}\,g^{\mu_1 \nu_1}\,g^{\mu_2 \nu_2}e^{\nu_{3}}(\,^kx)\wedge \,e^{\nu_{4}}(\,^kx)\\
 & = \,\alpha_{\mu_1\mu_2}(\,^kx)\epsilon_{\nu_1...\nu_4}\,\lambda^{-1}\eta^{\mu_1 \nu_1}\,\lambda^{-1}\eta^{\mu_2 \nu_2}\lambda \tilde{e}^{\nu_{3}}(\,^kx)\wedge \,\lambda \tilde{e}^{\nu_{4}}(\,^kx)\\
& = \,\alpha_{\mu_1\mu_2}(\,^kx)\epsilon_{\nu_1...\nu_4}\,\lambda^{-1}\eta^{\mu_1 \nu_1}\,\lambda^{-1}\eta^{\mu_2 \nu_2}\lambda \tilde{e}^{\nu_{3}}(\,^kx)\wedge \,\lambda \tilde{e}^{\nu_{4}}(\,^kx)\\
& = \,\alpha_{\mu_1\mu_2}(\,^kx)\epsilon_{\nu_1...\nu_4}\,\eta^{\mu_1 \nu_1}\,\eta^{\mu_2 \nu_2}\tilde{e}^{\nu_{3}}(\,^kx)\wedge \, \tilde{e}^{\nu_{4}}(\,^kx)=\,\star_\eta\,\alpha(\,^kx).
\end{align*}
\subsection{Length and proper time associated with generalized metrics}

Generalized metrics do not have well defined locally signature: given two tangent vectors $X_1,X_2\in \, T_xM$, the value of $g(X_1,X_2)$ will depend on the curve along  it is evaluated and therefore, the notion of signature depends at each point $x$ of the given curve. Therefore, it is interesting to restrict to sectors of curves, where the metric $g$ has a well defined signature along all the curves.
The main two relevant types of sectors are the following:
\begin{itemize}
\item If the generalized tensor $g$ is positive definite, a distance function on $M$ is defined in the following way.
 Let $x:[\sigma_1, \sigma_2]\to M$ be a smooth curve.
 Then the length $L[x]$ is
 \begin{align}
       L[x]=\int^{\sigma_2}_{\sigma_1} \sqrt{g|_{(x,\frac{dx}{d\sigma},\frac{d^2x}{d\sigma^2},...)}
       (\frac{dx(\sigma)}{d\sigma},\frac{dx(\sigma)}{d\sigma})}\,d\sigma.
\label{length}
       \end{align}
\begin{definicion} Let $M$ be a connected manifold.
The distance function between two points $p,\,q\in M$ is defined to be
\begin{align}
d(p,q)=\,\inf \big\{\,L[x],\,x:I\to M,\quad s.t.\quad \,x(\sigma_1)=p,\,x(\sigma_2)=q \big\}.
\label{distancepq}
\end{align}
Given a point $p$ and a submanifold $N\hookrightarrow M$, the distance between $p$ and $N$ is
\begin{align}
d(p,N):=\,\inf \big\{\,d(p,q),\,q\in N\,\big\}.
\end{align}
\label{distance}
\end{definicion}
The notions of boundness, metric balls, completeness etc... can be extended from Riemannian geometry to the category of generalized metrics in a direct way.

\item If $g$ is of Lorentzian signature, some of the standard notions of causal curves can be
 translated in complete analogy from Lorentzian geometry \cite{BeemEhrlichEasly}. For instance, there is a well defined notion of causal curve:
\begin{itemize}

\item A vector field $X\in\,\Gamma \,TM$ is timelike if $g(X,X)<0$. A curve $x:I\to\, M$
 is timelike if the tangent vector field is timelike.

\item A vector field is lightlike if $g(X,X)=0$; a curve on $M$ is lightlike if the tangent vector field is lightlike.

\item A vector field $X\in\,\Gamma \,TM$ is spacelike if $g(X,X)>0$. A curve on $M$ is spacelike if the tangent vector field is spacelike.
 \end{itemize}
 The proper-time of a non spacelike curve is
 \begin{align}
       L[x]=:\int^{\sigma_2}_{\sigma_1} \,\sqrt{-g_{\,^kx(s)}(\frac{dx}{d\sigma},\frac{dx}{d\sigma})}\,d\sigma.
\label{length}
       \end{align}
The Lorentzian distance between two causal connected points (that be connected by a non spacelike curve) is
\begin{align}
d(p,q)=\,\sup \big\{\,L[x],\,x(\sigma_1)=p,\,x(\sigma_2)=q \big\}.
\end{align}
The proper-time parameter along a non spacelike curve $x:[\sigma_1, \sigma_2]\to M$ is the function
 \begin{align}
       \tau[\xi]=:\int^{\xi}_{\sigma_1} \,\sqrt{-g_{\,^kx(\sigma)}(\frac{dx}{d\sigma},\frac{dx}{d\sigma})}\,d\sigma.
\label{arc}
       \end{align}
\end{itemize}
The notion of length of a curve in the positive case and of proper time in the Lorentzian case have geometric meaning,
       \begin{proposicion}
Given a generalized metric structure $(M,g)$, the Weierstrass functional acting on an arbitrary curve $x:[\sigma_1,\sigma_2]\to M$,
 \begin{align}
       W[x]:=\int^{\sigma}_{\sigma_1} \,\sqrt{\big|g_{\,^kx(\sigma)}(\frac{dx}{d\sigma},\frac{dx}{d\sigma})\big|}\,d\sigma,
\label{Weirestrass}
       \end{align}
 is independent of the parametrization.
 \end{proposicion}
 {\bf Proof.} It is a consequence of the homogeneity condition of $g$.\hfill$\Box$

 Given a non spacelike curve $x:I\to M$, the proper time function
 $\tau[x](\sigma)$ (or arc-length function $L[x](\sigma)$ in the positive case) depends only on the point $x(\sigma)$ and the initial point $x(\sigma_1)$.
Therefore, one can use $\tau(\sigma)$ as a parameter of the curve.
\begin{ejemplo} We can consider the following examples:
\begin{itemize}
\item Let us consider the generalized metric $g(\,^kx)=\,\lambda(\,^kx)\,\eta(x)$, with $\eta$ a Lorentzian metric. Then if $\lambda$ is
 re-parametrization invariant, $g$ is re-parametrization invariant too and the proper-time
  (or arc-length) is also re-parametrization invariant.

\item The fundamental tensor of a Finsler spacetime structure \cite{Beem1} is homogeneous of degree zero and is defined the expression
\begin{align}
g_{\mu\nu}=\,\frac{\partial ^2 L(x,y)}{\partial y^\mu \partial y^\nu}.
\label{fundamentaltensor}
\end{align}
 When one evaluates the tensor $g$ on each curve $x:I\to M$ it determines a generalized metric $g(x,\dot{x})$. Also, it is clear that the proper time parameter is reparameterization invariant.
\end{itemize}
\end{ejemplo}
\begin{ejemplo}
 Let $g$ be a generalized metric of signature $(-1,1,1,1)$. For generalized metrics of the
  type $g=\,\lambda\,\eta$, with $(M,\eta)$ a geodesic complete Lorentzian spacetime and $0<\,\lambda$
   bounded on $M$, the causal theory of $\eta$ and $g$ are related. For instance, it is not difficult
   to see that the notions of global hyperbolicity for $(M,g)$ and $(M,\eta)$ coincide.
    A notion of time orientable is natural in this case.
\end{ejemplo}

The investigation of the critical points of the functional \eqref{Weirestrass} can be
 performed by standard techniques of analysis. The critical points of  \eqref{Weirestrass} provide important information on the local and global properties
 of geometric spaces. In addition, for specific
 generalized metrics (for instance, for metrics of maximal acceleration studied later), the signature sector is the same for all the geodesics.
\subsection{Cartan calculus for generalized forms}

In this subsection we generalize the fundamental notions of the Cartan calculus to the algebra of
generalized differential forms viewed as sections of $\Lambda^p(M, \mathcal{F}(J^k_0 (M)))$.
\subsection*{The interior product homomorphism for generalized forms}

The interior product of the exterior algebra $\Lambda M$ is denoted by
\begin{align*}
\iota_X:\Gamma\Lambda^pM\to \Gamma\Lambda^{p-1}M,\quad \,\Lambda^pM\ni\,\alpha\mapsto \iota_X\alpha\in\,\Lambda^{p-1}M,
\end{align*}
with $X\in \Gamma TM$ and defined by the relation
\begin{align*}
(\iota_X\alpha)(X_1,...,X_{p-1}):=\,\alpha(X,X_1,...,X_{p-1}),\,X,X_1,...,X_{p-1}\in \Gamma TM.
\end{align*}
can be extended to generalized forms,
\begin{definicion}
 Given the algebra of generalized differential forms $ \Lambda\,\Gamma(M, \mathcal{F}(J^k_0 (M)))$ and $X\in \Gamma TM$,
the interior product is the $\mathcal{F}(M)$-endomorphism
\begin{align}
\iota_X:\Gamma\Lambda^p(M, \mathcal{F}(J^k_0 (M)))\to \Gamma\Lambda^{p-1}(M, \mathcal{F}(J^k_0 (M))),\quad
\bar{\alpha}\,\mapsto\,^k\zeta^{-1}\,\iota_{\,^hX}\,^k\zeta\bar{\alpha}.
\label{interiorproduct}
\end{align}
\label{definitionofinteriorproduct}
\end{definicion}
\begin{proposicion}
The interior product $\iota_X$ has the following properties,
\begin{enumerate}
\item It is linear,
\begin{align*}
\iota_X(\bar{\alpha}_1+\,f(x)\bar{\alpha}_2)\,=\,\iota_X\,\bar{\alpha}_1\,+
f(x)\,\iota_X\,f(x)\,\bar{\alpha}_2,\quad \bar{\alpha}_i\in\,\Gamma\,\Lambda^{p_i}J^k_0(M),\,f\in\,\mathcal{F}(M),
\end{align*}

\item For every $X\in\, \Gamma TM$, $\iota_X\circ\,\iota_X\,\bar{\alpha}=0,\,\,\forall \bar{\alpha}\in\,\Gamma\,\Lambda^p(M,\mathcal{F}(J^k_0(M))).$

\item It is an skew-derivation; for any $\bar{\alpha}_i\in\,\Gamma\Lambda^{p_i}(M, \mathcal{F}(J^k_0 (M))),\,i=1,2$,
\begin{align*}
\iota_X\,(\bar{\alpha}_1\,\wedge \bar{\alpha}_2)=\,\iota_X\,\bar{\alpha}_1\,\wedge\,\bar{\alpha}_2\,+(-1)^{p_1}\,\bar{\alpha}_1\,\wedge\iota_X\,\bar{\alpha}_2.
\end{align*}
\item It is zero when acting on sections of $\Gamma \,\Lambda^0(M,\mathcal{F}(J^k_0(M)))$.
\end{enumerate}
\label{propiedadesdelproducto interior}
\end{proposicion}
{\bf Proof}.
For the linearity property, the proof is a direct computation:
\begin{align*}
\iota_X( \bar{\alpha}_1\,+f(x)\,\bar{\alpha}_2)& =\,^k\zeta^{-1}\,\iota_{\,^hX}\,^k\zeta\big(\bar{\alpha}_1\,+\,f(x)\bar{\alpha}_2\big)\\
& =\,^k\zeta^{-1}\,\iota_{\,^hX}\big(\,^k\zeta\bar{\alpha}_1\,+\,^k\zeta\,f(x)\,^k\zeta\bar{\alpha}_2\big)\\
& =\,^k\zeta^{-1}\,\big(\iota_{\,^hX}\,^k\zeta\bar{\alpha}_1\,+\,^k\zeta\,f(x)\iota_{\,^hX}\,^k\zeta\bar{\alpha}_2\big)\\
& =\,^k\zeta^{-1}\iota_{\,^hX}\,^k\zeta\bar{\alpha}_1\,+\,^k\zeta^{-1}\,^k\zeta\,f(x)\,^k\zeta^{-1}\iota_{\,^hX}\,^k\zeta\bar{\alpha}_2\\
& =\,\iota_X\,\bar{\alpha}_1\,+f(x)\,\iota_X\,f(x)\,\bar{\alpha}_2.
\end{align*}
For the second property follows from the standard interior product,
\begin{align*}
\iota_X\,\iota_X\,\bar{\alpha}& =\,^k\zeta^{-1}\,\iota_{\,^hX}\,^k\zeta\,^k\zeta^{-1}\,\iota_{\,^hX}\,^k\zeta\,\bar{\alpha}\\
& =\,^k\zeta^{-1}\,\iota_{\,^hX}\,\iota_{\,^hX}\,^k\zeta\,\bar{\alpha}\\
& =0.
\end{align*}
For the third property, if $\bar{\alpha}\in\,\Lambda^{p_1}(M,\mathcal{F}(J^k_0(M)))$,
\begin{align*}
\iota_X(\bar{\alpha}_1\wedge\,\bar{\alpha}_2)& =\,^k\zeta^{-1}\,\iota_{\,^hX}\,^k\zeta\,(\bar{\alpha}_1\wedge\,\bar{\alpha}_2)\\
& =\,^k\zeta^{-1}\,\iota_{\,^hX}\big(\,^k\zeta\,\bar{\alpha}_1\wedge\,\,^k\zeta\bar{\alpha}_2\big)\\
&=\,^k\zeta^{-1}\,\Big(\,\iota_{\,^hX}\,^k\zeta\,\bar{\alpha}_1\wedge\,\,^k\zeta\bar{\alpha}_2
\,+(-1)^{p_1}\,\,^k\zeta\,\bar{\alpha}_1\wedge\,\iota_{\,^hX}\,\,^k\zeta\,\bar{\alpha}_2\Big)\\
& =\,\iota_X\,\bar{\alpha}_1\,\wedge\,\bar{\alpha}_2\, +(-1)^{p_1}\,\bar{\alpha}_1\,\wedge\iota_X\,\bar{\alpha}_2.
\end{align*}
Similarly, the proof of the last property is a direct calculation: if $\bar{\alpha}\in \,\Lambda^0(M,\mathcal{F}(J^k_0(M)))$,
\begin{align*}
\iota_X\,\bar{\alpha}=\,^k\zeta^{-1}(\iota_{^hX}\,^k\zeta\bar{\alpha})=\,^k\zeta^{-1}(0)=0.
\end{align*}
\hfill$\Box$
\begin{proposicion}
The definition of the interior products (\ref{definitionofinteriorproduct}) is connection independent.
\end{proposicion}
{\bf Proof}. Given two horizontal distributions $^k\hat{\mathcal{H}}_1$ and $^k\hat{\mathcal{H}}_2$ as in {\it Theorem} \ref{teoremasobreconexiones},
 there are two definitions of interior products of
generalized forms,
\begin{align*}
 ^k\hat{\mathcal{H}}_1\,\Rightarrow \,^1\iota_X\,\bar{\alpha}=\,^k\zeta^{-1}_1\iota_{\,^{h_1}X}\,^k\zeta_1\,\bar{\alpha},\quad
 \, ^k\hat{\mathcal{H}}_2\,\Rightarrow \,^2\iota_X\,\bar{\alpha}=\,^k\zeta^{-1}_2\iota_{\,^{h_2}X}\,^k\zeta_2\,\bar{\alpha}.
\end{align*}
Comparing the two interior products,
\begin{align*}
^1\iota_X\,\bar{\alpha}\,-\,^2\iota_X\,\bar{\alpha}  & =\,^k\zeta^{-1}_1\iota_{\,^{h_1}X}\,^k\zeta_1\,\bar{\alpha}
-\,^k\zeta^{-1}_2\iota_{\,^{h_2}X}\,^k\zeta_2\,\bar{\alpha}\\
& =\,\,^k\zeta^{-1}_1\Big(\iota_{\,^{h_1}X}\,^k\zeta_1\,\bar{\alpha}
-\,^k\zeta_1\,^k\zeta^{-1}_2\iota_{\,^{h_2}X}\,^k\zeta_2\,\bar{\alpha}\Big)
\end{align*}
Using local coordinates, one can see that $^k\zeta_1\,^k\zeta^{-1}_2=Id$. Therefore,
\begin{align*}
^1\iota_X\,\bar{\alpha}\,-\,^2\iota_X\,\bar{\alpha}  &  =\,^k\zeta^{-1}_1\Big(\iota_{\,^{h_1}X}\,^k\zeta_1\,\bar{\alpha}\Big)
-\iota_{\,^{h_2}X}\,^k\zeta_2\,\bar{\alpha}\Big)
\end{align*}
On the other hand, $^{h_1}X-\,^{h_2}X=\,{\bf V}$
  is vertical (the easiest way to check this is using local coordinates).
 Then one gets
\begin{align*}
\hat{\alpha}(\,^{h_2}X,\cdot)-\,\hat{\alpha}(\,^{h_1}X,\cdot)=\,\hat{\alpha}(\,^{h_1}X,\cdot)-\,\hat{\alpha}(\,^{h_1}X,\cdot)=0.
\end{align*}\hfill$\Box$

\subsection*{The exterior covariant derivative for generalized forms}

The generalization of the exterior differential operator $d$ on smooth forms over $M$
\begin{align*}
d:\Gamma \Lambda^p M\to \Gamma\Lambda^{p+1}M
\end{align*}
to sections of $\Lambda^p(M, \mathcal{F}(J^k_0 (M)))$ will be a
 $\mathcal{F}(M)$-anti-derivation of the algebra $\Lambda(M, \mathcal{F}(J^k_0(M)))$ of degree $1$,
\begin{align*}
d_4: & \Gamma\Lambda^p(M, \mathcal{F}(J^k_0 (M)))\to \Gamma\Lambda^{p+1}(M, \mathcal{F}(J^k_0 (M))).
\end{align*}
In order to introduce a natural operator with the properties of the usual exterior differentiation, we
will consider the exterior derivative of forms defined on the jet bundle $J^k_0(M)$,
\begin{align}
d_J:\Gamma\Lambda^p(J^k_0(M))\to \Gamma\Lambda^{p+1}_h(J^k_0(M)).
\end{align}
and the horizontal component
$\hat{h}_k\big(d_J\,^k\zeta\,\alpha\big).$
\begin{definicion}
The exterior derivative operator of generalized forms is the homomorphism
\begin{align}
d_4: & \Gamma\Lambda^p(M, \mathcal{F}(J^k_0 (M)))\to \Gamma\Lambda^{p+1}(M, \mathcal{F}(J^k_0 (M))),\,\quad
\bar{\alpha}\mapsto\,^k\zeta^{-1}\,\hat{h}_k\big(d_J\,^k\zeta\,\bar{\alpha}\big).
\label{definitionded4}
\end{align}
\end{definicion}
  \begin{comentario}
  The definition of the exterior derivative $d_4$ only makes sense for strongly defined
  generalized forms, as in {\it Definition} \ref{generalizedfield} and for $\Sigma\hookrightarrow M$ of codimension zero.
  \end{comentario}
\begin{lema}
The
isomorphism \ref{definiciondezeta} implies the relations
\begin{enumerate}
\item
For any $p+1$ vector fields $\{X_1,...,X_{p+1}\}$ on $M$,  it holds that
\begin{align}
\,^k\zeta^{-1}\,\hat{h}_k\big(d_J\,^k\zeta\,\bar{\alpha}\big)(X_1,...,X_{p+1})=
\hat{h}_k\big(d_J\,\hat{\alpha}\big)(\,^hX_1,...,\,^hX_{p+1}).
\label{equaciondJhatX}
\end{align}
\item For any form $\hat{\alpha}\in\,\Lambda^p(J^k_0(M))$  and $p+1$ vector fields $\{X_1,...,X_{p+1}\}$ on $M$, it holds that
\begin{align}
d_J\hat{h}_k\hat{\alpha}(\,^hX_1,...,\,^hX_{p+1})=\,\hat{h}_k d_J\hat{\alpha}(\,^hX_1,...,\,^hX_{p+1}).
\label{commutativityhathdJ}
\end{align}
\end{enumerate}
\label{propiedadesdehathdJ}
\end{lema}
{\bf Proof}. The first property follows directly from the homomorphism \eqref{definiciondezeta} and from the formula \eqref{deficiniconhatSbarS},
\begin{align*}
\,^k\zeta^{-1}\,\hat{h}_k\big(d_J\,^k\zeta\,\bar{\alpha}\big)(X_1,...,X_{p+1}) & =\,^k\zeta^{-1}\,\big(\hat{h}_k\big(d_J\,\hat{\alpha}\big)\big)(X_1,...,X_{p+1})\\
& =\,\hat{h}_k\big(d_J\,\hat{\alpha}\big)(\,^hX_1,...,\,^hX_{p+1}).
\end{align*}

The second property requires the use of Cartan's formula for the exterior differential $d_J \hat{\alpha}$,
\begin{eqnarray*}
d_J\hat{\alpha}(\,^hX_1,...,\,^hX_{p+1})\,=\,\sum^{p+1}_{i=1}\,(-1)^{i+1}\,^hX_i(\hat{\alpha}\,(\,^hX_1,...,\,^h\widehat{X_i},...,\,^hX_{p+1}))
\end{eqnarray*}
\begin{eqnarray}
+\sum^{p+1}_{i,j=1}\,(-1)^{i+j+1}\hat{\alpha}\,([\,^hX_i,^hX_j],\,^hX_1,...,\,^h\widehat{X_i},...,\,^h\widehat{X_j},...,\,^hX_{p+1}).
\label{Cartanformula}
\end{eqnarray}
Note that $\hat{h}_k$ is linear. Therefore,
\begin{align*}
d_J\big(\hat{h}_k\big(\hat{\alpha}\big) \big)(\,^hX_1,...,\,^hX_{p+1})=\,\sum^{p+1}_{i=1}\,(-1)^{i+1}\,^hX_i(\hat{h}_k\hat{\alpha}\,(\,^hX_1,...,\,^h\widehat{X_i},...,\,^hX_{p+1}))\\
+\sum^{p+1}_{i,j=1}\,(-1)^{i+j+1}\hat{h}_k\hat{\alpha}\,([\,^hX_i,\,^hX_j],\,^hX_1,...,\,^h\widehat{X_i},...\,^h\widehat{X_j},...,\,^hX_{p+1}).
\end{align*}
Also,
the curvature vector field is
\begin{align}
R(X_i,X_j):=\,[\,^hX_i,\,^hX_j]-\,^h[X_i,X_j]
\label{acurvaturevectorfield}
\end{align}is a vertical vector field. Therefore,
\begin{align*}
\hat{h}_k\hat{\alpha}\,([\,^hX_i,\,^hX_j]-\,^h[X_i,X_j],\,^hX_1,...,\,^h\widehat{X_i},...\,^h\widehat{X_j},...,\,^hX_{p+1})=0
\end{align*}
and since $^hX_i$ are horizontal,
\begin{align*}
d_J\big(\hat{h}_k\big(\hat{\alpha}\big) \big)(\,^hX_1,...,\,^hX_{p+1}) &
 =\sum^{p+1}_{i=1}\,(-1)^{i+1}\,^hX_i\big((\hat{h}_k\hat{\alpha})\,(\,^hX_1,...,\,^h\widehat{X_i},...,\,^hX_{p+1}))\big)\\
& +\sum^{p+1}_{i,j=1}\,(-1)^{i+j+1}\hat{h}_k\hat{\alpha}\,(\,[\,^hX_i,\,^hX_j],\,^hX_1,...,\,^h\widehat{X_i},...\,^h\widehat{X_j},...,\,^hX_{p+1})\\
& =\sum^{p+1}_{i=1}\,(-1)^{i+1}\,^hX_i(\hat{\alpha}\,(\,^hX_1,...,\,^h\widehat{X_i},...,\,^hX_p))\big)\\
& +\sum^{p+1}_{i,j=1}\,(-1)^{i+j+1}\hat{h}_k\hat{\alpha}\,(\,^h[X_i,X_j],\,^hX_1,...,\,^h\widehat{X_i},...\,^h\widehat{X_j},...,\,^hX_{p+1}).
\end{align*}

Note that since
\begin{align*}
d_J(\hat{\alpha}\,(\,^hX_1,...,\,^h\widehat{X_i},...,\,^hX_p))(\,^hX_i)=(\,^hX_i)\cdot(\hat{\alpha}\,(\,^hX_1,...,\,^h\widehat{X_i},...,\,^hX_{p+1})),
\end{align*}
Cartan's formula can be re-written as
\begin{eqnarray*}
d_J\hat{\alpha}(\,^hX_1,...,\,^hX_{p+1})\,=\,\sum^{p+1}_{i=1}\,(-1)^{i+1}\,d_J(\hat{\alpha}\,(\,^hX_1,...,\,^h\widehat{X_i},...,\,^hX_{p+1}))(\,^hX_i)\\
+\sum^{p+1}_{i,j=1}\,(-1)^{i+j+1}\hat{\alpha}\,(\,^h[\,X_i,X_j],\,^hX_1,...,\,^h\widehat{X_i},...,\,^h\widehat{X_j},...,\,^hX_{p+1}).
\end{eqnarray*}

Now we compute the right hand side of \eqref{commutativityhathdJ}. Since $\hat{h}_k$ is linear, one has that
\begin{eqnarray*}
\hat{h}_k \,d_J\hat{\alpha}(\,^hX_1,...,\,^hX_{p+1})\,=\,\sum^{p+1}_{i=1}\,(-1)^{i+1}\,\hat{h}_k ( d_J(\hat{\alpha}\,(\,^hX_1,...,\,^h\widehat{X_i},...,\,^hX_{p+1})))(\,^hX_i)\\
+\sum^{p+1}_{i,j=1}\,(-1)^{i+j+1}(\hat{h}_k \hat{\alpha})\,([\,^hX_i,^hX_j],\,^hX_1,...,\,^h\widehat{X_i},...,\,^h\widehat{X_j},...,\,^hX_{p+1}).
\end{eqnarray*}
Each of the terms in the first line can be computed as
\begin{eqnarray*}
\hat{h}_k d_J(\hat{\alpha}\,(\,^hX_1,...,\,^h\widehat{X_i},...,\,^hX_{p+1}))(\,^hX_i) =
 d_J(\hat{\alpha}\,(\,^hX_1,...,\,^h\widehat{X_i},...,\,^hX_p))(\,^hX_i).
\end{eqnarray*}
As before, we use that $[\,^hX_i,\,^hX_j]-\,^h[X_i,X_j]$ is vertical. Therefore,
\begin{align*}
&\hat{h}_k\big(\sum^{p+1}_{i,j=1}\,(-1)^{i+j+1}(\hat{h}_k\hat{\alpha})\,([\,^hX_i,^hX_j],
\,^hX_1,...,\,^h\widehat{X_i},...,\,^h\widehat{X_j},...,\,^hX_{p+1})\big)\\
& =\,\sum^{p+1}_{i,j=1}\,(-1)^{i+j+1}(\hat{h}_k\hat{\alpha})(\,^h[X_i,X_j],\,^hX_1,...,\,^h\widehat{X_i},...,\,^h\widehat{X_j},...,\,^hX_{p+1})\\
& =\,\sum^{p+1}_{i,j=1}\,(-1)^{i+j+1}(\hat{\alpha})(\,^h[X_i,X_j],\,^hX_1,...,\,^h\widehat{X_i},...,\,^h\widehat{X_j},...,\,^hX_{p+1}).
\end{align*}
Gluing all together,
\begin{align*}
 \hat{h}_k d_J \hat{\alpha}(\,^hX_1,...,\,^hX_{p+1}) & =\,\,\sum^{p+1}_{i=1}\,(-1)^{i+1}\,\hat{h}_k d_J(\hat{\alpha}\,(\,^hX_1,...,\,^h\widehat{X_i},...,\,^hX_{p+1}))(\,^hX_i)\\
& +\,\sum^{p+1}_{i,j=1}\,(-1)^{i+j+1}(\hat{h}_k\hat{\alpha})(\,^h[X_i,X_j],\,^hX_1,...,\,^h\widehat{X_i},...,\,^h\widehat{X_j},...,\,^hX_{p+1})\\
& = \,\sum^{p+1}_{i=1}\,(-1)^{i+1}\,(\,^hX_i)\big( d_J(\hat{\alpha}\,(\,^hX_1,...,\,^h\widehat{X_i},...,\,^hX_{p+1})\big)\\
& +\,\sum^{p+1}_{i,j=1}\,(-1)^{i+j+1}(\hat{h}_k\hat{\alpha})(\,^h[X_i,X_j],\,^hX_1,...,\,^h\widehat{X_i},...,\,^h\widehat{X_j},...,\,^hX_{p+1}).
\end{align*}
This calculation proves the second property.\hfill$\Box$

From the prove of {\it Lemma} \ref{propiedadesdehathdJ}, one obtains an explicit formula for the exterior differential
$d_4 \bar{\alpha}\in\,\Lambda^{p+1}(M,J^k_0(M))$ acting on $p+1$ vector fields $X_1,...,X_{p+1}$ on $M$,
\begin{align*}
d_4\bar{\alpha}(X_1,...,X_{p+1})& =\,\sum^{p+1}_{i=1}\,(-1)^{i+1}\,^hX_i(\hat{\alpha}\,(\,^hX_1,...,\,^h\widehat{X_i},...,\,^hX_{p+1}))\\
& +\sum^{p+1}_{i,j=1}\,(-1)^{i+j+1}(\hat{h}_k\hat{\alpha})(\,^h[X_i,X_j],\,^hX_1,...,\,^h\widehat{X_i},...,\,^h\widehat{X_j},...,\,^hX_{p+1}).
\end{align*}

The operator $d_4$ is related with the exterior covariant derivative $D_{^k\mathcal{H}}$ \cite{KolarMichorSlovak},
\begin{align*}
(D_{^k\mathcal{H}}\omega)(X_1,...,X_{p+1})& :=\,(d_J\omega)(\,^hX_1,...,\,^hX_{p+1})\\
& =\,\sum^{p+1}_{i=1}\,(-1)^{i+1}\,^hX_i({\omega}\,\cdot (\,^hX_1,...,\,^h\widehat{X_i},...,\,^hX_{p+1}))\\
& +\sum^{p+1}_{i,j=1}\,(-1)^{i+j+1}{\omega}\,([\,^hX_i,\,^hX_j],\,^hX_1,...,\,^h\widehat{X_i},...,\,^h\widehat{X_j},...,\,^hX_{p+1}),\\
&\quad \omega\in \, \Gamma\,\Lambda^p(J^k_0(M)).
\end{align*}
For horizontal forms $\hat{\alpha}\in\,\Lambda^p_h(J^k_0(M))$, this formula reduces to
\begin{align*}
(D_{^k\mathcal{H}}\hat\alpha)(X_1,...,X_{p+1})& :=\,\sum^{p+1}_{i=1}\,(-1)^{i+1}\,^hX_i(\hat{\alpha}\,\cdot (\,^hX_1,...,\,^h\widehat{X_i},...,\,^hX_{p+1}))\\
& +\sum^{p+1}_{i,j=1}\,(-1)^{i+j+1}{\hat\alpha}\,(\,^h[X_i,X_j],\,^hX_1,...,\,^h\widehat{X_i},...,\,^h\widehat{X_j},...,\,^hX_{p+1}),
\end{align*}
 from which follows the
\begin{proposicion}
Given $\bar{\alpha}=\,^k\xi^{-1}(\hat{\alpha})\in\,\Lambda^p(M,\mathcal{F}(M,J^k_0(M))$, then it holds that
\begin{align}
d_4\bar{\alpha}(X_1,...,X_{p+1})=\,^k\xi^{-1}(D_{^k\mathcal{H}}\hat{\alpha})(X_1,...,X_{p+1}).
\end{align}
\end{proposicion}

\begin{proposicion}
The exterior differential operator \ref{definitionded4} has the following properties:
\begin{enumerate}
\item It is linear,
\begin{align*}
d_4( \bar{\alpha}_1\,+\lambda\,\bar{\alpha}_2)=\,d_4(\bar{\alpha}_1)\,+\lambda\,d_4\bar{\alpha}_2,
\quad \bar{\alpha}_i\in\,\Gamma\,\Lambda^{p_i}J^k_0(M),\,\lambda\in\,R.
\end{align*}

\item It is an skew-derivation; for any $\bar{\alpha}_i\in\,\Gamma\Lambda^{p_i}(M, \mathcal{F}(J^k_0 (M))),\,i=1,2$,
\begin{align*}
d_4\,(\bar{\alpha}_1\,\wedge \bar{\alpha}_2)=\,d_4\,\bar{\alpha}_1\,\wedge\,\bar{\alpha}_2\,+(-1)^p\,\bar{\alpha}_1\,\wedge d_4\,\bar{\alpha}_2.
\end{align*}

\item Acting on an arbitrary function $f\in \mathcal{F}(J^k_0(M))=\Gamma \Lambda^0(M\mathcal{F}(J^k_0(M)))$ is
expressed in local coordinates as
\begin{displaymath}
d_4(f)=\, \frac{\partial f}{\partial x^{\mu}}\,d_4x^{\mu},\,\,f\in \mathcal{F}(J^k_0(M)).
\end{displaymath}

\item It is $2$-steps nilpotent,
\begin{align}
d_4\,d_4\,\bar{\alpha}=0,\quad \forall \bar{\alpha}\in\,\Gamma\,\Lambda^p(M,\mathcal{F}(J^k_0(M))).
\end{align}
\end{enumerate}
\label{propiedadesded4}
\end{proposicion}
{\bf Proof}. The first property is proven after a short calculation, completely analogous to the proof of the linearity of the interior product, except that
$\lambda$ is a scalar instead of an arbitrary smooth function.
The second property follows from the following calculation. If $\bar{\alpha}_1, \bar{\alpha}_2\,\in \Gamma \Lambda^{p}(M,\mathcal{F}(J^k_0(M)))$, then
\begin{align*}
d_4(\bar{\alpha}_1\wedge\,\bar{\alpha}_2)& =\,^k\zeta^{-1}d_J(\,^k\zeta(\bar{\alpha}_1\wedge\,\bar{\alpha}_2)\\
& = \,^k\zeta^{-1}\hat{h}_k\big(d_J(\,^k\zeta(\bar{\alpha}_1)\wedge\,^k\zeta(\bar{\alpha}_2)\big)\\
& = \,^k\zeta^{-1}\Big(\hat{h}_k\big(d_J\,^k\zeta(\bar{\alpha}_1)\big)\wedge\,^k\zeta(\bar{\alpha}_2)\,
+(-1)^p\,\bar{\alpha}_1\,\wedge\hat{h}_k\big( d_J\,^k\zeta(\bar{\alpha}_2)\big)\Big)\\
& =\,d_4( \bar{\alpha}_1)\,\wedge \bar{\alpha}_2\,+(-1)^p\,\bar{\alpha}_1\,\wedge d_4(\bar{\alpha}_2).
\end{align*}
The third property is proved using local coordinates.
If $(U,x)$ is a local coordinate system on $M$, we have that
\begin{align*}
d_4 f(^kx)=\,^k\zeta^{-1}\, \hat{h}_k\big(d_J\,^k\zeta (f(\,^kx))\big).
\end{align*}
This relation is equivalent to
\begin{align*}
d_4 f(^kx) & =\,^k\zeta^{-1}\hat{h}_k\big(\frac{\partial f}{\partial x^{\mu}}\,^k\zeta (dx^{\mu})+\,\frac{\partial f(\,^kx)}{\partial y^A}\,\delta y^A\big)\\
& =\,^k\zeta^{-1}\frac{\partial f}{\partial x^{\mu}}\,^k\zeta (dx^{\mu})=
\,\frac{\partial f}{\partial x^{\mu}}\,^k\zeta^{-1}\,\hat{h}_k\big(\,^k\zeta (dx^{\mu})\big)\\
& =\,\frac{\partial f}{\partial x^{\mu}}d_4x^{\mu}.
\end{align*}
After a short calculation, using {\it Lemma} \ref{propiedadesdehathdJ} one obtains
\begin{align*}
d_4\,d_4 \bar{\alpha}(X_1,...,\,X_{p+2})& =\,^k\zeta^{-1}\,\hat{h}_k
\big(d_J\,^k\hat{h}_k\big(\,d_J\,^k\zeta(\bar{\alpha})\big)\big)(\,X_1,...,\,X_{p+2})\\
& =\,\hat{h}_k \,\hat{h}_k d_J\,\big(\big(\,d_J\,\hat{\alpha})\big)\big)(\,^hX_1,...,\,^hX_{p+2})\\
& =\,\hat{h}_k\,d_J\,\,d_J\,(\hat{\alpha})(\,^hX_1,...,\,^hX_{p+2})=0,
\end{align*}
since $d_J\,\,d_J\,(\hat{\alpha})=0$.
\hfill$\Box$
\begin{corolario}
For any $^k\zeta(\bar{\alpha})=\hat{\alpha}\in \,
\Lambda^p_h(J^k_0(M))$ the relation
\begin{align}
D_{^k\mathcal{H}}\circ D_{^k\mathcal{H}}\,^k\zeta\hat{\alpha}(X_1,...,X_{p+2})=0
\label{nilpotenceforDkH}
\end{align}
 holds.
\end{corolario}
{\bf Proof}. The nilpotent condition $d_4\,d_4=0$ can be re-written as
\begin{align*}
0=\,d_4 d_4\bar{\alpha} &= \,\,^k\zeta^{-1}D_{^k\mathcal{H}}\,^k\zeta \,\,^k\zeta^{-1}D_{^k\mathcal{H}}\,^k\zeta(\hat{\alpha})\\
& =\,\,^k\zeta^{-1}D_{^k\mathcal{H}}\,D_{^k\mathcal{H}}\,^k\zeta(\hat{\alpha}).
\end{align*}
This implies $0=\,D_{^k\mathcal{H}}\,D_{^k\mathcal{H}}\,^k\zeta(\hat{\alpha})$.
\hfill$\Box$
\begin{comentario}
It is well known that in general the exterior covariant derivative is not nilpotent \cite{KolarMichorSlovak}. However, the nilpotent property \eqref{nilpotenceforDkH} holds for only horizontal forms and not for arbitrary forms.
\end{comentario}
\begin{proposicion}
The exterior differential operator (\ref{definitionded4}) is connection independent.
\end{proposicion}
{\bf Proof}.
Given to connections $^k\hat{\mathcal{H}}_i,\quad i=1,2$, one can define two anti-derivations,
\begin{align*}
\bar{\alpha}\mapsto\, ^k\zeta^{-1}_1\hat{h}^1_k\big(d_J\,^k\zeta_1(\bar{\alpha})\big),\quad\bar{\alpha}\mapsto\, ^k\zeta^{-1}_2\hat{h}^2_k\big(d_J\,^k\zeta_2(\bar{\alpha})\big).
\end{align*}
We compare both anti-derivations as follows,
\begin{align*}
& ^k\zeta^{-1}_1\hat{h}^1_k\big(d_J\,^k\zeta_1(\bar{\alpha})\big)(X_1,...,X_{p+1})
-\, ^k\zeta^{-1}_2\hat{h}^2_k\big(d_J\,^k\zeta_2(\bar{\alpha})\big)(X_1,...,X_{p+1})=\\
&\,\hat{h}^1_k\big(d_J\,^k\zeta_1(\bar{\alpha}_)\big)(\,^{h_1}X_1,...,\,^{h_1}X_{p+1}) -\hat{h}^2_k\big(d_J\,^k\zeta_2(\bar{\alpha})\big)(\,^{h_1}X_1+V_1,...,\,^{h_1}X_{p+1}+V_{p+2}),
\end{align*}
for some vertical vector fields $V_i$ such that $^{h_2}X_i=\,^{h_1}X_i+\,V_i,\,i=1,...,p+1$.
Then by {\it Proposition} \ref{proposiciondelisomorphismode tensores} the second pieces above is
\begin{align*}
\hat{h}^2_k\big(d_J\,^k\zeta_2(\bar{\alpha})\big)(\,^{h_1}X_1+V_1,...,\,^{h_1}X_{p+2}+V_{p+1})=
\hat{h}^2_k\big(d_J\,\hat{\alpha}\big)(\,^{h_1}X_1,...,\,^{h_1}X_{p+1}).
\end{align*}
Note the relation (by {\it Corollary} \ref{corolariosobreinvarianceofkh})
\begin{align*}
\hat{\alpha}_2:=\hat{h}^2_k\big(d_J\,\hat{\alpha}_2\big)=\,\hat{h}^1_k\big(d_J\,\hat{\alpha}_2\big)=\hat{\alpha}_1.
\end{align*}
Thus one has that
\begin{align*}
\,^k\zeta^{-1}_2\hat{h}^2_k\big(d_J\,^k\zeta_(\bar{\alpha}_2)\big)(\,^{h_1}X_1+V_1,...,\,^{h_1}X_{p+2}+V_{p+1})=
\hat{h}^1_k\big(d_J\,\hat{\alpha}_2\big)(\,^{h_1}X_1,...,\,^{h_1}X_{p+1}).
\end{align*}
Again, note that $\hat{\alpha}_1-\hat{\alpha}_2$ is zero,
\begin{align*}
\hat{h}^1_k\big(d_J\,\hat{\alpha}_2\big)(\,^{h_1}X_1,...,\,^{h_1}X_{p+1})=\hat{h}^1_k\big(d_J\,\hat{\alpha}\big)(\,^{h_1}X_1,...,\,^{h_1}X_{p+1})
\end{align*}
from where the result follows.\hfill$\Box$

Thus, the exterior derivative $d_4$ does not depend on the covariant derivative that we use in its definition.
\begin{ejemplo}
A $1$-form $\phi\in\,\Gamma \Lambda^1(M, \mathcal{F}(J^k_0 M))$ can be written in local coordinates as
\begin{align*}
\phi=\phi_i(x(s),\,\dot{x},\,\ddot{x},..., x^{(k)})\,d_4x^i.
\end{align*}
Then its exterior derivative is
\begin{align*}
d_4\phi & =d_4\big(\phi_i(x(s),\,\dot{x},\,\ddot{x},..., x^{(k)})\, d_4x^i\big)\\
& = d_4\big(\phi_i(x(s),\,\dot{x},\,\ddot{x},..., x^{(k)})\big)\,\wedge d_4x^i+\,\big(\phi_i(x(s),\,\dot{x},\,\ddot{x},..., x^{(k)})\big)\wedge\,(d_4)^2\,x^i\\
& =\partial_j\,\phi_i(x(s),\,\dot{x},\,\ddot{x},..., x^{(k)})\,d_4x^j\,\wedge d_4x^i.
\end{align*}
It is clear that $d_4\phi$ does not depend on the local coordinate system on $M$.
If we calculate $d^2_4$ one obtains the expression
\begin{align*}
d_4(d_4\phi)= \partial_j\partial_l\,\phi_i(x(s),\,\dot{x},\,\ddot{x},..., x^{(k)})\,d_4x^l\wedge\,d_4x^j\,\wedge d_4x^i=0.
\end{align*}
\end{ejemplo}

It is natural to ask wether there is a generalization of Cartan's formula. Indeed, there is such natural generalization:
\begin{proposicion}
 The following formula holds:
\begin{align}
(\iota_X\,\circ d_4\,+d_4\,\circ\iota_X)(\bar{\alpha})=\,^k\zeta^k\,\hat{h}_k\Big(\mathcal{L}_{\,^hX}\,^k\zeta(\bar{\alpha})\Big).
\label{Cartanformula}
\end{align}
for all $\bar{\alpha}\in\, \Gamma\,\Lambda^p(M,\mathcal{F}(J^k_0(M))), \,\,X\in\,\Gamma \,TM$.
\end{proposicion}
{\bf Proof}. By direct calculation we have that
\begin{align*}
(\iota_X\,\circ d_4\,+d_4\,\circ\iota_X)(\bar{\alpha}) & =\,^k\zeta^{-1}\iota_{\,^hX}\,^k\zeta\,^k\zeta^{-1}\,\hat{h}_k\Big(d_J\,^k\zeta(\bar{\alpha})\Big)
+\,^k\zeta^{-1}\,\hat{h}_k\Big(d_J\,^k\zeta\,^k\zeta^{-1}\,\iota_{\,^hX}\,^k\zeta(\bar{\alpha})\Big)\\
& =\,^k\zeta^{-1}\iota_{\,^hX}\,\hat{h}_k\Big(d_J\,^k\zeta(\bar{\alpha})\Big)+\,^k\zeta^{-1}\,\hat{h}_k\Big(d_J\,\iota_{\,^hX}\,^k\zeta(\bar{\alpha})\Big)\\
& =\,^k\zeta^{-1}\hat{h}_k\Big(\iota_{\,^hX}\,d_J\,^k\zeta(\bar{\alpha})\Big)+\,^k\zeta^{-1}\,\hat{h}_k\Big(d_J\,\iota_{\,^hX}\,^k\zeta(\bar{\alpha})\Big)\\
& =\,^k\zeta^{-1}\,\hat{h}_k\Big(\mathcal{L}_{\,^hX}\,^k\zeta(\bar{\alpha})\Big).
\end{align*}\hfill$\Box$

An additional property that will be of relevance in the generalization of exterior differential equations to higher order forms is the following,
\begin{corolario}
Let  $\bar{f}\in\,\mathcal{F}(J^k_0(M))$ be a function such that is constant on $M$. Then in any coordinate system it holds the relation
\begin{align}
d_4 \bar{f}(y^{(1)\mu},y^{(2)\nu},...,y^{(k)\rho})=0.
\label{derivativeofverticalfunctions}
\end{align}
\end{corolario}

\begin{ejemplo}
It is natural to ask whether the above formula \eqref{Cartanformula} reduces to the known case when we take $k=0$ and therefore $J^0_0(M)\simeq M$.
In this case, $^hX=X$ and $\hat{H}_0=Id.$ Therefore, it is clear that equation (\ref{Cartanformula}) coincides with the usual Cartan's formula.
\end{ejemplo}
Other standard relations of the Cartan calculus can be extended to the algebra of generalized forms by similar considerations. Therefore, we have a complete Cartan's calculus for generalized forms.
\subsection{Vertical volume forms}

Let us consider a connection $^k\hat{\mathcal{H}}$ on $J^k_0(M)$.
We can construct a non-zero vertical $nk$-form
\begin{align*}
dvol_V(\,^kx)\in\,\Gamma\,\Lambda^{nk}J^k_0(M)
\end{align*}
 at each point of a local
trivialization of $^k\pi^{-1}(U)$, with $U\subset M$ being an open subset on $M$.
Given the local natural coordinate system $(x^{\mu},\,y^A,\,\mu=1,...,k,\,A=1,...,kn)$ on a local trivialization $U\times \,R^{nk}$ of $J^k_0(M)$,
the vertical volume form can be constructed using the non-linear connection $^k\hat{\mathcal{H}}$ and
the corresponding covariant vertical frame $\{\delta y^A,\,A=1,...,nk\}$ (\ref{verticalforms}) by the expression
\begin{align}
dvol_V(\,^kx):=\,w(\,^kx)\,\delta y^1\wedge\cdot\cdot\cdot\wedge \delta y^{nk}.
\label{volumeform}
\end{align}
on the trivialization $U\times\,R^{nk}$ of $J^k_0(M).$
\begin{proposicion}
The $nk$-form  $dvol_V$ has the following properties:
\begin{enumerate}
\item $dvol_V$ is a vertical form, $dvol_V(...,X,...)=0$ for any $X$ horizontal.

\item It is natural and globally defined on $J^k_0(M)$.
\end{enumerate}
\end{proposicion}
{\bf Proof}. The first statement follows from the fact that $dvol_V$  is constructed with vertical forms and the notion of being vertical is connection independent.
The second property is immediate using the connection $^k\hat{\mathcal{H}}$ on the fiber bundle $^k\pi:J^k_0(M)\to M$. \hfill$\Box$

More generally, one can introduce a general notion of vertical volume form,
\begin{definicion}
Given a connection $^k\hat{\mathcal{H}}$, a vertical volume form is a section of $\Lambda^{nk}J^k_0(M)$
such that it is zero when evaluated on any horizontal vector $\hat{X}$ and is non-zero everywhere.
\label{verticalvolumeform}
\end{definicion}
\begin{comentario}
The vertical volume form (\ref{volumeform}) is not connection independent. However, given two connections
$^k\hat{\mathcal{H}}_1$ and $^k\hat{\mathcal{H}_2}$, if $dvol_V(1)$ is
the vertical volume form constructed with $^k\hat{\mathcal{H}}_2$
and $\hat{X}(2)\in\Gamma\,^k\hat{\mathcal{H}}_2$, one has
\begin{align*}
dvol_V(1)(...,\hat{X}(1),...)=0,\quad
dvol_V(1)(...,\hat{X}(2),...)=0.
\end{align*}
\end{comentario}
\begin{ejemplo}
The vertical vector
sub-bundle $^k\tilde{\pi}_V: \,^k\mathcal{V}\to {M}$ embedded in a fiber orientable tangent bundle
$TJ^k_0(M)$ admits a vertical volume form. Let $\sigma_V:\,^k\mathcal{V}\to TJ^k_0(M)$
be a fiber preserving embedding \footnote{The fact that the vertical and horizontal distributions are in general globally defined implies that the embedding $\sigma_V:\,^k\mathcal{V}\to TJ^k_0(M)$ is globally defined as well.}. If $J^k_0(M)$ is orientable, let $dvol_J$ be a volume form on $J^k_0(M)$ and $^k\hat{\mathcal{H}}$  a connection such that $\{\,^hX_1,...,\,^hX_{n}\}$ generates locally a horizontal distribution.
 Then the form
 \begin{align*}
 \,\iota_{\,^hX_1}\cdot\cdot\cdot\iota_{\,^hX_n}\,dvol_J
 \end{align*}
 is a vertical volume form and the form
 \begin{align}
\sigma^*_V\,(\,\iota_{\,^hX_1}\cdot\cdot\cdot\iota_{\,^hX_n}\,dvol_J)
 \end{align}
  is a $nk$-volume form on $^k\mathcal{V}$.
 \end{ejemplo}
\begin{ejemplo}
We can consider the following constructions for vertical forms,
\begin{enumerate}
\item Let $E=M\times V$, with $V$ be a {\it vertical vector space} of dimension $d$ and basis $\{e_1,...,e_d\}$, whose dual basis is  $\{e^1,...,e^d\}$.
 A vertical form along the fiber is given by the exterior product of $1$-forms
\begin{align*}
dvol_V=\,e^1\wedge\,\cdot\cdot\cdot\wedge\,e^d.
\end{align*}
Then any top-form $dvol_V$ of degree equal to $dim(V)$ and non-zero everywhere determines a vertical top form on $E$ by $(\pi_E)^*(dvol_V)$, with $\pi_E:E\to M$.

\item  Let $\pi_E:E\to M$ be an orientable
vector bundle of dimension $n+d$ with $dvol_E$ a $n+d$-volume form on $E$. Let us also assume that $(M,\eta)$ is a flat pseudo-Riemannian manifold.
Then there is a globally defined orthonormal frame $\{e_1,...,e_n\}$ on $M$ and then
\begin{align*}
\iota_{e_1}\,\cdot\cdot\cdot\iota_{e_n}dvol_E
\end{align*}
is a globally defined vertical volume form.
\item An orientable fiber bundle $\mathcal{E}$ with an horizontal distribution $\mathcal{H}_{\mathcal{E}}$
admits a vertical volume form, constructed in a similar way as the volume form \eqref{verticalvolumeform}.

\end{enumerate}
\end{ejemplo}

The normalization function $w(\,^kx)$ in equation  (\ref{volumeform}) is defined by the normalization condition,
\begin{align}
\int_{\,^k\tilde{\pi}^{-1}_V(x)}\,dvol_V=1.
\label{normalizationvolumeform}
\end{align}
\begin{proposicion}
Given the vertical volume form $dvol_V$ as in (\ref{volumeform}),
there is a local frame $\{V_1,...,V_{nk}\}$ such that $dvol_V(V_1,...,V_{nk})=\,1$.
\end{proposicion}
{\bf Proof}. Let $\{\tilde{V}_1,...,\tilde{V}_n\}$ be an arbitrary local frame for the vertical bundle.
Since the volume form acting on a basis is different than zero and finite,
the result is obtained dividing by the factor $dvol_V({V}_1,...,{V}_{nk})\neq 0$ the first element $\tilde{V}_1$ of the local frame.\hfill$\Box$
\begin{ejemplo}
 If $(M,\eta)$ is flat, the existence of a vertical volume form on the bundle $J^k_0(M)$ is assured.
 Then a global holonomic frame along the fiber $^k\pi^{-1}(x)$ is $\{\frac{\partial }{\partial y^1},\cdot,...,\cdot\frac{\partial}{\partial y^{nk}}\}$
  and the vertical volume form can be written locally as
 \begin{align*}
dvol_V(\,^kx):=\,w(\,^kx)\delta y^1\,\wedge\cdot\cdot\cdot\wedge\,\delta y^{nk}.
\end{align*}
\end{ejemplo}

\subsection{Closed vertical volume forms}
We will consider vertical volume forms that are closed in the sense that
\begin{align}
 d_J(dvol_V)=0.
\label{closeconditionfortheverticalvolume}
 \end{align}
The existence of a closed vertical volume form implies some constraints on the connections. If we write down the condition \eqref{closeconditionfortheverticalvolume}, one gets an expression of the form
\begin{align*}
 d_J(dvol_V)=\,\frac{1}{w} \frac{\partial w}{\partial x^\rho}\,dx^\rho\wedge\, dvol_V+\,w\,d_J(\delta y^1\wedge\cdot\cdot\cdot\wedge \delta y^{nk}).
\end{align*}
Each of the exterior derivatives of the vertical $1$-forms $\delta y^A$ can be written as
\begin{align*}
d_J(\delta y^A)=\, d_J(d_J y^A +\, N^A\,_\rho\, dx^\rho)=\,d_Jd_J y^A+\, d_J(N^A\,_\rho)\wedge \,dx^\rho =\, d_J(N^A\,_\rho)\wedge \,dx^\rho .
\end{align*}
Each of the differentials
\begin{align*}
d_J (N^A\,_\rho)=\, \frac{\partial N^A\,_{\rho}}{\partial x^\sigma}\, dx^\sigma+\, N^A\,_{B\rho}\,\delta y^B.
\end{align*}
After some elementary algebra one ends up with the following to constraints,
\begin{itemize}
\item A constraint involving only the connection,
\begin{align}
\sum^{kn}_{A=1}(-1)^A\,\frac{\partial N^A\,_{\rho}}{\partial x^\sigma}\, dx^\sigma\wedge\,dx^\rho \wedge \,dvol_V(\frac{\partial}{\partial y^A},\cdot)=0.
\end{align}
This constraint is satisfied if one imposes the strong condition of $^k\mathcal{H}$ being symmetric,
\begin{align}
\frac{\partial N^A\,_{\rho}}{\partial x^\sigma}\,=\frac{\partial N^A\,_{\sigma}}{\partial x^\rho}.
\label{symmetricconnection}
 \end{align}
 This is a generalization of the torsion-free condition for the Levi-Civita connection.

\item The second constraint is on the volume form: it must satisfied the condition
\begin{align}
\frac{1}{w} \frac{\partial w}{\partial x^\rho}+\,\sum^{kn}_{A=1}\,{ N^A\,_{A\rho}}=0,\quad \rho=1,...,n.
\label{constrainonw}
\end{align}
Note that because of the definition of $ N^A\,_{B\rho}$, the components $\frac{\partial w}{\partial x^\rho}$ and ${ N^A\,_{A\rho}}$ transform in the same way under local coordinate transformations. Thus the equation \eqref{constrainonw} is covariant.
\end{itemize}
\begin{ejemplo}
Let us consider the case that
\begin{enumerate}
\item The $1$-form $\varpi=N^A\,_{A\rho}\,dx^\rho$ lives on $M$,

 \item The $1$-form is closed in the sense that $d\,\varpi=0$.
\end{enumerate}
 Then the partial differential equation \eqref{constrainonw} has a (local) solution: locally there is also a $0$-form $\upsilon$ such that $\varpi=d_4 \,\upsilon$.
\end{ejemplo}
Closed vertical volume forms are $nk$-calibration forms,
 \begin{proposicion}
 A closed vertical volume form $dvol_V$ on $J^k_0(M)$ is a $nk$-calibration in the sense of Harvey and Lawson  \cite{HarveyLawson}.
 \end{proposicion}
{\bf Proof}. Let $T_{nk}(x)$ be any $nk$-dimensional vector sub-space of $T_uJ^k_0(M)$.
If there is a $\xi\in\,T_{nk}(u)\cap\,\mathcal\,^k\mathcal{H}_u\neq \,\empty\emptyset$, then $dvol_V|_{T_{nk}}=0$. If $T_{nk}=\,^k\hat{\mathcal{V}}$,
 then there is a $nk$-vector $\xi_{nk}$ with $dvol_V(\xi_{nk})=1$.
 The result follows by linearity of the action of $dvol_V$ on $nk$-vectors.\hfill$\Box$
\subsection{The de Rham cohomology of differential forms $\Lambda^p(M, \mathcal{F}(J^k_0 (M)))$}

The fact that the operator $d_4$ is nilpotent implies the existence of some interesting cohomology theories,
\begin{definicion}
A $p$-form $\alpha \in \Lambda^p(M,\mathcal{F}(J^k_0 (M)))$ is closed if $d_4\alpha=0$; it is exact if there is a $(p-1)$-form
$\beta\in\Lambda^{(p-1)}(M,\mathcal{F}(J^k_0 (M)))$ such that $d_4\beta =\alpha$.
\label{closedexactforms}
\end{definicion}
The vector space of $d_4$-closed forms is
\begin{align*}
Z^p(M,\mathcal{F}(J^k_0(M))):=\left\{ \alpha \in\,\Gamma \Lambda^p(M,\mathcal{F}(J^k_0 (M)))\,s.t.\, d_4\alpha=0\right\}.
\end{align*}
The vector space of $d_4$-exact forms is
\begin{align*}
B^p(\mathcal{F}(M,J^k_0(M)))& :=\{ \alpha \in \,\Gamma \Lambda^p(M,\mathcal{F}(J^k_0 (M)))\\
& |\exists\,\,\beta\in \Gamma \Lambda^{p-1}(M,\mathcal{F}(J^k_0 (M)))\,s.t.\,d_4\beta=\alpha\}.
\end{align*}
The $p$-cohomology group is
\begin{align*}
H^p(M,\mathcal{F}(J^k_0(M))):=\, Z^p(M,\mathcal{F}(J^k_0(M)))/B^p(M,\mathcal{F}(J^k_0(M))).
\end{align*}
Therefore, the cohomology $H^*(M,\mathcal{F}(J^k_0(M)))$ is defined from the differential structure of the manifold $M$.

One can also define
the vertical compact cohomology $p$-cohomology group of the $k$-jet bundle $J^k_0(M)$,
\begin{align*}
H^p_{cv}(M,\mathcal{F}(J^k_0(M))):=&\{\alpha\in \Gamma\Lambda H^p(M,\mathcal{F}(J^k_0(M)))\,\mid\, \alpha_x \\
&\textrm{has compact support on each fiber $^k\pi^{-1}(x)$}\,\}.
\end{align*}
Then the compact vertical cohomology $H^*_{cv}(M,\mathcal{F}(J^k_0(M)))$ is defined in a similar way.
\begin{comentario}
We will prove later (see below {subsection} 2.16) that $H^*_{cv}(M,\mathcal{F}(J^k_0(M)))$ is invariant under homotopy: if $M$ and $\tilde{M}$ are homotopy equivalent, then $H^*_{cv}(M,\mathcal{F}(J^k_0(M)))\simeq H^*_{cv}(\tilde{M},\mathcal{F}(J^k_0(\tilde{M})))$.
\end{comentario}

\subsection{Bounded vertical compact cohomology}

A further restriction in the cohomology forms is the following,
\begin{definicion}
Let $(M,\eta)$ be a Lorentzian manifold.
Then the covariantly bounded $k$-jet bundle is the fiber bundle $^k\pi:J^k_{0b}(M)\to M$ { such that}
\begin{align*}
& 1. \,\textrm{The curves}\,\,x:I\to M,\, x(0)=p\,\,\textrm{are smooth},\\
& 2. \textrm{ Each of the covarian derivatives is bounded},\\
& |\eta(D^i_{\dot{x}}\dot{x},D^i_{\dot{x}}\dot{x})|\leq\,c_{i+1}\in\,R^+,\,=1,...,k-1,
\end{align*} where $D$ is the Levi-Civita connection of $\eta$.
\label{definicionofboundedcohomology}
\end{definicion}
Note that the values of $c_i$ could be different for each $i=1,...,k$ and some of the values could be un-bounded, but by definition at least one of the constants $c_i$ should be finite. Thus the specification of the bounded jet bundle $J^k_{0b}(M)$ will be given by a finite collection of constants $\{(c_1,...,c_k)\}$.
\begin{ejemplo}The physical interpretation of {\it Definition} \ref{definicionofboundedcohomology}is better understood after the following two different examples,
\begin{enumerate}
\item For $k=1$, the second condition is equivalent to the requirement of a {\it maximal speed}. Note that it only indicates a bound on the value of the allowed covariant speed. In this case $c_1=1$ in natural coordinates and is $J^1_{0b}$.

\item For $k=2$ and for time like curves of $\eta$, the bound $c_2=A^2_{max}<\infty$ is equivalent to the requirement of maximal acceleration. Note that one can have $c_1=\infty$, dropping the requirement of maximal covariant speed. It corresponds to $J^2_0{0b}(M)$ with $(c_1=\,\infty,c_2<A^2_{max}$.

\item Other example of relevance is $k=3$ with $(c_1<\infty,c_2<\infty,c_3<\infty)$. The bounded jet bundle with such specification will be denoted by $J^3_{0b}(M)$. This bundle is of relevance for the electromagnetic theory that we will develop in this work.
\end{enumerate}
\end{ejemplo}

We can consider the following cohomology theory,
\begin{definicion}
Let $(M,\eta)$ be a Lorentzian manifold and consider the fiber bundle $J^k_{0b}(M)$.
The covariantly bounded compact vertical cohomology is
\begin{align*}
H^p_{cv}(M,\mathcal{F}(J^k_{0b}(M)))& :=\{\alpha\in \Gamma\Lambda^p(M,\mathcal{F}(J^k_{0b}(M)))\\
& \mid\, \alpha_x \,\textrm{has compact support on each fiber $^k\pi^{-1}(x)$}\,x\in M\}.
\end{align*}
\label{cohomologyofcompact verticalboundedforms}
\end{definicion}
In a similar way, one can consider the de Rham cohomologies
$H^*_{dR}(J^k_0(M))$ and $H^*_{cv}(J^k(M))$.
\begin{comentario}
We have the following remarks:
\begin{itemize}

\item  Note the restriction to $J^k_{0b}(M)$ and the associated cohomology
$H^*(M,\mathcal{F}(J^k_{0b}(M)))$, even for non-compact vertical forms.

\item The cohomologies $H^*(M,\mathcal{F}(J^k_0(M)))$,
 $H^*_{cv}(M,\mathcal{F}(J^k_{0}(M)))$, $H^*_{cv}(M,\mathcal{F}(J^k_{0b}(M)))$
 and $H^*(J^k_{dR}(M))$ are  in principle different from each other.

 \item From the different types of cohomologies that one can construct with generalized differential forms, the relevant four our applications in generalized field theory is $H^*_{cv}(M,\mathcal{F}(\tilde{J}^2_{0b}(M)))$. We will prove in subsection 2.16 that $H^*_{cv}(M,\mathcal{F}(\tilde{J}^2_{0b}(M)))$ is homotopy invariant of $M$, for the set of constants $\{c_1<\infty,c_2<\infty,c_3>\infty\}$.
\end{itemize}
\end{comentario}

The bounded cohomology is not only defined on jet bundles. For instance, it is possible to define it for homotopy equivalent bundles
\begin{align}
E^k\hookrightarrow J^{k'}_{0}(M),
\end{align}
for a given $k$, such that there is an $J^{k'}_{0b}(M)$ with the homotopy equivalence
\begin{align}
E^k\simeq\,J^{k}_{0b}(M).
\end{align}

For applications in classical field theory, the bundle $\tilde{J}^2_{0b}$ defined as follows will appear as the support where physical world-lines will life. Indeed, we will provide a formal definition of the bundle, for further use,
\begin{definicion}
The bundle  $\tilde{J}^3_{0b}(M)\to M$ is a subbundle of $J^3_{0}(M)\to M$ such that the fibers correspond to the lifts $^3x:I\to M$ of the smooth curves $x:I\to M$ to $J^3_0(M)$ with the following properties:
\begin{align*}
& 1. \,\textrm{The curves}\,\,x:I\to M,\, x(0)=p\,\,\textrm{are smooth},\\
& 2. |\eta(D_{\dot{x}}\dot{x},D_{\dot{x}}\dot{x})|\leq\,c_{i}\in R^+,\,\,i=1,2,\\
& 3. |\eta(D^3_{\dot{x}}\dot{x},D^2_{\dot{x}}\dot{x})|^{-1}\leq\,{c}_{3}\in\,R^+,
\end{align*} where $D$ is the Levi-Civita connection of $\eta$.
\label{definicionoftildej20b}
\end{definicion}
Later, we will see how to relate the cohomologies $H^*_{cv}(M,\mathcal{F}(J^2_{0b}(M)))$ and $H^*_{cv}(M,\mathcal{F}(\tilde{J}^3_{0b}(M)) )$.
\subsection{Thom's isomorphism theorem and homotopy invariance of\\
 $H^*_{cv}(M,\mathcal{F}(J^k_{0b}(M)))$}

In order to prove the homotopy invariance of $H^*_{cv}(M,\mathcal{F}(J^k_{0b}(M)))$, one strategy is to relate the cohomology theory $H^*_{cv}(M,\mathcal{F}(J^k_{0b}(M)))$ with the de
Rham cohomology  $H^*_{dR}(M,R)$.
This relation arises
as an application of {\it fiber integration} \cite{Botttu},
which is also fundamental for the proof of {\it Thom's isomorphism theorem}.
Before doing this, we note that  there is also an
 homotopy between $J^k_{0b}(M)$ and $J^{k-1}_0(M)$.
\begin{proposicion}
Let $(M,\eta)$ be a Lorentzian manifold and $c_k$ finite. Then $J^k_{0b}(M)$ and $J^{k-1}_0(M)$ have the same homotopy class.
\label{propositiononhomotopy}
\end{proposicion}
 {\bf Proof}. Since the fibers of $J^k_{0b}(M)$ are bounded, there is a finite time contraction: for each point $(x(\sigma),\, ^1y(\sigma),\,^2y(\sigma),...,\,^ky(\sigma))$, one defines the contraction
 \begin{align}
 (x(\sigma),\, ^1y(\sigma),\,^2y(\sigma),...,\,^ky(\sigma))\,\rightarrow (x(\sigma),\, ^1y(\sigma),\,^2y(\sigma),...,\,^{k-1}y,\lambda(\,^ky,\,^ky_{max}) \,^ky),
 \end{align}
with $\lambda(\,^ky,\,^ky_{max})\in [0,1]$ is a retraction to the origin of the $k$-derivatives coordinates. Since the retraction is in finite time and bounded, we can consider a globally defined retraction which is continuous and determined by the maximal retraction time for each of the fibers over $x\in \,M$.\hfill$\Box$
\begin{ejemplo}
{\it Proposition} \ref{propositiononhomotopy} admits straightforward generalizations, as the following example shows. Let us consider $\tilde{J}^3_{0b}(M)$ with the set o f constants $(c_1=\infty,c_2<\infty,c_3<\infty)$. Then the following homotopies hold,
\begin{align}
\tilde{J}^3_{0b}(M)\simeq\,J^2_{0b}(M)\simeq J^1_0(M)\simeq TM.
\label{exampleofhomotopyequivalence}
\end{align}
\end{ejemplo}
 As a direct consequence of {\it Proposition} \ref{propositiononhomotopy} we have the following
\begin{proposicion}
Let $(M,\eta)$ be a Lorentzian manifold. Then
\begin{enumerate}
\item It holds the following isomorphism:
\begin{align}
H^*(J^k_{0b}(M))\simeq\,H^*(J^{k-1}_{0}(M)).
\label{isomorphismJkJk-1}
\end{align}
\item It holds the following isomosphism:
\begin{align}
H^*_{cv}(J^k_{0b}(M))\simeq\,H^*_{cv}(J^{k-1}_{0}(M)).
\label{isomorphismJKJk-1verticalbounded}
\end{align}
\end{enumerate}
\label{homotopiaentrejetbundles}
\end{proposicion}
{\bf Proof}. We construct a homotopy between $J^k_{0b}M$ and $J^{k-1}_0(M)$.
Since both cohomologies $H^*_{cv}$ and $H^*$ are invariant under homotopy, the result follows.\hfill$\Box$
\begin{corolario}
The following relation holds,
\begin{align}
H^p_{cv}(J^k_{0b})=0,\,p=kn+1,...,(k+1)n.
\end{align}
\end{corolario}
{\bf Proof}. It is enough to consider the isomorphism \eqref{isomorphismJKJk-1verticalbounded} for the cohomology groups
\begin{align*}
H^p_{cv}(J^{k-1}_0(M)),\,\,\textrm{for}\, p>dim(J^{k-1}_0(M)),
\end{align*}
which are identically zero. \hfill $\Box$

As a consequence of \eqref{isomorphismJKJk-1verticalbounded} and {\it Example} \eqref{exampleofhomotopyequivalence} we have the following
\begin{corolario}
Let us consider $\tilde{J}^3_{0b}(M)$ with $(c_1=\infty,c_2<\infty,c_3<\infty)$. Then the following isomorphisms hold,
\begin{align}
H^*_{cv}(\tilde{J}^3_{0b}(M))\simeq H^*_{cv}(J^2_{0b}(M))\simeq\,H^{*}_{cv}(TM).
\end{align}
\label{corolariosobreisomorfismodecohomologias}
\end{corolario}
\begin{comentario}
{\it Corollary} \ref{corolariosobreisomorfismodecohomologias} can be extended easily to $k>1$.
Based on this {\it Corollary}, from now on we will contemplate vertical bounds of the form $\tilde{J}^3_{0b}(M)$ with $\{c_1=\infty,c_2<\infty,...,c_k<\infty\}$. This is indeed the jet bundle that we will find in the generalized electrodynamic theory of point charged particles.
\end{comentario}
\begin{definicion}
Let $\pi_E:E\to B$ be a vector bundle with orientable fiber over
with $\pi^{-1}_B(x)$ of dimension $l$. Then the averaging along the fiber of a $l+k$ form is the homomorphism
\begin{align}
\langle\cdot\rangle :  \Gamma \Lambda^{l+k}E \to \Gamma \Lambda^k B,\quad \alpha\mapsto \int_{\pi^{-1}_B(x)}\,\alpha.
\label{averagingforE}
\end{align}
\end{definicion}
Given a vector bundle $\pi_E:E \to B$, let us consider the cohomology
of the differential form $\Lambda^*(E)$ with compact vertical support $H^*_{cv}(E)$.
Then one has the following isomorphism \cite{Botttu},
\begin{teorema}
Let $\pi_E:E\to M$ be a vector bundle of finite type over a $n$-dimensional manifold $M$.
Then if the dimension of the fibers $\pi^{-1}_B(x)$ is $l$, one has the isomorphism
\begin{align}
H^*_{cv}(E)\,\simeq H^{*-l}_{dR}(M).
\label{Thomisomorphism forE}
\end{align}
\label{Thomstheorem}
\end{teorema}

Let us construct the following vector bundle structure on jet bundles. We consider a symmetric connection $\mathcal{N}$ on $J^k_0(M)$. Then to each $k$-jet $(x,\dot{x},\ddot{x},...,x^{(k)})$ one associates the corresponding $(x,D_{\dot{x}}\dot{x},D^2_{\dot{x}}\dot{x},...,D^{(k)}_{\dot{x}}\dot{x})$, where $D$ is covariant the covariant derivative associated with $\mathcal{N}$. Let $j^k_0(M)=\,^k\pi^{-1}(x)$ be the fiber over $x$. Then the map
\begin{align}
\psi:j^k_0(x)\to \oplus^k \,T_xM,\quad (x,\dot{x},\ddot{x},...,x^{(k)})\mapsto (x,D_{\dot{x}}\dot{x},D^2_{\dot{x}}\dot{x},...,D^{(k)}_{\dot{x}}\dot{x})
\label{mappsi}
\end{align}
is a bijection: if one considers normal coordinates containing $n$, $\psi$ can be written as
\begin{align*}
(x,\dot{x},\ddot{x},...,x^{(k)})\mapsto (x,\dot{x},\ddot{x},...,x^{(k)}),
\end{align*}
 which shows clearly that it is injective and surjective. Therefore, one can define a vector space structure $(j^k_0(x),+,\cdot)$, defined by
 \begin{itemize}
 \item $+:j^k_0(x)\times j^k_0(x)\to \,j^k_0(x),\quad j_1+j_2:=\psi^{-1}(\psi(j_1)+\psi(j_2)),$

 \item $\cdot:j^k_0(x)\times R\to \quad j^k_0(x), \lambda\cdot j_1:=\,\psi^{-1}(\lambda\,\psi(\,j_1)).$
 \end{itemize}
 Let us consider the Whitney sum $\oplus^k_W \,TM$ of $k$ tangent bundles $TM$.
\begin{proposicion}
Given a symmetric connection $\mathcal{N}$ on $M$, the fiber bundle $J^k_0(M)$ is furnished with a vector bundle structure,
which is fibrewise induced from \eqref{mappsi} and such that the following isomorphism holds:
\begin{align}
(J^k_0(M),+,\cdot)\simeq\,\oplus^k_W \,TM.
\end{align}
\label{proposiciononpsi}
\end{proposicion}
In particular, one can use the Levi-Civita connection $D$ associated with the metric $\eta$ and its extension of the Levi-Civita connection to higher order jet bundles. However, such a choice is not necessary and there are other possibilities for a natural choice of $\mathcal{N}$, for instance, with connection associated with a given spray.

As an application of {\it Theorem} \ref{Thomstheorem} and {\it Proposition} \ref{proposiciononpsi}, we found the following isomorphism:
\begin{proposicion}
Let $H^*_{cv}(J^k_{0}(M))$ be as before and $H^*_{dR}(M)$ be the real
de Rham cohomology of $M$. Then the following isomorphism holds:
\begin{align}
H^*_{cv}(J^k_{0}(M))\,\simeq H^{*-kn}_{dR}(M).
\label{HJ1=HdR}
\end{align}
\label{isomorphismthom}
\end{proposicion}
{\bf Proof}. One can apply Thom's isomorphism theorem \ref{Thomisomorphism forE} to the vector bundle
$^k\pi:J^k_0(M)\to M$, obtaining by fiber averaging the following isomorphism:
\begin{align*}
\xymatrix{ H^*_{cv}(J^{k}_0(M))
 \ar[r]^{\langle \cdot\rangle} & H^{*-nk}_{dR}(M)},
\end{align*}
from which follows the result.\hfill$\Box$
\begin{corolario}
The cohomology $H^*_{cv}(J^k_{0}(M))$ is an homotopy invariant of $M$.
\end{corolario}
\begin{proposicion}
Let $(M,\eta)$ be a Lorentzian manifold and $H^p_{cv}(M,\mathcal{F}(J^k_{0b}(M)))$ as before.
Then any element of $[\alpha]\in H^p_{cv}(M,\mathcal{F}(J^k_{0b}(M)))$ is bounded.
\end{proposicion}
{\bf Proof}. Using the corresponding operator norm, one obtains the result. \hfill$\Box$

Let us recall the injective homomorphism $\varphi$ defined by \eqref{embeddingforms},
\begin{align*}
\varphi:\Lambda^p M\to \Lambda^p(M,\mathcal{F}(J^k_{0}(M))),\,\,\varphi({\alpha})_u(X_1,...,X_p)=\alpha_x(X_1,...,X_p).
\end{align*}
 Then the following result is
a consequence from the isomorphism (\ref{isomorphismthom}),
\begin{lema}Let us assume the hypothesis as in Thom's theorem. Then
each element of a class $[\hat{\alpha}]\in H^{p}_{cv}(J^k_0(M))$ admits a decomposition,
\begin{align}
^3\zeta (\hat{\alpha})=\,^k\zeta(\varphi(\langle \hat{\alpha}\rangle))\,+d_J A, \quad A\in \Lambda^{p-1}_{cv}(J^k_0(M)),\quad \hat{\alpha}\in\,[\hat{\alpha}].
\label{decomposition}
\end{align}
This decomposition is unique.
\label{descomposicion de alpha}
\end{lema}
{\bf Proof}. By the invariance of integration operation along the fiber $\langle\cdot\rangle$, one has that
\begin{align*}
H^{p+kn}_{cv}(J^k_0(M))\,\simeq H^p_{dR}(M).
\end{align*}
Note that for each cohomology class $[\hat{\alpha}]\in \, H^{p+kn}(J^k_0(M))$, one has that
\begin{align*}
\langle \hat{\alpha}\rangle=\,
\langle\varphi(\langle \hat{\alpha}\rangle)\rangle.
\end{align*}
 Therefore, by Thom's isomorphism theorem \ref{isomorphismthom}, it follows that
\begin{align*}
[\hat{\alpha}]=[\varphi(\langle\hat{\alpha}\rangle)]
\end{align*}
 and one obtains the relation (\ref{decomposition}).\hfill$\Box$

 We prove now an isomorphism which is of relevance for physical interpretation of the theory of generalized electromagnetic fields,
\begin{proposicion}
Let $J^k_{0b}(M)$ be the bounded $k$-jet bundle over $M$ with constants $(c_1=\infty,\,c_i<\infty,\,1<i\leq k)$.
Then there is the following isomorphism,
\begin{align}
\varphi^*_b:H^*_{dR}(M)\to H^*_{cv}(M,\mathcal{F}(J^k_{0b}(M))).
\label{isomorphismcohomology}
\end{align}
\label{isomorphism}
\end{proposicion}
{\bf Proof}.
Let us consider the homomorphism $\varphi_b:\Gamma \Lambda^p(M)\to \Gamma \Lambda^p_{cv}(M,\mathcal{F}(J^k_{0b}(M)))$, as a restriction of $\varphi$.
The homomorphism $\varphi_b$ commutes with $d_4$,
\begin{align}
\varphi_b(d\alpha)=d_4(\varphi_b(\alpha)).
\end{align}
 Therefore, it induces an homomorphism of cohomologies:
\begin{align*}
\varphi^*_b:H^*_{dR}(M)\to H^*_{cv}(M,\mathcal{F}(J^k_{0b}(M))).
\end{align*}
$\varphi^*_b$ is injective, which follows from the injectivity of $\varphi_b$. The isomorphism (\ref{isomorphismcohomology}) is also surjective.
Let $[\bar\alpha]\in H^p_{cv}(M,\mathcal{F}(J^k_{0b}(M)))$ be a class of cohomology
such that for each $\bar\alpha\in [\alpha]$, $d_4(\alpha)=0$ holds. Let us define the homomorphism
\begin{align}
\tilde{\varphi}_b  :H^*_{cv}(M,& \mathcal{F}(J^k_{0b}(M)))\to H^{*+kn}_{cv}(J^k_0(M)),\quad [\bar\alpha]\mapsto [\,^k\xi(\bar\alpha)\wedge dvol_V].
\label{anotherisomorphism}
\end{align}
Then $\langle \,^k\xi(\bar\alpha)\wedge dvol_V\rangle$ is in some cohomology class in $H^*_{dR}(M)$. Thus the image $\varphi^*_b(\langle \,^k\xi(\bar\alpha)\wedge dvol_V\rangle)$ is in the same cohomology class of $\bar\alpha$
and the result is proved.\hfill$\Box$
\begin{teorema}
The cohomology $H^*_{cv}(M,\mathcal{F}(J^k_{0b}(M)))$ is a homotopy invariant of $M$.
\end{teorema}
{\bf Proof}. Since it is related with the usual de Rham cohomology $H^*_{dR}(M)$, it is an homotopy invariant of $M$.\hfill$\Box$

As a consequence of {\it Proposition} (\ref{isomorphism}) and the homotopy invariance of the cohomology, we have that
\begin{corolario}
Let $\pi_B:B\to M$ be a sub-bundle of $\pi_V:\mathcal{V}\to M$ with fibers $\pi^{-1}_B(x)$ homotopic to $\pi^{-1}_V(x)$. Then
\begin{align*}
H^*_{cv}(B)\simeq H^{*}_{cv}(M,\mathcal{F}(J^k_0(M))).
\end{align*}
\end{corolario}

As a relevant example we have the following,
\begin{proposicion}
One has the following isomorphism
\begin{align}
H^*_{cv}(M,\mathcal{F}(\tilde{J}^3_{0b}(M)))\simeq\,H^*_{dR}(M).
\label{isomorphismtildej3j2}
\end{align}
\label{proposicionsobretildej}
\end{proposicion}
{\bf Proof}. Note that the fiber $\tilde{j}^3_{0b}(x)$ has one extra dimension more than $j^2_{0b}(M)$, corresponding to the direction of the globally defined coordinate
\begin{align}
\beta^{-1}:\,\tilde{j}^3_{0b}(x)\to R,\quad ^3x\mapsto \beta_2(\,^kx):=|\eta(D^2_{\dot{x}}\dot{x},D_{\dot{x}}\dot{x})|.
\end{align}
Then there is the following isomorphism
\begin{align}
\psi:\Lambda^*J^2_{0b}(M)\to \,\Lambda^{*+1}\tilde{J}^3_{0b}(M),\quad \hat{\alpha}\mapsto \hat{\alpha}\wedge \beta\,d(\beta^{-1}).
\end{align}
This induces the following isomorphism,
\begin{align}
\psi^*:H^{*+2n}J^2_{0b}(M)\to \,\Lambda^{*+2n+1}\tilde{J}^3_{0b}(M),\quad [\hat{\alpha}\wedge\,dvol_V]\mapsto [\hat{\alpha}\wedge dvol_{V}\wedge\,\beta d(\beta^{-1})].
\end{align}
From which follows the following isomorphism (by homotopy and fiber integration),
\begin{align}
\vartheta:H^{*+2n+1}_{cv}\tilde{J}^3_{0b}(M)\to  H^*_{cv}(M,\mathcal{F}(\tilde{J}^3_{0b}(M))),\quad [\hat{\alpha}\wedge dvol_{V}\wedge\,\beta d(\beta^{-1})]\mapsto [\beta\,\bar{\alpha}].
\end{align}
Thus we have the following isomorphisms
\begin{align*}
H^{*+2n+1}_{cv}\tilde{J}^3_{0b}(M)\simeq H^{*+2n}_{cv}(J^2_{0b})\simeq H^*_{cv}(J^1_0(M))\simeq H^*_{dR}(M).
\end{align*}\hfill$\Box$
\begin{comentario} We have the following remarks:
\begin{itemize}
\item The closed vertical volume form in $\tilde{J}^3_{0b}(M)$ is $dvol_V\wedge\,\beta d(\beta^{-1}).$

\item There is a natural map (not depending on coordinates)
\begin{align}
\Lambda^p(M,\mathcal{F}(\tilde{J}^3_{0b}(M)))\to \Lambda^p(M,\mathcal{F}({J}^3_{0}(M))),\quad \hat{\alpha}\wedge\,\beta d(\beta^{-1})\mapsto \beta\,\hat{\alpha}
\end{align}
which is injective.
 \end{itemize}
\end{comentario}

\begin{corolario}
There is the following isomorphism,
\begin{align}
 H^*_{cv}(M,\mathcal{F}(J^1_{0}(M)))\simeq\, H^{*}_{dR}(M).
\end{align}
\label{isomorfismoj0bdr}
\end{corolario}
We remark that the two ingredients, vertical compact support and bounded jet bundle manifold are of relevance to our application to electrodynamics: the compact vertical domain is relevant to have consistence with bounded jet bundle base manifold and to avoid infinite kinetic world-lines; bounded cohomology is useful to avoid run-away solutions.

\subsection{Integration theory of generalized forms}
\begin{definicion}
The integral of a generalized form $\bar{\alpha}\in \, \Gamma\Lambda^p (M,\mathcal{F}(J^k_0(M)))$ on the $p$-dimensional submanifold $M_p$ of $M$ is
\begin{align}
\int_{M_p}\,\bar{\alpha}:=\,\int_{M_p}\,\langle\,\hat{\alpha}\rangle.
\label{integralofgeneralizedform}
\end{align}
\label{definiciondeformaintegral}
\end{definicion}
We note the commutativity of the integration operations,
\begin{align}
\int_{M_p}\,\langle\,\hat{\alpha}\rangle=\,\langle\int_{M_p}\,\,\hat{\alpha}\rangle.
\end{align}
 There is a direct formula for the integral \eqref{integralofgeneralizedform} in terms of integral of forms $\hat{\alpha}$,
\begin{align}
\int_{M_p}\,\langle\hat{\alpha}\rangle =\,\int_{^k\pi^{-1}(M_p)}\,\hat{\alpha}\wedge \,dvol_v,
\end{align}
where $dvol_v$ is a vertical volume form. At this point, it is not required that $dvol_v$ is closed, but we will require later such property to prove the corresponding Stokes' theorem for generalized forms.

From the definition \ref{definiciondeformaintegral}, the following is direct,
\begin{proposicion}
For the integral operation \eqref{integralofgeneralizedform} the following properties hold:
\begin{itemize}
\item It reduces to the standard definition for the case of standard differential forms.

\item It is direct that the integral of generalized forms is linear,
\begin{align}
\int_{M_p}\,\big(\bar{\alpha}+\,\lambda \bar{\beta}\big)=\,\int_{M_p}\,\bar{\alpha}+\,\lambda\,\int_{M_p}\,\bar{\alpha},
\end{align}
for any $\hat{\alpha}, \,\hat{\beta}\in \,\Gamma\,\Lambda(M,\mathcal{F}(J^k_0(M)))$ and $\lambda\in \,R$.
\end{itemize}
 \end{proposicion}

\begin{lema}
If $\zeta{\bar{\alpha}}=\,\hat{\alpha}$, then
\begin{align}
\int_{M_p}\,d_4\bar{\alpha}=\,\int_{M_p}\,\langle d_J\hat{\alpha}\rangle.
\end{align}
\end{lema}
{\bf Proof}. We make use of {\it Lemma} \ref{propiedadesdehathdJ},
\begin{align}
\int_{M_p}\,d_4\bar{\alpha}=\,\int_{M_p}\,\langle\hat{h}_k d_J\hat{\alpha}\rangle=\,\int_{M_p}\,\langle d_J\hat{\alpha}\rangle.
\end{align}\hfill$\Box$

\subsection*{Generalization of Stokes' Theorem}

There is a version of Stokes' theorem
 for generalized forms and an invariance under diffeomorphisms of differential forms,
\begin{proposicion}
Let $\bar{\alpha}\in \, \Gamma\Lambda^p (M,J^k_0(M))$ such that $d_4(\bar{\alpha})=0$ and consider a
$p$-dimensional submanifold $M_p$ as before. Let also $d_J(dvol_v=0)$. Then the following formula holds,
\begin{align}
\int_{M_p}\,d_4\bar{\alpha}=\,\int_{\partial M_p}\,\bar{\alpha}.
\label{stokestheoremforintegralforms}
\end{align}
\end{proposicion}
{\bf Proof}. From the definition,
\begin{align*}
\int_{M_p}\,d_4\bar{\alpha}=\,\langle\Big(\int_{M_p}\,d_J\hat{\alpha}\Big)\rangle=\,\int_{^k\pi^{-1}(M_p)}\,(d_J\hat{\alpha})\wedge \,dvol_v=*.
\end{align*}
Then note that the exterior derivative of $dvol_J$ is zero by hypothesis. Thus we have that $d_J(\hat{\alpha}\wedge\,dvol_v)=\,d_J\hat{\alpha}\wedge\,dvol_v$. Therefore,
\begin{align*}
*=\,\int_{^k\pi^{-1}(M_p)}\,(d_J\hat{\alpha}\wedge \,dvol_v).
\end{align*}
Then one can use Stokes' theorem, showing that
\begin{align*}
\int_{^k\pi^{-1}(M_p)}\,d_J(\hat{\alpha}\wedge \,dvol_v)=\,\int_{\,^k\pi^{-1}(\partial(M_p))}\,\hat{\alpha}\wedge \,dvol_v.
\end{align*}
Taking into account the relation
\begin{align*}
 \partial(\,^k\pi^{-1}(M_p))=\,^k\pi^{-1}(\partial(M_p))
 \end{align*}
  and again the definition of the integral \ref{definiciondeformaintegral}, we have that
\begin{align*}
\int_{\,^k\pi^{-1}(\partial(M_p))}\,\hat{\alpha}\wedge \,dvol_v=\int_{\partial M_p}\,\langle\hat{\alpha}\rangle=\,\,\int_{\partial M_p}\,\bar{\alpha}.
\end{align*}
\hfill$\Box$
\subsection*{Invariance under diffeomorphisms of the integral}

Let $f:\tilde{M}\to M$ be a differential function. One can define the pull-back of a generalized form,
\begin{align}
f^*\bar{\alpha}(X_1,...,X_p):=\bar{\alpha}(f_* (X_1),...,f_*(X_p)).
\label{pullbackgeneralizedforms}
\end{align}
One has the following relations,
\begin{align*}
f^*\langle\hat{\alpha}\rangle(X_1,...,X_p)& =\,\langle\hat{\alpha}\rangle(f_* (X_1),...,f_*(X_p))\rangle\\
& =\,\,\langle \hat{\alpha}(f_* (X_1),...,f_*(X_p))\rangle\\
& =\,\langle f^*\hat{\alpha}(X_1,...,X_p)\rangle\\
& =\,\langle f^*\hat{\alpha}\rangle(X_1,...,X_p).
\end{align*}
Then one obtains the following invariance under diffeomorphism,
\begin{proposicion}
Let $f:\tilde{M}_p\to M_p$ be a diffeomorphism between $p$-dimensional manifolds. Then
\begin{align}
\int_{\tilde{M}_p}\,f^*\bar{\alpha}=\,\int_{M_p}\,\bar{\alpha}.
\end{align}
\label{invarianceoftheintegrak}
\end{proposicion}
{\bf Proof}. Using the above computation, one has
\begin{align*}
\int_{\tilde{M}_p}\,f^*\bar{\alpha}=\,\int_{\tilde{M}_p}\,f^*\langle\hat{\alpha}\rangle=
\,\int_{{M}_p}\langle\hat{\alpha}\rangle =\,\int_{{M}_p}\,\bar{\alpha}.
\end{align*}\hfill $\Box$

\section{Elements of maximal acceleration spacetimes}

The second principal assumption adopted in this work is that the $n$-acceleration of a point charged
particle is {\it bounded} (in a covariant way, that is, independent of local coordinate system).
 Let us briefly summarize the idea of maximal acceleration, prior to see where our point of view departs from previous ones.
 The original idea of maximal acceleration starts with E. Caianiello's work (see for instance the review \cite{Caianiello} and references therein and also the works
 \cite{Brandt1989} and \cite{Toller} for original developments of the idea of maximal acceleration)
in the contest of a geometrization of quantum mechanical systems. Thus, in Cainiello's theory, uncertainty in the observables is related with the curvature of a Sasaki-type metric in the relevant phase space. It was as a consequence of this and that speed of light in vacuum is an upper limit for causal interactions that one arrives to maximal proper acceleration.

Caianiello's theory was not general covariance of the theory: it was formulated on flat spacetime, using coordinates and for the proper acceleration.
 A covariant formalism for geometries of
maximal acceleration was developed in \cite{GallegoTorrome} and explores instead the possibility of {\it bounded covariant $n$-acceleration}.
 Although it was motivated by the non-covariance problem  of Caianiello's
 quantum geometry \cite{Caianiello}, the
 formalism developed in \cite{Caianiello} was independent of the details of the mechanism generating the bound in
 the covariant  acceleration and also from the particular value that the maximal acceleration $A_{max}$ can take. Thus,
 the framework developed in \cite{GallegoTorrome} can be used in a more general context than Caianiello's quantum geometry. That the maximal proper acceleration is bounded is a direct consequence from the fact that $n$-acceleration is bounded in our formalism.

 However, \cite{GallegoTorrome} presupposed the existence of the Lorentzian metric $\eta$,
 from where one constructs  the metric of maximal acceleration.
  In the theory presented in this {\it section} we revers the situation. We first introduce a generalized metric
  $\bar{g}\in\,\Gamma\, T^{(0,2)}(M\mathcal{F}(J^2(M)))$ as the natural object associated with measurements of the proper time of physical
  clocks at rest with the generic point particle $x:I\to M$. Then the {\it Lorentzian proper time} along a timelike curve
   is associated with the proper time of the average of the generalized metric $\bar{g}$. Thus, the Lorentzian spacetime metric emerges as averaged description of the geometry of maximal acceleration.

A main distinction between our approach and other theories and models of maximal acceleration
is that in our theory the kinematical formalism contains a maximal $n$-acceleration
respect to a given Lorentzian metric from the beginning. We will not attempt to provide the mechanism for the bound of the $n$-acceleration of
 point charged particles. Such mechanism should be based on a deepest description of spacetime as a discrete, quantum spacetime and will be developed elsewhere. Instead, we provide an heuristic argument for the existence of a maximal acceleration, based on maximal speed and minimal {\it characteristic length} Moreover, we speak of maximal bound $n$-acceleration, in contrast with Caianiello and others approach, that consider bound of the proper acceleration. Our approach to maximal acceleration, has the benefit that it is covariant respect to the Lorentz group and by the introduction of a connection, can be made general covariant, as reference \cite{GallegoTorrome} showed.

The notion of maximal acceleration appears in other theories in Physics.
In string theory, it appears as a consequence of the formation of {\it Jean's} instability when the strings
 reach a critical temperature, that makes the strings disconnected \cite{ParentaniPotting, BowickGiddins}. Very recently, maximal proper
 acceleration was found as a natural consequence in covariant loop gravity \cite{RovelliVidotto}.
  A dramatic consequence of maximal acceleration for the gravitational theory and the equivalence principle is that
  the corresponding theories in such geometries should be free of singularities, a point first noted by Caianiello.

\subsection*{Maximal acceleration in classical electrodynamics}
  In classical electrodynamics,
there are several scenarios were maximal acceleration appears:
\begin{itemize}
\item The  Abraham and Lorentz's electron models
are theoretically valid under the assumption that acceleration have a value less than a threshold value (see reference
\cite{Spohn2} for a discussion of those models), in order to preserve causality.

\item  In the extended model of the electron proposed by P. Caldirola \cite{Caldirola1956}
 it appears a maximal acceleration, when a maximal speed of interaction  and the hypothesis of the minimal unit of time
 {\it chronon} are used in the definition of acceleration \cite{Caldirola2}.

\item The existence of a maximal value for the electric field in Born-Infeld electrodynamics
\cite{BornInfeld} suggests that the four dimensional force on a point electron is also bounded
in such theory. Thus, it is natural to conjecture the existence of a bounded proper acceleration in such a case.
\end{itemize}

These three examples bring to light that the hypothesis of maximal acceleration could be of relevance for the solution of some of the important problems in the foundations of classical electrodynamics.

\subsection{An heuristic argument for maximal acceleration}

There is an heuristic argument for the existence of a maximal acceleration
based on the assumption that there is a minimal length $L_{min} $ and
a maximal speed. This argument, first used in the context of classical electrodynamics by Caldirola \cite{Caldirola2}, is here expressed in complete generality. The minimal length is assumed to be scale of the domain in the spacetime that affects the individual particle in changing the dynamical state.
This idea is not necessarily related with a quantification of spacetime,
but requires a notion of extended local domain where cause-effect relations are originated.
Therefore, the maximal acceleration could be relational, depending on the
physical system. This is in contrast with universal maximal acceleration.
However, we will require that the maximal acceleration is
very large compare with the acceleration of the probe particle.
In this way, our perturbative scheme will be perfectly applicable.

By adopting the above hypotheses, {\it the effect} on a particle done by its surrounding is bounded by a maximal work
\begin{align*}
L_{min} m \,a \sim \delta {m}\,v_{max}^2 ,
\end{align*}
where $a$ is the value of the acceleration in the direction of the total exterior effort is done.
 Then one associates this value to the work over any fundamental degree of freedom evolving in $M$, caused by rest of the system.
Since  the speed must be bounded, $v_{max}\leq
c$. Also, the maximal work produced by the system on a point particle is
 $\delta {m}=\, -m$. Thus, there is a bound for the value of the acceleration,
\begin{align}
a_{max}\,\simeq \frac{c^2}{L_{min}}.
\label{maximalacceleration}
\end{align}
Thus, the existence of maximal propagation speed and minimal length are of
fundamental relevance in the argument for maximal acceleration.
This will be of relevance when we investigate the possibility of superluminical motion in this {\it section}.

\subsection{Principle of maximal $n$-acceleration and the clock hypothesis in classical mechanics}

The notion of physical clock should be related with the general form of the principles of classical dynamics. The first dynamical principle corresponds to the {\it principle  of inertia} and the notion of {\it inertial coordinate system},
\begin{definicion}
An inertial coordinate system $(U,R^n)$ is such that the world-line of any free particle is described by a parameterized  straight line of $R^n$.
\label{principle of inertia}
\end{definicion}
Thus, if the world-line describing a point particle is not a straight line in some given coordinate system, it must be because the point particle is not moving free or the coordinate system is not inertial or both. Note that the geodesic motion is not free motion, even if there is a coordinate system (normal coordinate system) where the geodesic is described by straight line in $R^n$.

Even if strictly speaking, there are not perfect inertial coordinate system in real experiments, {\it Definition} \ref{principle of inertia} is not empty of content, as it serves as a natural foundation for the second law of dynamics in terms of differential equations. The simplest possibility of a equation of motion for  point particles, in concordance with experience and the principle of inertia, is that such differential equations are of second order respect to the time parameter in the inertial coordinate system. Moreover, one can consider approximate inertial systems, that relate with experimental settings. This shows that the notion of inertial coordinate system is useful.

There is some arbitrariness in the choice of the time parameter in the differential equation for a given point particle. However, there are some parameters that appear more natural that others, in the context of the discussion before.
\begin{definicion}
Given a world-line $\tilde{x}:I\to M$ corresponding to a physical point particle, a {\it physical time parameter} $s$ is such that the curve $\tilde{x}:\tilde{I}\to M$ satisfies a second order differential equation respect to $s$.
\label{definitionofphysicalclock}
\end{definicion}
If the world-line $x:I\to M$ satisfies a second order differential equation respect to a physical time parameter, a natural definition of {\it co-moving physical clock} is the following:
\begin{definicion}
A {\it co-moving physical clock} associated with the world-line $\tilde{x}:I\to M$ of a point particle is a map $\tau:\tilde{x}(\tilde{I})\to R,\tau\mapsto x(\tau)$ such that the following diagram
\begin{align*}
\xymatrix{ &
{ J^k_0(M)} \ar[d]^{^2\pi}  \ar[r]^{\tau} & R  \ar[dl]^{^2x}\\
\tilde{ I} \ar[ur]^{^2\tilde{x}}  \ar[r]^{\tilde{x}} & { M} & I \ar[u]^i \ar[l]^{x} &.}
\end{align*}
commutes for any physical time parameterization $\tilde{x}:\tilde{I}\to M$ of the un-parameterized curve $x(I)$.
\label{definitionofphysicalco-movingclock}
\end{definicion}
 We need to prove the existence of co-moving physical clock parameters. In Special and General Relativity, the existence of a co-moving physical time parameter is guaranteed by {\it the clock hypotheses}, {\it the Principle of Relativity} and the {\it principle of constancy of the speed of light in vacuum}. As a result, the parameter $\tau$ can be  chosen to be the proper-time of a Lorentzian metric $\eta$. Such proper-time of $\eta$ corresponds to the time of a co-moving physical clock as in {\it Definition} \ref{definitionofphysicalclock} that is independent of the acceleration respect to an inertial coordinates system.

 \subsection*{The clock hypothesis}

 For world-lines whose acceleration vector is different than zero, it was found useful to make an additional assumption about the rate of clocks associated with the world-line. The point was, to consider the clocks as un-affected by the acceleration of the curve, in such a way that they work exactly as a relativistic clock will do (in the corresponding instantaneous inertial frame).
 The {clock hypothesis} can be formulated as follows (see \cite{Einstein1922}, p. 64 and \cite{Rindler}, p. 65):

{\it There exist ideal clocks, that is, clock that are completely
unaffected by acceleration; that is, as one whose instantaneous rate depends only
on its instantaneous speed in accordance with the time dilatation formula of Special Relativity. Thus, one can adopt such clocks as the co-moving proper clocks.}

This hypothesis is of fundamental relevance for the foundation of gravity theory in terms of a metric structure and in Special Relativity to define instantaneously inertial systems and its relation with inertial coordinate systems.
  However, the clock hypothesis it does not necessarily hold for arbitrary co-moving physical clocks
  Indeed, for clocks that measure time according to a generalized metric, the clock hypothesis does not necessarily holds, since the rate of the proper clocks could depend on the acceleration. Indeed, we will see later that there are other options that the proper time associated with a Lorentzian metric that are compatible with {\it Definition} \ref{definitionofphysicalco-movingclock}.

\begin{definicion}
Given a physical world-line $\tilde{x}:\tilde{I}\to M$, parameterized by a physical parameter $t$, an instantaneous inertial coordinate system associated to the world-line $\tilde{x}$ is a local coordinate system such that, respect to a given inertial system $(U,R^n)$, it moves with constant speed $\frac{d\tilde{x}^\mu}{dt}$.
\end{definicion}

\subsection*{Principle of maximal $n$-acceleration}
 We can formulate the principle in the following way,

 {\it For each physical world-line, there is associated a physical co-moving clock such that for the time parameter $\tau$ of such clock, the $n$-acceleration measured is bounded by a constant value $A_{max}$ that does not depend on the world-line in the sense that}
\begin{align}
\eta(\ddot{x},\ddot{x})<A_{max}.
\label{boundproperacceleration}
\end{align}
{\it for the Lorentzian metric} $\eta$.
Therefore, we reverse the previous theories of maximal acceleration and introduce the principle of maximal $n$-acceleration from the beginning. The principle, stated in such a way, will convey the modifications of the field theories as well as the classical equation of motion of point particles.

In the next {\it section} we will prove that it is possible to construct a generalized metrics $g$ from a Lorentzian metric $\eta$, such that principle of maximal $n$-acceleration holds. On the other hand, the clock hypothesis does not hold for geometries of maximal acceleration. This is a form of converse result that the one in  \cite{FriedmanGofman}, where the starting point is the non-validity of the clock hypothesis, they concluded the existence of a maximal (universal) acceleration \cite{FriedmanGofman}.

\subsection{General covariant formulation of the metric of maximal $n$-acceleration}

Since we are assuming the existence of the Lorentzian metric $\eta$, we can construct its Levi-Civita connection $D$ and its derivative operator along $x:I\to M$. Then the four covariant acceleration vector $D_{\dot{x}}\dot{x}\in \,T_{x(\tau)}M$
is bounded by the Lorentzian metric $\eta$ (note that the
covariant acceleration is spatial like vector). This is defined on the bundle $TM\setminus NC$, where
$\pi_{NC}: NC\to M$ is the null-cone bundle determined by $\eta$ and
\begin{align*}
NC:=\,\bigsqcup_{x\in M}\,NC_x,\quad NC_x:=\,\big\{y\in\,T_xM\,s.t.\,\, g(y,y)=0\,\big\}.
\end{align*}

 Let $(M,\eta)$ be a Lorentzian $n$-dimensional spacetime.
Then there is defined a Sasaki-type metric on the bundle $TM\setminus NC$,
\begin{align}
g_S= \,\eta_{\mu \nu} dx^\mu \otimes dx^\nu + \frac{1}{A^2_{max}}\eta_{\mu \nu}\Big( {\delta 
y^{\mu}}\otimes {\delta y^{\nu}}\Big).
\label{sasakitypemetric}
\end{align}

\begin{teorema}
Let $x:I\to { M}$ be a curve such that
\begin{enumerate}
\item If the tangent vector at the point $^1x\in \,J^1_0(M)$ is $T=(\dot{x},\ddot{x})$ and

\item It holds that $\eta(\dot{x},\dot{x})\neq 0$.
\end{enumerate}
Then there is a non-degenerate, symmetric form $g$ obtained from the embedding
{\it isometric embedding} $e:x(I) \hookrightarrow  TM$ from $(TM,g_S)$ such that acting on the tangent vector $\dot{x}$ the value is
\begin{align}
g_{\mu \nu}(x(\tau)) =\Big(1+ \frac{  \eta(D_{\dot{x}}\dot{x}(\tau)
D_{\dot{x}}\dot{x}(\tau))}{A^2 _{max}\,\eta(\dot{x},\dot{x})}\Big)\, \eta_{\mu \nu}.
\label{maximalaccelerationmetric}
\end{align}
\label{teoremasobremaximaacceleration}
\end{teorema}
{\bf Proof.} The tangent vector at the point $(x(\tau),\dot{x}(\tau))=\,^1x(\tau)\in TM$
is $(\dot{x},\ddot{x})\in \,TTM$. The metric $g_S$ acting on the vector field $T=(\dot{x},\ddot{x})\in\,T_{(x(\tau),\dot{x}(\tau))}N$ has the value
\begin{align*}
g_S(T,T)& =\,\Big( \eta_{\mu \nu}\, dx^\mu \otimes dx^\nu + \frac{1}{A^2_{max}}\eta_{\mu \nu}\Big( {\delta 
y^{\mu}}\otimes {\delta y^{\nu}}\Big)\Big)\,\big(T,T\big)\\
& =\,\Big( \eta_{\mu \nu}\, \dot{x}^\mu \dot{x}^\nu + \frac{1}{A^2_{max}}
\eta_{\mu \nu}\big(\ddot{x}^{\mu}\,-N^{\mu}\,_{\rho}(x,\dot{x})\,\dot{x}^{\rho}\big) \,
\big(\ddot{x}^{\nu}\,-N^{\nu}\,_{\lambda}(x,\dot{x})\,\dot{x}^{\lambda}\big)\Big)=*.
\end{align*}
By the second hypothesis, one obtains that the following expression holds:
\begin{align*}
* =\,\Big( 1 + \frac{1}{A^2_{max}\,\eta(\dot{x},\dot{x})}
\eta_{\mu\nu}\,D_{\dot{x}}\dot{x}^{\mu}\,D_{\dot{x}}\,\dot{x}^{\nu}\Big)\eta(\dot{x},\dot{x}),
\end{align*}
that coincides with the value $(g_{\mu\nu}\,d{x}^{\mu}\otimes d{x}^{\nu})(\dot{x},\dot{x})$,
with the components $g_{\mu\nu}$ given by the formula \eqref{maximalaccelerationmetric}.
\hfill$\Box$

The bilinear form
\begin{align}
g(\,^2x)=\,g_{\mu\nu}(\,^2x)\,dx^\mu\otimes dx^\nu
\end{align}
 with components given by the covariant formula \eqref{maximalaccelerationmetric} is the {\it metric of
maximal acceleration} in general coordinates. It determines the proper
time along the curve $x:I\to M$ and also provides a generalization of the notion of
angle.

 \begin{corolario} Let $x:I to M$ such that:
 \begin{itemize}
 \item It holds that
 $g(\dot{x},\dot{x})<0$, $\eta(\dot x,\dot x)<0$,

\item The covariant condition
\begin{align}
\eta(D_{\dot{x}}\,\dot{x},\,D_{\dot{x}}\,\dot{x})\,\geq 0.
\label{spacelikeaccelerations}
\end{align}
holds.
Then
 one has the natural bound
 \begin{align}
  0 \leq \eta(D_{\dot{x}}\dot{x},D_{\dot{x}}\dot{x})<\,A^2_{max}.
  \label{boundedconditionforacceleration}
 \end{align}
 \end{itemize}
 \end{corolario}
Because this property, the bilinear form \eqref{maximalaccelerationmetric} is a metric
 of maximal acceleration, since the covariant $n$-acceleration it turns out to be bounded by the maximal value $A_{max}$.
\begin{definicion}
A curve of maximal acceleration is a curve $x:I\to M$ such that
\begin{align}
\eta(D_{\dot{x}}\dot{x},D_{\dot{x}}\dot{x})=A^2_{max}.
\label{curveofmaximalacceleration}
\end{align}
\end{definicion}
\begin{corolario}
For a curve of maximal acceleration $x:I\to M$, one has the relation
\begin{align}
 g(\dot{x},\dot{x})=1+\eta(\dot{x},\dot{x}).
\end{align}
\label{lightlikecurvesofmaximal acceleration}
\end{corolario}
This result indicates that for maximal acceleration curves, the proper parameters associated to $\eta$ and $g$ differ considerably.

If the covariant derivative $D$ is the induced connection on $TM\setminus NC$ induced by
 the Levi-Civita connection of the Minkowski metric $\eta$ as it was defined in {\it section 2},
there is defined globally a coordinates system where the connection coefficients $\gamma^{\mu}\,_{\nu\rho}$
are zero. In such coordinate system $N^{\mu}\,_{\rho}=0$ holds and therefore $D_{\dot{x}}\,\dot{x}=\ddot{x}$.
In such coordinate system it also holds $\gamma^{\mu}_{\nu\rho}=0$.
Therefore, the metric coefficients(\ref{maximalaccelerationmetric}) can be written as
\begin{align}
g_{\mu\nu}(x(\tau)):=\Big(1+ \frac{  \eta_{\sigma\lambda}\,\ddot{x}^{\sigma}
(\tau)\ddot{x}^{\lambda}(\tau)}{A^2 _{max}\,\eta(\dot{x},\dot{x})}\Big)\, \eta_{\mu \nu}\,dx^{\mu}\otimes dx^{\nu},
\label{gtau}
\end{align}
that defines an element $g\in\,\Gamma T^{(0,2)}(M,\mathcal{F}J^2_0(M))$ in the following way: given
a curve $^2x:I\to J^2_0(M)$, the value of $g$ along the curve $s\mapsto x(s)$ is $g(\,^2x(\tau))=g(\tau)$.

\begin{comentario}
The appearance of the metric $\eta$ avoids  a fully general invariant theory of geometry of maximal acceleration. It should be much more natural to obtain $\eta$ from the basic fundamental object that is a generalized metric $\bar{g}$. However, we hope that a future form of the theory could provide a natural origin to the metric $\eta$ from fundamental principles.
\end{comentario}

\subsection{Perturbation scheme}

 Let us denote by $\tau$ the proper-time parameter along a given curve $x:I\to M$ respect to $g$.
 The acceleration square function is defined by the expression
\begin{align}
 a^2(\tau):=\,\eta_{\mu\rho}\,\ddot{x}^{\mu}\,\ddot{x}^{\rho}.
 \label{accelerationsquarefunction}
 \end{align}
 In our considerations the curves are {\it far} from the maximal acceleration condition,
$a^2(\tau)\ll\, A^2_{max}$. To make precise meaning to such statement we need of a perturbation theory. Let us consider the difference
\begin{align*}
\delta(\tau):=\eta(\dot{x},\dot{x})-g(\dot{x},\dot{x}).
\end{align*}
Then the following approximations hold,
 \begin{align*}
 \Big(1+ \frac{  \eta_{\sigma\lambda}\,   \ddot{x}^{\sigma}(\tau)
\ddot{x}^{\lambda}(\tau)}{A^2 _{max}\eta(\dot{x},\dot{x})}\Big) \eta(\dot{x},\dot{x})& = \,\Big(1+ \frac{  \eta_{\sigma\lambda}\,   \ddot{x}^{\sigma}(\tau)
\ddot{x}^{\lambda}(\tau)}{A^2 _{max}(g(\dot{x},\dot{x})+\delta (\dot{x},\dot{x}))}\Big)\eta(\dot{x},\dot{x})\\
& \simeq\,\Big(1- \frac{  \eta_{\sigma\lambda}\,   \ddot{x}^{\sigma}(\tau)
\ddot{x}^{\lambda}(\tau)}{A^2 _{max}}\big( 1-\delta (\dot{x},\dot{x})\big)\Big)\eta(\dot{x},\dot{x})\\
& \simeq \Big(1- \frac{  \eta_{\sigma\lambda}\,   \ddot{x}^{\sigma}(\tau)
\ddot{x}^{\lambda}(\tau)}{A^2 _{max}}\Big)\eta(\dot{x},\dot{x})+\,\frac{  \eta (\ddot{x}(\tau),
\ddot{x}(\tau))}{A^2 _{max}}\delta\,\eta(\dot{x},\dot{x}).
 \end{align*}
Since the form $\delta$ is {\it small} on curves far from the maximal acceleration curves, we neglect the second term
\begin{align}
\Big|\Big(1-\frac{  \eta_{\sigma\lambda}\,   \ddot{x}^{\sigma}(\tau)
\ddot{x}^{\lambda}(\tau)}{A^2 _{max}}\Big)\eta(\dot{x},\dot{x})\Big|>> \,\Big|\frac{  \eta (\ddot{x}(\tau)
\ddot{x})}{A^2 _{max}}\delta\,\eta(\dot{x},\dot{x})\Big|.
\end{align}
We will keep the above approximation $\delta\simeq 0$ when $a^2<<A^2_{max}$ in the
calculations performed in this work. This is because the function $\delta(\dot{x},\dot{x})$
adds a higher order term in the perturbative expressions that we will consider. Thus for
instance, the expression of the generalized metric without such approximation will be
\eqref{maximalaccelerationmetric}; with the approximation, the metric is instead
\begin{align*}
g_{\mu\nu}=\, (1-\,\frac{\eta(D_{\dot{ x}}\dot{ x},D_{\dot{ x}} \dot{ x})}{A^2_{max }})\eta(\dot{ x},\dot{ x})\,+\textrm{higher order terms}.
\end{align*}

 In a normal coordinate system for $\eta$, the function
${\epsilon}$ is defined by the relation
\begin{align}
{\epsilon}(\tau):=\,\big(\frac{ \eta_{\sigma\lambda}\,\ddot{x}^{\sigma}(\tau)
\ddot{x}^{\lambda}(\tau)}{A^2_{max}}\big).
\label{definiciondeepsilon}
\end{align}
The covariant definition of the function $\epsilon$ requires of a non-linear connection,
\begin{align}
{\epsilon}(\tau):=\,\big(\frac{\eta(D_{\dot{x}}\dot{x},D_{\dot{x}}\dot{x})}{A^2_{max}}\big).
\label{covariantdefiniciondeepsilon}
\end{align}
 For timelike trajectories and from the relation \eqref{gtau}, the relation between $g$ and $\eta$
\begin{align}
g(\tau)=\,(1-\epsilon(\tau))\eta.
\label{relationgeta}
\end{align}
It follows that the relation between the proper parameter of $\eta$ and $g$ is
 \begin{align}
 d s=\,(1-\dot{\epsilon})^{-1}\,d \tau.
 \label{relationtaus}
 \end{align}
 Therefore, given a timelike curve respect to $\tau$ the relation between $s$ and $\tau$ is not
  strictly speaking a re-parameterization $r:R\to R$. This is related with the fact that $g$ measures the physical proper times, while $\eta$ appears as a derived, although convenient object.

The function $\epsilon(\tau)$ determines a bookkeeping parameter
${\epsilon}_0$ by
where
\begin{align}
{\epsilon}(\tau)={\epsilon}_0\,h(\tau),\quad \epsilon_0=\,\max\{\,\epsilon(\tau),\,\tau\in I\}.
\end{align}
For compact curves,  the bookkeeping parameter always exists. However, we will need to
bound the value of higher order derivatives in order to keep such parameter bound for
non compact curves. A generator set for asymptotic expansions is $\{\epsilon^l_0,\, l=-\infty,...,-1,0,1,...\}$.
Now we can make sense of the statement that we will consider curves that are far from
 the curves of maximal acceleration. What that means is that the effects of order
  $\epsilon^2_0$ or higher in any analytical function on $\epsilon_0$ are negligible compared with first order.
Then the monomials in powers of the derivatives of $\epsilon$ define a basis for asymptotic expansions.

Also, for curves of maximal acceleration the generalized metric \eqref{relationgeta} is degenerated,
\begin{align*}
g(Z,Z)|_{a=A_{max}}=\,(1-\epsilon)|_{a=A_{max}}\eta(Z,Z)=0,
\end{align*}
for all vector field along $^2x:I\to M$. Therefore, more of our considerations will not be applicable to such curves. Indeed,
we will assume that all the derivatives $(\epsilon, \dot{\epsilon}$, $\ddot{\epsilon}...,)$ are small.

The generalized metric $g$ defines different kinematical relations than $\eta$.
The parametrization of the world-line curves are such that  $g(\dot{x},\dot{x})=-1$, implying that $\eta(\dot{x},\dot{x})\neq -1$ in general. Indeed, we have the following
\begin{proposicion}
In a maximal acceleration geometry space $(M,g)$ the following kinematical conditions hold:
\begin{align}
& g(\dot{x},\dot{x})= -1,\label{covariantkineticconstrain1}\\
& g(\dot{x},\ddot{x})=\,\frac{\dot{\epsilon}}{2}\eta(\dot{x}\dot{x})\label{covariantkineticconstrain2}+\,\textrm{higher order terms}\\
&g(\dddot{x},\dot{x})+\,g(\ddot{x},\ddot{x})=
\,\frac{d}{d\tau}\Big(\frac{\dot{\epsilon}\eta(\dot{x},\dot{x})}{2}\Big)+\,\dot{\epsilon}\eta(\ddot{x},\dot{x})+\,\textrm{higher order terms}.
\label{covariantkineticconstrain3}
\end{align}
and analogous conditions hold for higher derivatives obtained by derivation of the previous ones.
\label{proposiciononkineticconstraints}
\end{proposicion}
{\bf Proof}. The first condition holds by definition.
The second relation is obtained by taking the derivative of \eqref{covariantkineticconstrain1}:
\begin{align*}
2g_{\mu\rho}\,\ddot{x}^{\mu}\dot{x}^{\rho}+\,\frac{d}{d\tau}\big(g_{\mu\nu}\big)\,\dot{x}^\mu\dot{x}^\nu=0.
\end{align*}
From the definition of $\dot{\epsilon}$ and $g_{\mu\nu}$, one has that
\begin{align*}
\frac{d}{d\tau}(g_{\mu\nu}) & = \dot{\epsilon}\eta_{\mu\nu}+\,\textrm{higher order terms}.
\end{align*}
Therefore,
\begin{align*}
2g_{\mu\rho}\,\ddot{x}^{\mu}\dot{x}^{\rho}-\,\dot{\epsilon}\eta_{\mu\nu}\,\dot{x}^\mu\dot{x}^\nu+\,\textrm{higher order terms}=\,0,
\end{align*}
from which follows \eqref{covariantkineticconstrain2}.
The third relation is obtained by deriving \eqref{covariantkineticconstrain2} and taking into account \eqref{covariantkineticconstrain1}:
\begin{align*}
\frac{d}{d\tau}(g_{\mu\nu})\ddot{x}^{\mu}\dot{x}^{\nu}\,+g_{\mu\nu}\dddot{x}^{\mu}\dot{x}^{\nu}\,+ g_{\mu\nu}\ddot{x}^{\mu}\ddot{x}^{\nu}=\,\frac{d}{d\tau}\Big(\frac{{\epsilon}}{2}\eta(\dot{x},\dot{x})\Big)+\,\textrm{higher order terms},
\end{align*}
from which follows the third relation \eqref{covariantkineticconstrain3},
\begin{align*}
g(\dddot{x}^{\rho}\,\dot{x})+\,g(\ddot{x},\ddot{x})& =g_{\mu\nu}\dddot{x}^{\mu}\dot{x}^{\nu}\,+ g_{\mu\nu}\ddot{x}^{\mu}\ddot{x}^{\nu}\\
& =\frac{d}{d\tau}\Big(\frac{{\epsilon}}{2}
\eta(\dot{x},\dot{x})\Big)-\frac{d}{d\tau}(g_{\mu\nu})\ddot{x}^{\mu}\dot{x}^{\nu}\\
& =\,\frac{d}{d\tau}\Big(\frac{\dot{\epsilon}\eta(\dot{x},\dot{x})}{2}\Big)+\dot{\epsilon}\eta(\ddot{x},\dot{x})+\,\textrm{higher order terms}.
\end{align*}
The conditions for higher derivatives are obtained from previous ones by derivation and algebraic manipulations. \hfill$\Box$

As a consequence one has the following approximate coordinate expressions,
\begin{corolario} For a geometry of maximal acceleration $(M,g)$, given the normalization $g(\dot{x},\dot{x})=-1$, the following approximate expressions hold:
\begin{align}
& \dot{x}^{\rho}\,\dot{x}_{\rho} :=\,g_{\mu\rho}\,\dot{x}^{\mu}\dot{x}^{\rho}=-1,\label{kineticconstrain1}\\
&\ddot{x}^{\rho}\,\dot{x}_{\rho}:=\,g_{\mu\rho}\,\ddot{x}^{\mu}\dot{x}^{\rho}=\,
\frac{\dot{\epsilon}}{2}+\,\mathcal{O}(\epsilon^2_0),\label{kineticconstrain2}\\
& \dddot{x}^{\rho}\,\dot{x}_{\rho}+
\,\ddot{x}^{\rho}\,\ddot{x}_{\rho}=\,g_{\mu\rho}\,\dddot{x}^{\mu}\dot{x}^{\rho}+\,
\,g_{\mu\rho}\,\ddot{x}^{\mu}\ddot{x}^{\rho}=
\,\frac{d}{d\tau}\Big(\frac{\dot{\epsilon}\eta(\dot{x},\dot{x})}{2}\Big)-\,\dot{\epsilon}+\,\mathcal{O}(\epsilon^2_0).
\label{kineticconstrain3}
\end{align}
\label{kineticrelatiocoordinate expressions}
\end{corolario}
{\bf Proof}. From Proposition \ref{proposiciononkineticconstraints} and the fact that $g=(1-\epsilon)\eta$, one gets the above expressions as the first order approximation in $\epsilon_0$.\hfill$\Box$
\begin{comentario}
Since we have disregarder from the beginning to consider higher orders in $\epsilon$ and $\delta$, our theory is a linear theory and not a complete perturbative theory and only with validity up to first order in $\epsilon$.
\end{comentario}
\subsection{Causal structure of the metric of maximal acceleration}

For the analysis  of the null sectors of the metric $g$, the expression
\eqref{maximalaccelerationmetric} cannot be used directly, since the factor in the denominator  $\eta(\dot{x},\dot{x})$
is not allowed to be zero. It is more natural to use the relation \eqref{sasakitypemetric}, where $\eta$ can be zero. Thus,
the natural structure of a metric of maximal acceleration is the {\it null set}, that we define as follows,

\begin{proposicion}
The {\it null bundle} $\pi_{NC}:NC\to M$ of the maximal acceleration metric $g$ is characterized by the following type of curves,
\begin{enumerate}
\item Null geodesics, characterized by  $\eta(\dot{x},\dot{x})=0$ and $\eta(D_{\dot{x}}\dot{x},D_{\dot{x}}\dot{x})=0.$
\item Curves $x:I\to M$ characterized by the condition
\begin{displaymath}
\Big(1+ \frac{  \eta_{\sigma\lambda}\,   \ddot{x}^{\sigma}(\tau)
\ddot{x}^{\lambda}(\tau)}{A^2 _{max}\eta(\dot{ x}, \dot{ x})}\Big)=0,\quad\eta(\dot{ x}, \dot{ x})\neq 0.
\end{displaymath}
\end{enumerate}
\label{relationnullconditions}
\end{proposicion}
{\bf Proof}. From the formula of the Sasaki metric \eqref{sasakitypemetric},
it follows that $g(\dot{x},\dot{x})=\,\eta(\dot{x},\dot{x})+\,\frac{1}{A^2_{max}}\eta(D_{\dot{ x}}\dot{ x}, D_{\dot{ x}}\dot{ x})$. Thus, if $g(\dot{x},\dot{x})=0$ and $\eta(\dot{ x},\dot{ x})=0$, it is necessary that:
 \begin{enumerate}
 \item $\eta(D_{\dot{x}}\dot{x},D_{\dot{x}}\dot{x})=0$ or

 \item $A^2_{\max}$ is not bounded.

 \end{enumerate}
 Since we are assuming that the $n$-acceleration is bounded, only the first possibility is applicable, that corresponds to lightlike geodesics.
On the other hand, if $\eta(\dot{ x},\dot{ x})\neq 0$, then the condition $g(\dot{ x},\dot{ x})=0$ is
 equivalent to
 \begin{align*}
 \Big(1+ \frac{  \eta_{\sigma\lambda}\,   \ddot{x}^{\sigma}(\tau)
\ddot{x}^{\lambda}(\tau)}{A^2 _{max}\eta(\dot{ x}, \dot{ x})}\Big)=0,
\end{align*}
 that corresponds to a curve of maximal acceleration.\hfill$\Box$

This  implies that apart from the curves of maximal acceleration, the only null curves compatible with a maximal acceleration geometry, parameterized by the proper-time of $g$, are the null geodesics only.

The natural definition of time orientation in the space $(M,g)$ is the following,
\begin{definicion}
A spacetime $(M,g)$ is time oriented if there is a vector field $W\in \,\Gamma TM$ such that at each point $x\in\,M$ and for each integral curve $x_W:I\to M$ of $W$ with initial condition $x_W(0)=x$, the vector field $W$ is timelike in the sense that $g(W,W)<0$ along $x_W:I\to M$.
\end{definicion}
Thus, a future oriented timelike vector $Z$ is such that if the corresponding integral curve is $x_Z:I\to M$ and $W:I\to Tx_Z$ is the restriction along the curve $x_Z$, then
\begin{align*}
g(W,Z):=g_{\mu\nu}(\,^kx_Z)\,Z^\mu\,W^\nu\,<0.
\end{align*}
In a similar way, a curve $x:I\to M$ is future oriented if it is timelike and the tangent vector is future oriented respect to $g$.

The following result is direct,
\begin{proposicion}
For curves $x:I\to M$ such that the condition \eqref{boundedconditionforacceleration} holds. Then:
 \begin{itemize}
\item Any time orientation $T\in \,\Gamma TM$ of $g$ is a time orientation of $\eta$.

\item The causal character of any vector $Z\in T_xM$ respect to $g$ and $\eta$ is the same.
\end{itemize}
\end{proposicion}

\begin{proposicion}
The subset $T^{--}_xM \subset TM$ of timelike, future oriented tangent vectors respect to $g$ is an open set.
\end{proposicion}
{\bf Proof}. The future pointed sets respect to $\eta$ is
 a strictly convex open  cone on $T_xM$.  The condition $\frac{a^2}{A^2_{max}}$ is open. Therefore, the intersections of both conditions, that is the constraint for $g(\dot{x},\dot{x})\eta(\dot{x},\dot{x})$ is positive and therefore $T^{--}_xM $ is an open set.\hfill$\Box$.

The sector with $\frac{a^2}{A^2_{max}}>1$ corresponds to a change in the causal structure of $g$ respect to the causal structure of the averaged metric $\eta$: any timelike vector with the metric $\eta$ is spacelike with the metric $g$ and viceversa, any timelike vector of $g$ is a spacelike vector of $\eta$. Thus curves of maximal acceleration define the boundary of different signature sectors.

\subsection{Measurable Euclidean length and $n$-velocity in a geometry of maximal acceleration}
By an observer we will mean a timelike, future oriented curve $\mathcal{O}:I\to M$. The
 notion of {\it speed vector} for a timelike trajectory can be defined un-ambiguously as follows. First, one considers the {\it measurable Euclidean distance} between an observer $\mathcal{O}$ and a point $q$ as follows: when the observer $\mathcal{O}$ is at the spacetime point $p$ sends a light signal. The signal reaches the point $q$ and is reflected back to $p'$. The reflected signal light is detected by the observer. The distance $d(\mathcal{O},q)$ between the observer $\mathcal{O}$ and the point $q$
 is defined as one half times the speed of light in vacuum multiplied by the elapsed time $T_{pp'}$ measured by the observer $\mathcal{O}$.

 Now consider two points $p,q\in\,M$.
 \begin{definicion}
 The {\it measurable Euclidean distance} between the points $p,q\in M$ measured by the observer $\mathcal{O}:I\to M$ is defined as
 \begin{align}
 d_E(p,q)=\,|d(\mathcal{O},p)-d(\mathcal{O},q)|.
 \label{Euclideandistanceformula}
  \end{align}
  \label{Euclideandistancedefinition}
\end{definicion}
This definition is applicable even if $p$ and $q$ are not simultaneous respect to $\mathcal{O}$.
The formal consistency requirement for this for this operational definition of distance is that the speed of light is universally constant for all  observer.

As a direct result we have that
\begin{proposicion}
For each observer $\mathcal{O}:I\to M$, the function
\begin{align*}
d_E:M\times M\to R,\quad (p,q)\mapsto |d(\mathcal{O},p)-d(\mathcal{O},q)|
\end{align*}
determines a metric function on $M$.
\end{proposicion}
Until now, the above definitions does not make use of any metric structure and only of the  constancy of the speed of light. Therefore, we adopt the above as the definitions happening in a spacetime of maximal acceleration $(M,g)$.

In a maximal acceleration spacetime $(M,g)$, along a timelike curve $x:I\to M$ the proper time elapsed from $x(\tau)$ to $x(\tau+\delta)$ is given by the formula
 \begin{align*}
 \delta=\,\int^{\tau+\delta}_{\tau}\sqrt{-\eta_{\mu\nu}\dot{x}^\mu\dot{x}^\nu}\,d\tilde{\tau}.
 \end{align*}
 The fact that the geometry is of maximal acceleration changes the notion of proper time respect to the averaged Lorentzian geometry determined by the metric $\eta=\langle g\rangle$. Under the assumption that the observable geometry is the geometry of maximal acceleration, the {\it observable averaged speed} of a curve $x:I\to M$ at the instant $\tau$ measured by an observer $\mathcal{O}$  is defined to be
\begin{align}
v(\tau)=\lim_{\delta\to 0}\frac{1}{\delta}\,{T_{x(\tau)x(\tau+\delta)}},
\end{align}
where the distance ${T_{x(\tau)x(\tau+\delta)}}$ is operationally defined as before for the observer $\mathcal{O}$.
Thus for a geometry of maximal acceleration $(M,g)$, we adopt an analogous definition but replacing $\eta$ by $g$ in the new rule to calculate physical proper times. In this way, we define
\begin{definicion}
Let $(M,g)$ be a maximal acceleration geometry and let $x:I\to M$ be a timelike curve. Then the averaged speed between the points $x(\tau)$ and $x(\tau+\delta)$ along the trajectory $x:I\to M$ is
\begin{align}
{v}(\tau)=\,\lim_{\delta\to 0}\frac{1}{\int^{\tau+\delta}_{\tau}\sqrt{-g_{\mu\nu}\dot{x}^\mu\dot{x}^\nu}\,d\tilde{\tau}}\,{T_{x(\tau)x(\tau+\delta)}}.
\label{avergaedspeed}
\end{align}
\end{definicion}
Using the formula for the maximal acceleration metric \eqref{relationgeta}, one can re-write the averaged speed as
\begin{align*}
v(\tau):=\lim_{\delta\to 0}\frac{1}{\int^{\tau+\delta}_{\tau}\,\Big(\sqrt{1-\frac{a^2}{A^2_{max}}}\,\Big)\,d\tilde{\tau}}\,\tilde{v}(s),
\end{align*}
where $\tilde{v}(\tau)$ stands for the standard definition relativistic speed and with the condition that $x(\tau)=\tilde{x}(s)$\footnote{Note that the map $t\mapsto \tau$ is not a diffeomorphism form $R$ to $R$. Therefore, although $s$ and $\tau$ are valid parameterizations of a curve, they are not reparameterizations of each other.}, using the proper-time of the Lorentzian metric $\eta$,
\begin{align*}
\tilde{v}(s):=\lim_{\delta\to 0}\,\frac{T_{x(s)x(s+\delta)}}{\int^{s+\delta}_{s}\,\sqrt{-\eta_{\mu\nu}\dot{x}^\mu\dot{x}^\nu}\,d\tilde{s}}.
\end{align*}
In order to simplify the treatment, let us consider $a^2(\tau)$ to be constant. Then the above expression leads to
\begin{align}
v(\tau):=\,\frac{1}{\sqrt{1-\frac{a^2}{A^2_{max}}}}\tilde{v}(s).
\label{vtildev}
\end{align}

The vector form of the notion of measurable averaged speed is the {\it measurable $n$-velocity}
\begin{align}
v^\mu(\tau)=\,\frac{1}{\sqrt{1-\frac{a^2}{A^2_{max}}}}\,\tilde{v}^\mu(s),\quad \mu=1,...,n.
\label{fourvelocity}
\end{align}
The limit of Special Relativity is recovered by $A_{max}\to \infty$.

The decomposition of $v$ in temporal and spatial components respect to an inertial observer $\mathcal{O}$ is
\begin{align*}
\frac{v}{c}=\,\frac{1}{\sqrt{1-\frac{a^2}{A^2_{max}}}}\,\frac{1}{\sqrt{1-\frac{\tilde{\vec{v}}^2}{c^2}}}\,(1,\frac{\tilde{\vec{v}}}{c}).
\end{align*}
Taking into account the relation \eqref{vtildev}, we have that for the measurable $n$-velocity the formula
\begin{align}
\frac{v}{c}=\,\frac{1}{\sqrt{1-\frac{a^2}{A^2_{max}}}}\,\frac{1}{\sqrt{1-({1-\frac{a^2}{A^2_{max}})\frac{
\,{\vec{v}}^2}{c^2}}}}\,(1,\sqrt{1-\frac{a^2}{A^2_{max}}}\frac{{\vec{v}}}{c}).
\label{measurablenvelocity}
\end{align}

 \subsection{On the possibility of superluminical motion in spacetimes with metrics of maximal acceleration}
 From an abstract point of view, it seems possible to have superluminical motion in geometries of maximal acceleration. However, when the physical arguments in faubour of maximal acceleration are considered, the impossibility of superluminical motion is concluded.
 \subsection*{The mathematical conditions for superluminical motion}
\begin{proposicion}
Let $(M,g)$ be a metric of maximal acceleration. Then
\begin{itemize}
\item Along a timelike trajectory with parameterizations $x:I\to M,\,\tau\mapsto x(\tau)$ and $\tilde{x}:\tilde{I}\to M,\,s\mapsto \tilde{x}(s)=\,x(\tau)$ as before, one has that
\begin{align}
 v(\tau)\geq \,\tilde{v}(s)
 \label{mayorv}
 \end{align}
 and equality only holds iff $D_{\dot{x}}\dot{x}=0$.
 \item The definition of definition of the measurable speed $v$ is the domain $a^2<A^2_{max}$.
\end{itemize}
\end{proposicion}

Condition \eqref{mayorv} suggests the possibility of superluminical motion in spacetimes with maximal acceleration,
\begin{proposicion}
Let $(M,g)$ be a spacetime of maximal acceleration. Then a necessary condition for superluminical motion $\vec{v}(\tau)>1=c$ consistent with the principle of maximal acceleration in  the physical domain $\eta(\ddot{x},\ddot{x})<A^2_{max}$, $g(\dot{x},\dot{x})=-1$ is
\begin{align}
a|\vec{v}|>\,A_{max} \,c.
\label{conditionofsuperluminical}
\end{align}
 \label{conditionsforsuperluminicalmotion}
\end{proposicion}
{\bf Proof}.
By the vector components in formula \eqref{measurablenvelocity}, one has that
\begin{align*}
\big|\frac{\vec{v}}{c}\big|=\,\frac{1}{\sqrt{1-({1-\frac{a^2}{A^2_{max}})\frac{
\,{\vec{v}}^2}{c^2}}}}\,\big|\frac{\vec{v}}{c}\big|>1.
\end{align*}
After a bit of algebra, this condition can be written as equation \eqref{conditionofsuperluminical}.\hfill$\Box$

This condition is compatible with $\frac{a^2}{A^2_{max}}<1$ if $|\vec{v}|>c$; it is compatible with $\frac{a^2}{A^2_{max}}<1$ iff $|\vec{v}|>>c$. It is not possible for the case $c=|\vec{v}|$ in the perturbative regime.

 Some direct consequences of such kinematics are the following:
  \begin{itemize}
 \item It is not possible to have superluminical motion for geodesic motion.  Thus a point  particle in a pure gravitational field will not be superluminical.

\item In the perturbative regime $a^2<<A^2_{max}$, superluminical motion is possible iff $\beta^{-1}<<1$.
 \end{itemize}
 Finally, let us mention that the quantity $P_{c}=\,mA_{max}\,c$ has the dimensions of a power, that we can call {\it critical power}. Therefore, superluminical condition \eqref{conditionofsuperluminical} requires the existence of powers on an electron bigger than the critical power.
\subsection*{The superluminical motion from the heuristic point of view}
The heuristic argument in the beginning of the {\it section}, if one needs to consider it seriously, implies the existence of a maximal acceleration from two assumptions:
\begin{itemize}
\item Existence of a maximal speed for the propagation of physical interactions,

\item Existence of a minimal length for the characteristic phenomena (a kind of generalized Debye  length).
\end{itemize}
Thus, if maximal acceleration is a consistent requisite from a minimal length and maximal speed, it is not possible to have superluminical world-lines in a geometry of maximal acceleration. A further argument in the context of electromagnetic radiation will be given in {\it section 6}, when we obtain the value of the maximal acceleration for point charged particles.

\subsection{Isometry group of the metric of maximal acceleration and null-structure}

 The assumption that the speed of light is constant and independent of the observer $\mathcal{O}$ has been of fundamental relevance to obtain a covariant notion of $n$-velocity vector which is independent of parameterizations. The proper parameterization of physical world-line curves is determined by the metric of maximal acceleration $g$. Doing this, the principle of maximal acceleration is full-filled. The speed of light has been used as the standard to measure distances. Therefore, the {\it Principle of constancy of the speed of light in vacuum} also holds in our theory. Therefore, the generalized metric structure $g$ should have an isometry group that leaves invariant the null-cone structure. By {\it Proposition} \ref{relationnullconditions}, the null-cone of $g$ is the null-cone of $\eta$. Therefore,
   \begin{proposicion}
   Given a metric of maximal acceleration, the group leaving the null-cone structure of $g$ invariant is a Lie group (conformal group in dimension $n>0$).
   \end{proposicion}
 Also, the isometry group of the metric $g$ contains the isometry group of the Lorentzian metric $\eta$,
 \begin{align}
 Iso(g)\supset Iso(\eta).
 \end{align}
  The factor $(1-\frac{a^2}{A^2_{max}})$ explicitly breaks the group of conformal transformations (that leaves invariant the light-cone). In the case that $A_{max}\to \infty$, the group leaving invariant the cone structure is the conformal group in $n$-dimensions.
  In the case that $\eta$ is the Minkowski metric, if $Iso(\eta)=\, O(1,n-1)$, then  $Iso(g)\supset O(1,n-1)$ and then the proper-time of a metric of maximal acceleration is Lorentz invariant. Indeed, this is the largest group possible, since the acceleration factor

\section{Higher order generalization of the electromagnetic field and current}

In this {\it section} we introduce the notion of generalized electromagnetic fields as sections of a  bundle $\Lambda^2(M,\mathcal{F}({J}^k_0(M)))$ for an integer $k\geq 2$. The value of the integer $k$ will be fixed later. We will use the results on exterior algebra and cohomology from {\it section 2}. Due to the isomorphisms  (\ref{anotherisomorphism}) and \eqref{isomorphismtildej3j2}, the fields can also be considered as sections of $\Lambda^{2+nk}(J^{k}_0(M))$ and $\Lambda^{2+nk+1}(\tilde{J}^{k+1}_0(M))$. However, in order to keep a general formulation of generalized higher order field, we will develop in this {\it section} the formalism at the level of forms in $\Lambda^2(M,\mathcal{F}({J}^k_0(M)))$.

It is well known that in standard Maxwell electrodynamics, the electromagnetic field of a point charged particle is divergent, with a singularity of Coulomb type. For the generalized higher order fields we will assume that also contain singularities of Coulomb type and that these are the only singularities that appear. Under such hypothesis, one can regularize the fields, with the result that the fields appearing in the equation of motion of a point charged probing particle are finite.

 In Maxwell's electrodynamics the fields are divergent along the world-line of the point charged particle, which is the submanifold
the submanifold $e(I)=S\hookrightarrow {M}$, with $e:I\to {M}$ the world-line. $S$ is dynamically determined by the equation of motion of the point particle. In this case, the spacetime is $\tilde{M}:=M\setminus S$ is by definition the region where the fields are finite.

After renormalization of mass procedure, all the fields that appear in the equations of motion are finite. Therefore, there is not need to subtract the submanifold $S$ from the domain of definition of the fields. This is the procedure that we will follow in this {\it section}. We first start with fields that are divergent on $S$ and after the introduction of a regulation procedure (that coincides with Dirac's procedure), the fields will defined in the whole manifold $M$.

 Finally, since we will use the exterior derivative of generalized forms, the notion of generalized higher order field that we will use is in the strong sense of {\it definition} \ref{generalizedfield}. This will be implicitly understood in the rest of the work.

 \subsection{Generalization of the electromagnetic field as sections of \\$\Lambda^2(M, \mathcal{F}(J^k_0 (M)))$}

 The generalization of the Faraday form is given by the following
\begin{definicion}
Given a curve $x:I\to M$, the electromagnetic field $\bar{F}$ along the lift $^kx:I\to J^k_0(M)$ is a closed $2$-form
$\bar{F}\,\in \Gamma\Lambda^2 (M\setminus{S}, \mathcal{F}(J^k_0( M\setminus{S})))$.
\end{definicion}
Thus, in a local natural coordinate system, the generalized Faraday form $\bar{F}$ can be written as
\begin{equation}
\bar{F}(\,^kx)=\bar{F}(x,\dot{x},\ddot{x},\dddot{x},...)=\Big(\big(F_{\mu\nu}(x)\big)+ \Upsilon_{\mu \nu}(x,\dot{x},\ddot{x},\dddot{x},...)\big)dx^{\mu}\wedge dx^{\nu}\Big),
\label{electromagneticfield}
\end{equation}
with $F(x)\in \Gamma\Lambda^2 M$.
\begin{comentario} We have the following remarks:
\begin{itemize}
\item The field $\bar{F} \in\,\Gamma\Lambda^{2} (M\setminus{S}, \mathcal{F}(J^k_0( M\setminus{S})))$. Other related fields appearing in the electromagnetic theory will  be  generalized to sections higher order jet bundles.

\item $\tilde{\varphi}(\bar{F})$ is closed, $d_4 \,\tilde{\varphi}(\bar{F})=0$ and therefore it defines an element of the cohomology $H^*_{cv}(J^k_0(M))$,
\begin{align*}
[\bar{F}]\mapsto [\tilde{\varphi}(\bar{F})]=[\,^k\xi(\bar{F})\wedge dvol_V].
\end{align*}
These two notions for the electromagnetic field are equivalent.

\item $x:I\to M$ must not be interpreted physically as the world-line of  the particle generating the field. Indeed the curve $x:I\to M$ corresponds to the world-line of a point charged particle whose motion can experimentally be observed. Comparing such trajectories with the free particle world-lines, one should be able to identify the effect of the full electromagnetic field.

\item Since the field $F(x)$  lives on $M$, it is clear that $d_4(\varphi(F))=\,\varphi(dF)$.
\end{itemize}
\end{comentario}

Using the cohomology theory of forms on $\Lambda^p(M,\mathcal{F}(J^k_0(M)))$ one can establish the following result,
\begin{teorema}
Let $\bar{F}\,\in\,\Gamma \Lambda^{2} (M\setminus{S}, \mathcal{F}(J^k_0( M\setminus{S}))),$ $d_4\bar{F}(\,^kx)=0$ be a generalized higher order electromagnetic field. Then the decomposition \eqref{electromagneticfield} with $dF=0$ is unique.
\label{decompositionofelectromagnetic field}
\end{teorema}
{\bf Proof}. By the uniqueness in (\ref{descomposicion de alpha}), the field $F(x):=F_{\mu\nu}dx^{\mu}\wedge\, dx^\nu$ in \eqref{electromagneticfield} can be identified with $\varphi(\langle \bar{F}\rangle)$ and since $\varphi^*$ is an isomorphism, $F(x)$ can also be identified with $\langle\bar{F}\rangle$ up to gauge in $H^{2+kn}_{cv}(M\setminus{S},\mathcal{F}(J^k_0(M\setminus{S})))$. This implies that
\begin{align*}
\bar{F}(x)=\langle \bar{F}(x)\rangle \,+\langle d_J A\rangle,\, A\in \Lambda^{1+kn}(M\setminus{S},\mathcal{F}(J^k_0(M\setminus{S}))).
\end{align*}
Since $\langle \bar{F}(x)\rangle $ is unique, then $\langle d_J A\rangle =0$,
\begin{align*}
\langle \bar{F}\rangle=\,\langle \varphi(\langle \bar{F}\rangle)\rangle\,+\langle d_JA \rangle =\,\langle\bar{F}\rangle\,+\langle d_JA\rangle =F(x).
\end{align*}
\hfill$\Box$
\begin{proposicion}
For a $n$-dimensional manifold,
if $[\bar{F}\wedge dvol_V]\in\, \,H^{2+kn}_{cv}(J^k_0(M\setminus{S})$, then $[F]\in\,H^2_{dR}(M\setminus{S})$.
\end{proposicion}
{\bf Proof}. It is a consequence of Thom's isomorphism theorem  \ref{isomorphismthom}. \hfill$\Box$

From the above considerations, it follows the following formula for $\bar{F}$,
\begin{align}
\bar{F}=\,\varphi \langle \bar{F}\rangle+\,\Upsilon.
\label{estructuradebarF}
\end{align}
Let us consider $^k\zeta(\bar{F}) \wedge dvol_V$. Then it follows that
\begin{align*}
0=\,d_J (\,^k\zeta\bar{F}\wedge dvol_V) & =\,d_J\big((\,^k\zeta\varphi(F)\,+\,^k\zeta(\Upsilon))\wedge dvol_V\big)\\
& =\big(\,^k\zeta\,d_4 (\varphi(F))\,\wedge dvol_V\,+d_J(\,^k\zeta\Upsilon)\big)\\
& =\,^k\zeta d_4 (\varphi(F))\,\wedge dvol_V\,+d_J(\Upsilon_{\mu \nu}dx^{\mu}\wedge dx^{\nu}\,\wedge dvol_V).
\end{align*}
Since $d_4\bar{F}=0,$ it follows the relation
\begin{align}
d_4\varphi(F)=-\,d_4\Upsilon.
\end{align}

\subsection{The generalized excitation tensor $\bar{G}$}

Similarly to the case of the generalized Faraday tensor $\bar{F}$, the generalized excitation tensor is described by a generalized two form $\bar{G}$,
\begin{definicion}
The excitation tensor along the charged point particle $^kx$ is a $2$-form $\bar{G}\in \Gamma \Lambda^2 (M\setminus{S}, \mathcal{F}(J^k_0( M\setminus{S}))),$
\begin{equation}
\bar{G}(\,^kx)=\bar{G}(x,\dot{x},\ddot{x},\dddot{x},...)=\big(G_{\mu\nu}(x)+ \Xi_{\mu \nu}(x,\dot{x},\ddot{x},\dddot{x},...)\big)dx^{\mu}\wedge dx^{\nu}.
\label{excitationtensor}
\end{equation}
\end{definicion}
The determination of the form $\bar{G}$ from the generalized Faraday tensor $\bar{F}$ is provided by the {\it constitutive relation},
\begin{align}
\bar{G}=\,\star \bar{F}.
\label{constitutiverelationinvacuum}
\end{align}
The choice of this constitutive \eqref{constitutiverelationinvacuum}  is for the vacuum. Other choices can be possible but will not be explored in this work. Note that
the form $\bar{G}$ is not necessarily closed with $d_4$.

Let us consider the star operator associated with $g$. In this case, the field $G(x)=G_{\mu\nu}dx^{\mu}\wedge\, dx^{\nu}$ is defined such that
\begin{align*}
G(x):=\star \varphi(\langle \,\bar{G}(u)\rangle ).
\end{align*}
Therefore, the decomposition in the excitation tensor is unique. Note that a equivalent description of the excitation field is
\begin{align*}
\tilde{\varphi}(\bar{G})=\bar{G}\wedge dvol_V.
\end{align*}

There are some advantages in generalizing the electromagnetic field in this way. For instance, defining the tensors $\Upsilon$ and $\Xi$ as differential forms allow us to perform integrals that are diffeomorphism invariant and to use the machinery of exterior calculus in an analogous way as in standard classical electrodynamics. Also, it is straightforward to generalize Maxwell's equations for the fields $\bar{F}$ and $\bar{G}$.

\subsection{Higher order charge current density}

The density current in electrodynamics is represented by a $d_4$-closed $3$-form
\begin{align*}
J\in \,\Lambda^3 (M\setminus{S}, \mathcal{F}(J^k_0( M\setminus{S}))).
\end{align*}
This can be generalized to
\begin{equation}
\bar{J}(\,^kx)=\bar{J}(x,\dot{x},\ddot{x},\dddot{x},...)=\,J(x)\,+\Phi(x,\dot{x},\ddot{x},\dddot{x},...).
\end{equation}
For sources corresponding to point charged particles, the current density $J$  is a distribution with support on the embedding $S:I\to M$. Therefore, one requires that the support of $\bar{J}$ also lives on the lift $^k{S}:I\to J^k_0(M)$.

We will see that the current density $\bar{J}$ satisfies
\begin{align}
d_4 \bar{J}=0
\label{conservationofcharge}
\end{align}
is a consequence of the minimal extension of the fields (see {\it section 7}).
This relation generalizes the charge conservation law in electrodynamics.

\subsection{Geometric description of a point charged particle and other charge configurations}

Let $(M,\eta)$ be a time orientable spacetime and consider $g$ the metric of maximal acceleration, defined as in {\it Section 3}. The space of world-line curves is
\begin{align*}
 \mathcal{C}(M):=\left\{x:I\to M,\quad g(\frac{dx}{d\tau},\frac{dx}{d\tau})<0,\quad \,\dot{x}\, \textrm{ future oriented}\right\}.
\end{align*}
 One point particle is described by a curve $x:I\to M$. In the case that a system has more than one point particle, each of them are described by the disjoint union of curves, describing the evolution of each of such particles.
The world-line $x:I\to M$  of a point charged particle has an associated lift
\begin{align*}
^kx:I_x\to J^k_0(M).
\end{align*}
Thus the set of lift to $J^k_0(M)$ corresponding to  $\mathcal{C}(M)$  is
\begin{align}
 ^k\mathcal{C}(M):=\left\{\,^kx:I\to J^k_0(M),\quad g(\frac{dx}{d\tau},\frac{dx}{d\tau})<0,\quad \dot{x}\,\textrm{ future oriented}\right\}.
\end{align}
\label{pointchargedparticle}

\subsection{Operator norms associated with generalized metrics}

Let us first consider the operator norm  associated with the generalized metric $g$ of the endomorphism $ \bar{O}\in \,\Gamma T^{(1,1)}(M,\mathcal{F}(J^k_0(M)))$,
\begin{align}
\|\bar{O}\|_{g}:=\,\sup \Big\{\,\frac{\|\bar{O}\,u\|_{g}}{\|u\|_{g}}, \, u\neq 0\Big\}.
\label{operatornorm}
\end{align}
Defined in this way, $\|\bar{O}\|_{g}$ is intrinsic in the sense that is independent of the local frame or local coordinates that we use and only depends on $\bar{O}$, the symmetric form ${g}$ and the vector $V$.
\begin{proposicion}
For generalized metrics $g(\,k^x)=\,\lambda(\,^kx) \eta(x)$, it holds that
\begin{align*}
\|\bar{O}\|_{g}=\|\bar{O}\|_{{\eta}}.
\end{align*}
\end{proposicion}
\begin{definicion}
Given two operators $ \bar{O}_1,\,\bar{O}_2\, \in \,\Gamma T^{(1,1)}(M,\mathcal{F}(J^k_0(M)))$. Then we say that
\begin{align*}
O_1\,<O_2\,\,
 \textrm{iff} \,\,\|\bar{O}_1\|_{g}\,<\|\bar{O}_2\|_{g}.
 \end{align*}
\label{orderoperators}
\end{definicion}
Let us consider the isomorphism induced by the generalized metric $g$ given by (\ref{definiciondekappa}).
 In particular this expression determines the norm of the generalized $2$-form $\bar{F}\in \Lambda^2(M,\mathcal{F}(J^k_0(M\setminus{S})))$,
\begin{align}
\|\bar{F}\|_{g}:=\|\kappa(\bar{F})\|_{g}.
\label{normF}
\end{align}
Also, note that the norm $\|,\|$ induces a pre-order relation in $\bar{F}\in \Lambda^2(M,\mathcal{F}(J^k_0(M\setminus{S})))$. Such relation is useful to compare the strength of the $2$-forms.
\subsection{Definition of generalized electric and magnetic fields}

 Given a timelike vector field $W\in\,\Gamma TM$ associated to an observer,
 the generalized electric field is defined in a similar way as in the standard case,
 \begin{align}
 \bar{E}:= \iota_W \bar{F}.
 \label{definitionofelectricfield}
 \end{align}
 The norm of $\|\kappa(\bar{F})|_{g}$ coincides with the norm $\|\cdot\|_{g}$ of the {\it electric field},
\begin{align*}
\|\bar{E}(u)\|_{g}=\|\iota_W \bar{F}(u)\|_{g}=\|\kappa(\bar{F})(u)\|_{g},\,\quad u\in \,^k\pi^{-1}(x).
\end{align*}
The interpretation of $E=\iota_W \bar{F}$ as the electric field depends on the local coordinate system (determined by the integral curves of $W$): if we change the observer to $\tilde{W}$, the electric field will be in principle different.
\begin{ejemplo}
For a point charged particle, outside the world-line $S=x(\sigma)$, the {\it electric field} corresponding to $\bar{F}$ measured in the coordinate system adapted to the motion of the charged point particle is
\end{ejemplo}
\begin{align*}
\|\bar{E}_u\|_{g}:=\|\bar{F}\|_{g}(u)=\frac{1}{r^2},
\end{align*}
where $r=d_g(x,S)$ is the distance from $x$ and $S$ measured with the positive definite generalized metric $g_+$ defined as in \eqref{positivebilinearform}.

In a similar way, the definition of the generalized magnetic field is
\begin{align}
\bar{B}=\iota_W\,\star \bar{F}.
\end{align}
\begin{proposicion}
An electromagnetic field is zero iff there is an observer $W$ such that the electric and magnetic field are both zero.
\end{proposicion}
\subsection{Short distance behavior of the electromagnetic field}

Let  $r=d(x,S)$ be the distance function from the point $x\in M$ to the world-line $S\hookrightarrow M$ using the Riemannian function $g_+$ as in equation \eqref{positivebilinearform}.
The expression of the Coulomb field for a point charged particle suggests that we should consider electromagnetic fields $\bar{F}$ with the following development in powers of $r$,
\begin{align}
\| \bar{F}\|_{\hat{g}}(x)=a_{-2}\frac{1}{r^2}\,+a_{-1}\frac{1}{r}\,+\,\sum^{+\infty}_{k=0} a_k \,r^k,
\label{Expansionoftheelectricfieldinr}
\end{align}
where each of the functions
\begin{align*}
a_i:J^k_0(M)\to M,\,i=-2,-1,0,1,...,+\infty
\end{align*}
 are homogeneous of degree zero in $r$ and smooth functions on $J^k_0(M)$. For such generalized higher order fields, given a point $x(\sigma_0)\in S$ that is an isolated singularity of the field $\| \bar{F}\|_{\bar{g}}$ and a small sphere $S^2_{(x(\sigma_0),r)}$ surrounding $x(\sigma_0)$, one has the relation
\begin{align}
4 \pi Q=\int_{S^2_{(x(\sigma_0),r)}} \star \bar{F}
\label{definitionofcharge}
\end{align}
 If one imposes that it must not depend on $r$ for $r$ finite and small enough, one has constraints on the {\it averaged values} of the coefficients,
 \begin{align}
 \int_{S^2_{(x(\sigma_0),r)}} \,a_i(\,^kx)\, d\Omega =0,\,i=-1,...,\infty.
 \end{align}
 If we impose the stronger condition that the integrals along any closed surface surrounding $x(\sigma_0)$ must be  $4 \pi Q$, one obtains the conditions
 \begin{align}
 a_i(\,^kx)=0, \,i=-1,...,\infty.
 \end{align}
We adopt the convention that the quantity $Q$ in equation \eqref{definitionofcharge}
 corresponds to the total charge of the point charged particle whose world-line is the singularity $S$.

\subsection{Analytic structure of the generalized electromagnetic fields}

Given a charged point particle with world-line $x:I\to M$ with $k$-lift  $^kx:I\to M$, the
 electromagnetic field outside from the singularity region $S$ can be decomposed as
\begin{align*}
\bar{F}(x,\dot{x},\ddot{x},\dddot{x},...)& =\,\varphi\big(\big(F^C_{\mu\nu}(x)+\,F^{D}_{\mu\nu}(x)\big)dx^{\mu}\wedge dx^{\nu}+\,F^{ext}_{\mu\nu}(x)dx^{\mu}\wedge dx^{\nu}\big)
\end{align*}
\begin{align}
& +\big(\Upsilon^{div}_{\mu \nu}(x,\dot{x},\ddot{x},\dddot{x},...)\,+\Upsilon^{reg}_{\mu \nu}(x,\dot{x},\ddot{x},\dddot{x},...)\big)dx^{\mu}\wedge dx^{\nu}.
\label{Expansionoftheelectricfield}
\end{align}
The piece $F^{ext}(x)$ corresponds  to the contribution to the  field not generated by the singularity on $S$.
The piece $\Upsilon^{reg}_{\mu \nu}(x,\dot{x},\ddot{x},\dddot{x},...)$ and $\Upsilon^{div}_{\mu \nu}(x,\dot{x},\ddot{x},\dddot{x},...)$
are the regular and divergent pieces of the field $\langle \bar{F}\rangle$ on the singularity region $S$.
 $F^C_{\mu\nu}(x)$ is the divergent field on the world-line.

The behavior for the divergent field at short distances is
\begin{align}
\|F^C_{\mu\nu}(x)\|_{\hat{g}}=\,a_{-2}\,\frac{1}{r^2}.
\end{align}
On the other hand, the {\it field} $F^{D}_{\mu\nu}$ is regular at short distances.
\begin{proposicion}
The following relations hold
\begin{enumerate}
\item
The singularity of $\bar{F}$ is such that
\begin{align*}
a_{-2}=\lim_{r\rightarrow 0} r^2\, \bar{F}=\,\lim_{r\rightarrow 0} r^2\, {F}^C,
\end{align*}
\item The singularities of $\Upsilon$ are of the form
\begin{align*}
\lim_{r\rightarrow 0}r^2\,\Upsilon (\,^kx)=0.
\end{align*}
\end{enumerate}
\end{proposicion}
Then one has the relation
\begin{align}
F^C_{\mu\nu}(x)=a_{-2}\, \frac{1}{r^2} \theta_{\mu \nu}(x),\,r\neq 0,
\label{Coulombfield}
\end{align}
where the tensor $\theta_{\mu\nu}(x)$ is skew-symmetric and homogeneous of degree zero in $r$.

The higher order piece  $\Upsilon^{div}_{\mu \nu}$ of $\bar{F}$ must be such that asymptotically
\begin{equation}
 \big(F^C_{\mu\nu}(x)\,+\Upsilon^{div}_{\mu \nu}(x,\dot{x},\ddot{x},\dddot{x},...)\big)\to F^C_{\mu\nu}(x),
\end{equation}
when $r\to 0$.
Therefore, we can assume that
\begin{align}
\Upsilon^{div}_{\mu \nu}(x,\dot{x},\ddot{x},\dddot{x},...)=0.
\label{align}
\end{align}
This is in concordance with Dirac's result on the regularity of the radiation field.
\subsection*{The generalized radiation field}

The standard
radiation field holds the relation
\begin{align}
\,F^{rad}(\,^kx)=\,\bar{F}(\,^kx)-F^C(x)\,-F^{ext}(x)-\Upsilon^{reg}(\,^kx).
\label{radiationfield}
\end{align}
Following Dirac \cite{Dirac}, for the standard Maxwell-Lorentz theory the {\it radiation field} is finite on the whole spacetime $M$. In such theory, the  radiation reaction field,
\begin{align}
F^{rad}(x)=\,F^{ret}\,-F^{adv}.
\end{align}
with $F^{ret}$ and $F^{adv}$ are the retarded and advanced fields, obtained from the corresponding Li\'{e}enard-Wiechert potentials.
A balance equation analysis provides the following value for the radiation field on the world-line $x:I\to M$ $\cite{Dirac}$,
\begin{align}
F^{D}_{\mu\nu}=\frac{4}{3}\,\Big(\frac{d^3x_{\mu}}{ds^3}\frac{d^2x_{\nu}}{ds^2}
-\,\frac{d^3x_{\nu}}{ds^3}\frac{d^2x_{\mu}}{ds^2}\Big),
\label{Diracfields}
\end{align}
where $s$ is the proper-time along $S$ calculated using the metric $\eta$ and $\frac{dx_{\mu}}{d s}=\,\eta_{\mu\nu}\,\frac{dx^{\nu}}{d s}$. This field is finite. However, the parameter $s$ must be substituted by the parameter $\tau$ in the expression of the radiating field. This is because we are considering that the standard clocks associated with a world-line are measured with $g$ and not with $\eta$. This provides the following expression for the radiation field,
 \begin{align}
 F^{rad}(\,^3x)=\frac{4}{3}\,\Big(\frac{d^3x_{\mu}}{d\tau^3}\frac{d^2x_{\nu}}{d\tau^2}
-\,\frac{d^3x_{\nu}}{d\tau^3}\frac{d^2x_{\mu}}{d\tau^2}\Big)+\mathcal{O}(\epsilon^2_0).
\label{Diracfield}
 \end{align}
 We will adopt this value of the field for the generalized radiation field. Note that adopting this value of the field, we are assuring conservation of energy momentum and that the field is a solution of the corresponding Maxwell's equations. Also, note that considering $F^{rad}\in \Gamma \, \Lambda^2(M,\mathcal{F}(J^3_0(M)))$, we are reinterpreting the standard radiation field as a generalized field with values in a higher order jet bundle.

For the rest of the paper the regular part $\Upsilon^{reg}$ will be just denoted by $\Upsilon$ (and similarly for the regular part for $\Xi$, $\Xi^{reg}$).
Let us consider the field generated by a point charged particle, which is also considered
to be the {\it absorber}. If $F^{ext}=0$, out from the world-line $x:I\to M$
all the fields in (\ref{electromagneticfield}) are finite and one can write
\begin{equation}
\bar{F}(x,\dot{x},\ddot{x},\dddot{x},...)=\big(F^C_{\mu\nu}(x)+\,F^{rad}_{\mu\nu}(\,^3x)+\,\Upsilon_{\mu \nu}(x,\dot{x},\ddot{x},\dddot{x},...)\big)dx^{\mu}\wedge dx^{\nu}.
\end{equation}
The pieces $\,F^{rad}_{\mu\nu}(\,^3x)$ and $\Upsilon_{\mu \nu}(x,\dot{x},\ddot{x},\dddot{x},...)$ are analytic functions of the Euclidean distance to the world-line. Thus, at zero order approximation in the distance to the particle world line, one can write the relations for the regularized radiation field (the Coulomb field does not contribute to the radiation),
\begin{align}
\bar{F}^{rad}(x,\dot{x},\ddot{x},\dddot{x},...):=\big(F^{rad}_{\mu\nu}(\,^3x)+ \Upsilon_{\mu \nu}(x,\dot{x},\ddot{x},\dddot{x},...)\big)dx^{\mu}\wedge dx^{\nu}.
\label{regularizedradiationfield}
\end{align}
The constitutive relation $G=\star F$ makes natural to consider
\begin{align}
G^{rad}:=\,\star\,F^{rad},
\end{align}
where the star operator is associated with the metric of maximal acceleration $g$.
Thus, for the excitation tensor one has the relation
\begin{align}
\bar{G}^{rad}(x,\dot{x},\ddot{x},\dddot{x},...):=\Big(G^{rad}_{\mu\nu}(\,^3x)+ \Xi_{\mu \nu}(x,\dot{x},\ddot{x},\dddot{x},...)\Big)dx^{\mu}\wedge dx^{\nu}.
\label{regularizationexcitationfield}
\end{align}
All the fields in (\ref{regularizedradiationfield}) and (\ref{regularizationexcitationfield}) are smooth on $M$.

Note that although $F^{rad},G^{rad}\in \Gamma \, \Lambda^2(M,\mathcal{F}(J^3_0(M)))$, we did not fix the value of $k$ where the full fields \eqref{regularizedradiationfield} and \eqref{regularizationexcitationfield} are defined. However, the fact that $F^{rad}$ and $G^{rad}$ depends on
the third order jet bundle strongly suggests that $k=3$. Since there is not physical distinction (that is, an operational identification by experiment using prove particles) on the field $\bar{F}$ between the contributions coming from $F^{rad}(\,^3x)$ and $\Upsilon(\,^kx)$. This can be generalized to the case with external fields. On the other hand, the Coulomb singularity will be renormalized in the mass. The value $k=3$ will be confirmed in {\it section 3} after we obtain the consistent equation of motion for a point charged particle interacting with generalized higher order fields.

\subsection{Physical interpretation of the generalized higher order fields}

In {\it section 1} we motivated the introduction of the generalized higher order fields based on the criticism of the concepts of external field-test particle system. The solution that we suggest is to substitute both notions by a new notion of field and probe particle such that the field
depends on the state of motion of the probe particle. The notion of generalized higher order field  accommodates naturally to the ideal of an {\it operational field theory}, in the sense that the mathematical structures must be linked with observables in a minimal economical way.

Thus, in the contest of Electrodynamics, generalized higher order fields are associated with the physical field that a probe particle will interact with. This is particularly clear for the generalized Faraday form $\bar{F}$, since it is that field which will entry on the equation of motion of the particle. On the other hand, the generalized excitation tensor $\bar{G}$ being also a generalized field is a postulate, which obviously natural from the point of view of the symmetries of the theory. However, that the excitation tensor $\bar{G}$ depends on the state of motion of the probe particle suggests that not only the electrodynamics but also matter fields should be considered from the perspective of generalized higher order fields.

 Finally, let us clarify that our notion of generalized higher order field is not a generalization of the notion of Finsler field theories (see for instance \cite{VacaruStavrinos, Voicu} among the recent contributions). While our theory is originated from an aim of maximal economical postulate in the formalism of field theory, the introduction of  Finsler field theory is mainly based on arguments coming from quantum gravity phenomenology.

\subsection{Physical interpretation of the relation $F=\,\langle\,^k\zeta(\bar{F})\rangle $}

 Since $F(x)$ lives on the spacetime manifold $M$, it is natural to view $F$ as the standard  electromagnetic field, outside the world-line $S$. However, this interpretation is not appropriated. First of all, such interpretation is in conflict with the original motivation to introduce the field $\bar{F}$. In the framework of generalized higher order fields, it is not longer valid the notion of external field and test particle. Only in special physical cases, when the field $\bar{F}-\varphi{\langle \bar{F}\rangle }$ is small compared with $\bar{F}$ in the sense of the norm (\ref{operatornorm}), one can approximate the pair $(\bar{F},\,^kx(s))$ by the pair $(F, x(s))$ and ascribe to both independent physical reality. In such case, the language of external fields and test particles is useful to describe the physical phenomena.

An alternative interpretation of the field $F$ comes from its definition as a result of an integration along the fiber. If $F=\langle\,^k\zeta(\bar{F})\rangle $, then it can be thought as the expected value of a physical measurement. In addition, the statistical distributions that one considers are compact distributions along the fibers. The integration is performed using a solution of a kinetic model but defined by probability distribution functions living in higher order jet bundles.  Therefore, the physical interpretation of the relation $F=\,\langle\,^k\zeta(\bar{F})\rangle $ is statistical, related with the measurement of fields using {\it bunch of particles} instead of individual point particles.

\section{The Lorentz-Dirac equation in the framework of higher order electromagnetic fields and maximal acceleration geometry}

In this {\it section} we derive the standard Lorentz-Dirac equation using a simple argument by F. Rorhlich. Then we reconsider the derivation in the framework of generalized higher order fields. In this two derivations the spacetime $M$ will be four dimensional and the metric is the Minkowski metric $\eta$ with signature $(-1,1,1,1)$. Thus, all the contractions and lowering indices operations performed in this {\it section} are performed with $\eta$. The parameter $\tau$ is the proper time along a given curve respect to the Minkowski metric $\eta$. The calculations are performed in a normal coordinate system of $\eta$, where it has the diagonal form $(-1,1,1,1)$. Finally, we repeat the calculation but with a spacetime endowed with a metric of maximal acceleration $g$.

\subsection{A simple derivation of the Lorentz-Dirac force equation}

The Lorentz-Dirac equation describes the motion of a point charged particle interacting with its own electromagnetic and with an {external field} or force. In a normal coordinate system of $\eta$, if the external force is the Lorentz force on a particle, the Lorentz-Dirac is the third order differential equation
\begin{equation}
m\,\ddot{x}^{\mu}=\,eF^{\mu}\,_{\nu}\,\dot{x}^{\nu}+\,\frac{2}{3} e^2\,\big(\dddot{x}^{\mu}-(\ddot{x}^{\rho}\,\ddot{x}_{\rho})\dot{x}^{\mu}\big),\quad \ddot{x}^{\mu}\ddot{x}_{\mu}:=\ddot{x}^{\mu}\ddot{x}^{\sigma}\eta_{\mu\sigma}.
\label{lorentzdiracequation}
\end{equation}
This equation contains run-away and
pre-accelerated solutions \cite{Dirac}, both against what is observed in everyday experience and in contradiction
with Newton's fist law of classical dynamics.

We  present a
simple derivation of the equation (\ref{lorentzdiracequation}) based on the geometric method of adapted local frames. This derivation is what we have called {\it Rohrlich's argument} \cite{Rohrlich}.
It illustrates a method that we will use in the next {\it section} in the context of generalized higher order fields.
One starts with the Lorentz force equation for a point particle interacting with an electromagnetic field $F_{\mu\nu}$ \footnote{Indeed one can use the same argument if instead of the Lorentz force there is an external force orthogonal to the $4$-velocity.},
\begin{align}
m_b \,\ddot{x}^{\mu}\,=\, e F_{\mu \nu}\, \dot{x}^{\nu},
\label{Rohrlichrelation1}
\end{align}
where $m_b$ is the {\it bare mass} and $e$ the electric charge of the particle.
Both sides of \eqref{Rohrlichrelation1} are consistently orthogonal to $\dot{x}$.
If one wants to generalize the equation \eqref{Rohrlichrelation1} to have into account the radiation reaction,
one can add to the right side a vector field along the curve $x:R\to {\bf M}$. This is a map $Z:I\to J^k_0(M)$ such that the diagram
\begin{align*}
\xymatrix{ &
{ J^k_0(M)} \ar[d]^{^k\pi}\\
{ I} \ar[ur]^{Z}  \ar[r]^{x} & { M}}
\end{align*}
commutes.
The orthogonality condition
\begin{align}
\eta (Z(\tau),\dot{z}(\tau))=0
\label{ortogonalityconditionfoZeta}
\end{align}
 implies the following general expression for $Z$,
\begin{align}
Z^{\mu}(\tau)=P^{\mu}\,_{ \nu}(\tau)\big(a_1\dot{x}^{\nu}(\tau)+a_2\,\ddot{x}^{\nu}(\tau)+a_3\,\dddot{x}^{\nu}), \quad P_{\mu \nu}=\eta_{\mu \nu}+\dot{x}_{\mu}(\tau)\dot{x}_{\nu}(\tau),\quad \dot{x}_{\mu}=\eta_{\mu\nu}\,\dot{x}^{\nu}
\label{newcontribution}
\end{align}
with $a_1$, $a_2$ and $a_3$ a priori arbitrary.
Using the orthogonality \eqref{ortogonalityconditionfoZeta} we can write $a_1=0$, since this contribution will not appear in the right hand side if $P^\mu\,_\nu\dot{x}^\nu=0$. Then using the kinematical relations, one obtains
\begin{align*}
\dddot{x}^{\rho}\,\dot{x}^{\sigma}\eta_{\rho\sigma}\,=\,- \ddot{x}^{\rho}\,\ddot{x}^{\sigma}\eta_{\rho\sigma},\quad \ddot{x}^{\mu}\dot{x}_{\mu}=0,
\end{align*}
from which follows the relations
\begin{align*}
Z^{\mu}(\tau)=\,a_2 \ddot{x}^{\mu}(\tau)+\, a_3(\dddot {x}^{\mu}\,-(\ddot{x}^{\rho}\,\ddot{x}^{\sigma}\eta_{\rho\sigma})\,\dot{x}^{\mu})(\tau).
\end{align*}
The term $a_2 \ddot{x}^\mu$ combines with the left hand side to {\it renormalize} the mass
\begin{align}
(m_b -a_2)\ddot{x}^{\mu}=\,m\ddot{x}^{\mu}.
\end{align}
The argument from Rohrlich is completed after realizing that in order to obtain the Lorentz-Dirac equation,
 one needs $a_3=2/3\,e^2$.
 The same equation is obtained
  if instead of searching for a term containing the whole piece
$\, a_3(\dddot {x}^{\mu}(\tau)\,-(\ddot{x}^{\rho}\,\ddot{x}^{\sigma}\eta_{\rho\sigma})\,\dot{x}^{\mu})(\tau) $,
one requirers that the right hand to be compatible with the relativistic
    Larmor's law \cite{Jackson, Rohrlich},
\begin{equation}
\dot{P}^{\mu}_{rad}(\tau)=\, \frac{2}{3}\,e^2 (\ddot{x}^{\rho}\,\ddot{x}^{\sigma}\eta_{\rho\sigma})(\tau)\,\dot{x}^{\mu}(\tau).
\label{larmor}
\end{equation}
In order to fulfill this constraint, the minimal piece required in the equation of motion of a charged particle is $-2/3 e^2 (\ddot{x}^{\rho}\,\ddot{x}^{\rho}\eta_{\rho\sigma})\dot{x}^{\mu}$.
The Schott term $\frac{2}{3}\,e^2\dddot{x}$ is a total derivative. It does not contribute to the averaged power emission of radiation. However, in Rohrlich's argument, the radiation reaction term and the Schott term are necessary, due to the kinematical constraints of $\eta$ \cite{Dirac}.

Rohrlich's argument provides the  Lorentz-Dirac equation in a short and elegant way, without introducing complicated integrations and balance relations. However,
there are some points that make the Rohrlich argument not completely satisfactory. One difficulty is related with the structure of the vector $Z^{\mu}(\tau)$. In principle one can add pieces with higher derivatives and there is a lack of  justification for the absence of such pieces. For instance, it is possible to introduce an additional term in the right hand side of the form
\begin{equation}
Y^{\mu}(\tau)=\, (B^{\mu}\dot{x}^{\nu}-\,\dot{x}^{\mu}B^{\nu})(\tau)\dot{x}^{\rho}(\tau)\eta_{\nu\rho},
\label{Yterm}
\end{equation}
The vector field $B^{\mu}(\tau)$ along $x(\tau)$ can be arbitrary. However, this new contribution can be written in the form \eqref{newcontribution}.
Therefore, Rohrlich's argument is incomplete and several questions arise:
\begin{itemize}
\item It is not such that $\dddot{x}=0$ in the whole interval $I$,

\item Why one needs to start with the Lorentz force in defining the procedure?

\item What is the origin of the additional terms like $Z^{\mu}(\tau)$ or ${Y^{\mu}}(\tau)$ to the Lorentz force in the equation of motion?

\item Why are not there  derivative terms higher than three?
\end{itemize}
One can consider an adapted (Frenet) frame to the curve. For curves embedded in $R^n$ with the canonical flat connection, this method works if $\{\dot{x}(s),\ddot{x}(s),...,x^{(4)}(s)\}$ is a frame along $x:[0,1]\to M$; singularities must be treated individually. Thus, we do not expect higher order derivatives than four for describing arbitrary vector fields along $x:I\to M$. Finally, in order to keep the equation compatible with Larmor's law, it is necessary that
\begin{align}
a_4=0.
 \end{align}
We have proved the following
\begin{proposicion}
In a four dimensional spacetime $(M,F)$, the only Lorentz covariant
differential equation such that
 \begin{itemize}

 \item It is compatible with the covariant Larmor's law,
  \item The kinematical constraint $\eta(\dot{x},\dot{x})=-1$ holds,
  \item The constant observable mass condition $\dot{m}=0$ holds,
  \item The external electromagnetic forces are not dissipative, $\eta(F^{ext}_L,\dot{x})=0$
  \end{itemize}
  is the Lorentz-Dirac equation.
  \label{characterizationofLorentzDirac}
\end{proposicion}

\subsection{Rohrlich's derivation of the Lorentz-Dirac equation in the framework of generalized higher order fields}

Let us consider Rohrlich's argument in the framework of generalized electromagnetic fields introduced in {\it section 4}. In particular, the field defined by equation (\ref{regularizedradiationfield}), when evaluated on the lift $^kx(s)$ of the world-line $x(\tau)$ of a charged particle is
\begin{align*}
\bar{F}^{rad}(x,\, \dot{x},\,\ddot{x},\,\dddot{x},...,x^{(k)})=\big(F^{rad}_{\mu\nu}(\,^3x)+ \Upsilon_{\mu \nu}(x,\, \dot{x},\,\ddot{x},\,\dddot{x},...,x^{(k)})\big)dx^{\mu}\wedge dx^{\nu}.
\end{align*}

In the previous subsection we provided a solution to the required compatibility with the equation \eqref{larmor} that was $F^{rad}(S)=F^D(S)$, $\Upsilon=0$. In this subsection we consider $F^{rad}(S)=0$ and we find a solution for $\lim_{^kx\rightarrow \,^kS}\Upsilon(\,^kx)$ which is consistent with (\ref{larmor}). The piece $\Upsilon_{\mu \nu}(x,\, \dot{x},\,\ddot{x},\,\dddot{x},...)\dot{x}^{\nu}$ must be formally like $Z^{\mu}(s)$.
Let us write a formal series for $\Upsilon$ contracted with $\dot{x}$ (using $\eta$ to down indices),
\begin{align*}
\Upsilon_{\mu \nu}(x,\, \dot{x},\,\ddot{x},\,\dddot{x},...)\dot{x}^{\nu}=\upsilon_1\dot{x}_{\mu}\,
+\upsilon_2\ddot{x}_{\mu}\,+\upsilon_3\dddot{x}_{\mu}\,+...
\end{align*}
The minimal choice of the coefficients in the expansion that much such compatibility are
\begin{align}
\upsilon_{1}=\,-\frac{2}{3}(e^2)\ddot{x}^{\rho}\,\ddot{x}^{\sigma}\eta_{\rho\sigma},\quad, \upsilon_2=0,\quad \upsilon_3=\,\frac{2}{3}(e^2),\quad \upsilon_k=0,\,\forall k\geq 4.
\end{align}
The value of $\upsilon_1$ is necessary to recover the standard radiation reaction term. The value of $\upsilon_2$ is zero by simplicity. Indeed, if one considers the electromagnetic mass originated by the Coulomb field, it will be compensate with a convenient $\upsilon_2$ term, producing a renormalization of the mass. $\upsilon_3$ has this value in order that the constraint $Z^{\mu}\,\dot{x}^{\rho}\eta_{\mu\rho}=0$ holds.

Apart from the arbitrary election $\upsilon_k=0,\,k\geq 4$, there are more extra terms that we can add. For instance, a term like $(\ref{Yterm})$ it also possible to write down in the equation of motion. If we do this, the equation of motion has the following form
\begin{align*}
m\,\ddot{x}^{\mu}=\,\,eF^{\mu}\,_{\nu}\,\dot{x}^{\nu}+\,\frac{2}{3} e^2\,\big(\dddot{x}^{\mu}+(\ddot{x}^{\rho}\,\ddot{x}_{\rho})\dot{x}^{\mu}\big)+\, \big(B^{\mu}\dot{x}^{\nu}-\,\dot{x}^{\mu}B^{\nu}\big)\dot{x}_{\nu}\,
\end{align*}
\begin{align}
+\big(C^{\mu}\ddot{x}^{\nu}-\,\ddot{x}^{\mu}C^{\nu}\big)\dot{x}_{\nu}
\,+\big(D^{\mu}\dddot{x}^{\nu}-\,\dddot{x}^{\mu}D^{\nu}\big)\dot{x}_{\nu},
\end{align}
with $B(x,\dot{x},\ddot{x},\dddot{x},...)$, $C(x,\dot{x},\ddot{x},\dddot{x},...)$ and $D(x,\dot{x},\ddot{x},\dddot{x},...)$ elements of the jet bundle $J^k_0(M)$.
This equation can be written as
\begin{align*}
m\,\ddot{x}^{\mu}&\,=\,eF^{\mu}\,_{\nu}\,\dot{x}^{\nu}+\, P^{\mu}\,_{\nu}\,Z^{\nu}\,+\, \big(B^{\mu}\dot{x}_{\nu}-\,\dot{x}^{\mu}B_{\nu}\big)\dot{x}^{\nu}\,+\big(C^{\mu}\ddot{x}^{\nu}-\,\ddot{x}^{\mu}C^{\nu}\big)\dot{x}_{\nu}
\\
& +\big(D^{\mu}\dddot{x}^{\nu}-\,\dddot{x}^{\mu}D^{\nu}\big)\dot{x}_{\nu}.
\end{align*}
However, the orthogonality condition $\eta(\ddot{x},\dot{x})=0$ implies that
\begin{align}
B^\mu=0,\quad C^\mu=0,\quad, D^\mu=0.
\end{align}
Therefore, in the same conditions than in \ref{characterizationofLorentzDirac}, one has
\begin{proposicion}
With the same hypothesis than in \ref{characterizationofLorentzDirac},
 for generalized higher order fields, the differential equation for a point charged particle is the Lorentz-Dirac equation.
\end{proposicion}
This result implies that only using  generalized higher order fields is not enough
 to obtain a second order differential equation for point charged particles.

\subsection{Rohrlich's derivation of the Lorentz-Dirac equation with maximal acceleration}

Using the kinematical constraints for metrics of maximal acceleration, we can repeat Rohrlich's argument to obtain the Lorentz-Dirac equation. It will have a small modification due to the bound in the acceleration. Indeed, using the same notation as in {\it section 4} and with $\Upsilon=\,\sum^k_{i=1}\, \lambda_i\,x^{(i)}$, one obtains
\begin{align*}
Z^{\mu}(\tau)=P^{\mu}\,_{ \nu}(\tau)\big(\lambda_1\dot{x}^{\nu}(\tau)+\lambda_2\,\ddot{x}^{\nu}(\tau)+\lambda_3\,\dddot{x}^{\nu}),
\end{align*}
and by an analogous procedure as before, $\lambda_1$ is arbitrary and we can prescribe $\lambda_1=0$. The equation is consistent with Larmor's law if
\begin{equation}
\frac{\dot{\epsilon}}{2}\lambda_{2}+\,\lambda_3\frac{d}{d\tau}(g(\dot{x},\ddot{x}))=0,\quad \lambda_3=\,
\frac{2}{3}\,e^2,\quad \lambda_k=0,\quad\forall k\geq 4.
\end{equation}
The corresponding modified ALD equation is
\begin{equation}
m\,\ddot{x}^{\mu}=\,eF^{\mu}\,_{\nu}\,\dot{x}^{\nu}+\,\frac{2}{3} e^2\,\big(\dddot{x}^{\mu}-
(\ddot{x}^{\rho}\,\ddot{x}^{\sigma}\eta_{\rho\sigma})\dot{x}^{\mu}\big)+\,
\mathcal{O}({\epsilon^2_0}),
\end{equation}
with $F^{\mu}\,_{\nu}:=\,\eta^{\mu\rho}F_{\rho\sigma}$.
This equation is formally identical to the ALD equation,
\begin{proposicion}
With the same hypothesis than in \ref{characterizationofLorentzDirac} with a maximal acceleration geometry, the differential equation for a point charged particle is the Lorentz-Dirac equation.
\end{proposicion}
 This fact implies that only maximal
acceleration hypothesis is not enough to solve the {\it problem} of the Schott term in the ALD equation.
We have not considered the electrostatic contribution to the mass coming from the Coulomb field. If one wants to consider such contribution, an additional contribution $\lambda_2\neq 0$ and a renormalization of the mass to eliminate the divergence that appears are needed.

\section{A differential equation for point charged particles}

In this {\it section} we obtain a relativistic dynamical model for point charged particles. The dynamics of the probe charged particle will be described by an implicit second order differential equation compatible with Larmor's covariant formula $(\ref{larmor})$. The fact that we impose the requirement of being a second order differential equation is to be in accordance with Newton's first and second law. This requirement is a mayor difference with standard approaches to the electrodynamics of classical charged particles. Also, we do not require of the Lorentz-Dirac as first step in our theory. This is in sharp contrast the Landau-Lifshitz theory, where one finds an equivalent second order differential equation from the starting Lorentz-Dirac equation.

We will need to introduce maximal acceleration geometry. The justification to introduce such geometric framework is to be consistent with our philosophy of generalized higher order fields. Thus, if physical fields are generalized higher order fields and gravity should couple to them by a generalization of Einstein's equations, then the spacetime metric structure should be in the same category of generalized higher order fields. Moreover, metrics of maximal acceleration provide a perturbative parameter. Also, since in such frameworks the acceleration is bounded, it is natural that there are no pre-accelerated solutions.
\subsection{Derivation of the equation of motion}

Let us assume that the spacetime metric is of the type of maximal acceleration $g$ and that the metric $\eta$, obtained by averaging, is the Minkowski metric. We will perform our calculations in a normal coordinate system $(x,U_n)$of the Minkowski metric $\eta$. In such coordinate system the metric $\eta$ is globally constant and therefore
\begin{align*}
N^{\mu}\,_{\nu}=0,\,\,\mu,\nu=1,...,4,\quad x\in\,U_n.
\end{align*}
This is not a fundamental condition and will make easier some computations.
  Also, let us assume that the physical world-line of a point charged particle is a smooth curve of class $\mathcal{C}^k$ such that $g(\dot{x},\dot{x})=-1$, $\dot{x}^1>0$ and such that the acceleration field is bounded from above. Using the generalized tensor fields introduced in {\it section 2}, one obtains the following general form for the differential equation,
\begin{align*}
\Upsilon_{\mu\nu}(x,\,\dot{x},\,\ddot{x},\,\ddot{x},...)=
\,B_{\mu}\dot{x}_{\nu}\,-B_{\nu}\dot{x}_{\mu}\,+
C_{\mu}\ddot{x}_{\nu}\,-C_{\nu}\ddot{x}_{\mu}\,
+D_{\mu}\dddot{x}_{\nu}\,-D_{\nu}\dddot{x}_{\mu}+...\,\quad \dot{x}_{\mu}=\,g_{\mu\nu}\,\dot{x}^{\nu}.
\end{align*}
This implies that we will have the expression
\begin{align*}
m_b\,\ddot{x}^{\mu} & =\,e\,F^{\mu}\,_{\nu}+ \big(B^{\mu}\dot{x}_{\nu}-\,\dot{x}^{\mu}B_{\nu}\big)\dot{x}^{\nu}\\
& +\big(C^{\mu}\ddot{x}_{\nu}-\,\ddot{x}^{\mu}C_{\nu}\big)\dot{x}^{\nu} +\big(D^{\mu}\dddot{x}_{\nu}-\,\dddot{x}^{\mu}D_{\nu}\big)\dot{x}^{\nu}\,+...,
\end{align*}
and with $F^{\mu}_{\nu}:=\,g^{\mu\rho}F_{\rho\sigma}=\eta^{\mu\rho}F_{\rho\sigma}$.
 On the right hand side of the above expression all the contractions that appear in expressions like $\big(B^{\mu}\dot{x}_{\nu}-\,\dot{x}^{\mu}B_{\nu}\big)\dot{x}^{\nu}$, etc,  are performed with the metric $g$ instead of the Minkowski metric $\eta$.
The term $\delta m\,\ddot{x}^{\mu}$ corresponds to the electrostatic mass due to the Coulomb force \cite{Jackson}.
The other terms come from the higher order terms of the expression $(\ref{electromagneticfield})$ of the electromagnetic field.

The general form of a $k$-jet along a smooth curve $x:R\to M$ implies the relations
\begin{align*}
& B^{\mu}(s)=\,\beta_1\,\dot{x}^{\mu}(s)\,+\,\beta_2\,\ddot{x}^{\mu}(s)\,+
\,\beta_3\,\dddot{x}^{\mu}(s)\,+\,\beta_4\,\ddddot{x}^{\mu}(s)+\cdot\cdot\cdot,\\
& C^{\mu}(s)=\,\gamma_1\,\dot{x}^{\mu}(s)\,+\,\gamma_2\,\ddot{x}^{\mu}(s)\,+
\,\gamma_3\,\dddot{x}^{\mu}(s)\,+\,\gamma_4\,\ddddot{x}^{\mu}(s)+\cdot\cdot\cdot,\\
& D^{\mu}(s)=\,\delta_1\,\dot{x}^{\mu}(s)\,+\,\delta_2\,\ddot{x}^{\mu}(s)\,+
\,\delta_3\,\dddot{x}^{\mu}(s)\,+\,\delta_4\,\ddddot{x}^{\mu}(s)+\cdot\cdot\cdot.
\end{align*}
In principle we can add higher derivative terms. However,
the treatment of them will be the same than the fourth derivative term and eventually all of them will vanish. Also, for a $4$-dimensional smooth manifold $\{\dot{x},\ddot{x},\dddot{x},\ddddot{x}\}$ defines a local frame along the curve $x:I\to M$, except for the singularity points where some of the above derivatives are zero. We will assume that singularities only happens for isolated points. Otherwise, one needs to make a separate analysis (like we will see later for the uniform covariant motion).

Using the above expressions, one obtains
\begin{align*}
m_b\,\ddot{x}^{\mu} & =\, \big((\,\beta_1\,\dot{x}^{\mu}(s)\dot{x}_{\nu}\,+\,\beta_2\,\ddot{x}^{\mu}(s)\dot{x}_{\nu}\,+
\,\beta_3\,\dddot{x}^{\mu}(s)\dot{x}_{\nu}\,+\,\beta_4\,\ddddot{x}^{\mu}(s)\dot{x}_{\nu})\dot{x}^{\nu}\\
& -\,(\,\beta_1\,\dot{x}_{\nu}(s)\,+\,\beta_2\,\ddot{x}_{\nu}(s)\,+
\,\beta_3\,\dddot{x}_{\nu}(s)\,+\,\beta_4\,\ddddot{x}_{\nu}(s))\dot{x}^{\mu}\big)
\dot{x}^{\nu}\\
& + \big((\,\gamma_1\,\dot{x}^{\mu}(s)\,+\,\gamma_2\,\ddot{x}^{\mu}(s)\,+
\,\gamma_3\,\dddot{x}^{\mu}(s)\,+\,\gamma_4\,\ddddot{x}^{\mu}(s))\ddot{x}_{\nu})\dot{x}^{\nu}\\
& -\,(\,\gamma_1\,\dot{x}_{\nu}(s)\,+\,\gamma_2\,\ddot{x}_{\nu}(s)\,+
\,\gamma_3\,\dddot{x}_{\nu}(s)\,+\,\gamma_4\,\ddddot{x}_{\nu}(s))\ddot{x}^{\mu}\big)
\dot{x}^{\nu}\\
& + \big((\,\delta_1\,\dot{x}^{\mu}(s)\,+\,\delta_2\,\ddot{x}^{\mu}(s)\,+
\,\delta_3\,\dddot{x}^{\mu}(s)\,+\,\delta_4\,\ddddot{x}^{\mu}(s))\dddot{x}_{\nu})\dot{x}^{\nu}\\
& -\,(\,\delta_1\,\dot{x}_{\nu}(s)\,+\,\delta_2\,\ddot{x}_{\nu}(s)\,+
\,\delta_3\,\dddot{x}_{\nu}(s)\,+\,\delta_4\,\ddddot{x}_{\nu}(s))\dddot{x}^{\mu}\big)
\dot{x}^{\nu}.
\end{align*}
Let us assume that there are not derivatives higher than $2$ in the differential equation
of a point charged particle. One way to achieve this is to impose that all the coefficients for higher derivations are equal to zero,
\begin{align}
\gamma_k\,=\delta_k\,=0,\, k\geq 0,\quad \beta_k=0, k>3.
\end{align}
With this choice and using the kinetic relations for $g$, one obtains the expression
\begin{align}
m_b\,\ddot{x}^{\mu} =\,e\,F^{\mu}\,_{\nu}\,-\beta_2 \ddot{x}^{\mu}-\,\beta_3\,\dddot{x}^{\mu} -\,\frac{1}{2}\beta_2\,\dot{\epsilon}\,\dot{x}^{\mu} -\,\beta_3\,(-a^2(\tau)+\,\dddot{\epsilon})
\dot{x}^{\mu}.
\label{equationofmotionbeforerenormalization}
\end{align}
The differential equation governing the motion of a charge particle must be of second order and compatible with power radiation formula \eqref{larmor}.Thus, if $\dot{\epsilon}(\tau)$ is different than zero, one obtains the relations
\begin{align}
&\beta_2=\, \frac{4}{3}e^2\,a^2(s)\,\frac{1}{\dot{\epsilon}},\\
& \beta_k =0,\,\forall k\geq 3.
\end{align}
Therefore, at leading order in $\epsilon_0$ we obtain the differential
equation for a charged particle in a higher order field in the case that $a^2\neq 0$ and $(A^2_{max})^{-1}\neq 0$ to be
\begin{align*}
m_b\,\ddot{x}^{\mu} =\,\frac{2}{3}\,e^2\,a^2(s)\,\dot{x}^{\mu}-
\,\frac{2}{3}\,e^2\,a^2(s)\,\frac{1}{\dot{\epsilon}}\,\ddot{x}^{\mu}.
\end{align*}
There is a re-normalization of the {\it bare mass}.  For $\dot{\epsilon}\neq 0$ and the renormalization of mass reads,
\begin{equation}
 m_b+\,\frac{2}{3}\,e^2\,a^2(s)\,\frac{1}{\dot{\epsilon}}=\,m,\quad \dot{\epsilon} \neq 0.
\label{renormalizationofmass}
\end{equation}
If we add an external field $F^{\mu}\,_{\nu}$ interacting with the particle, we find the differential equation
\begin{equation}
m\,\ddot{x}^{\mu} =\,e\,F^{\mu}\,_{\nu}\,\dot{x}^{\nu}-
\,\frac{2}{3}\,{e^2}\,\eta_{\rho\sigma}\,\ddot{x}^{\rho}\ddot{x}^{\sigma}\,\dot{x}^{\mu},\quad F^{\mu}\,_{\nu}=\,g^{\rho\nu}F_{\rho\nu}.
\label{equationofmotion}
\end{equation}
In the case $a^2=0$, we postulate the same differential equation, which is equivalent to Newton's first law.

One can express the equation of motion (\ref{equationofmotion}) in a covariant way as
\begin{align}
m\,D_{\dot{x}}\,\dot{x}=\,e\widetilde{\iota_{\dot{x}}F}(x(s))\,-\frac{2}{3}\,e^2\eta(D_{\dot{x}}\,\dot{x},D_{\dot{x}}\,\dot{x}),
\label{covariantequationofmotion}
\end{align}
where $D_{\dot{x}}$ is the non-linear covariant derivative along $X=(\dot{x},0)\in T_{(x(s),\dot{x}(s))}TM$ and $\widetilde{\iota_{\dot{x}}F}=\,g^{-1}(\iota_{\dot{x}}F,\cdot)$.
We postulate that $(\ref{covariantequationofmotion})$ is the differential equation that the particle follows.

\begin{comentario} We have the following remarks:
\begin{itemize}
\item The factor $\beta_3$ is multiplying a factor $\dddot{x}^{\mu}-(a^2)\dot{x}^\mu$. However, we have seen that with the requirement of maximal acceleration, only $k=2$ is needed  to obtain a second order differential equation compatible with the power radiation formula. Note that although in the development in the derivatives $k=2$, the corresponding generalized higher order fields are sections of $\Lambda^2(M,\mathcal{F}(J^3_0(M)))$.

    \item Note that in the above derivation, $\dot\epsilon$ must be bounded and that the coefficient $\beta_2\in \,\mathcal{F}(J^3_0(M))$.

\item The derivation of equation  \eqref{equationofmotion} is not valid when $\dot{\epsilon}$ is zero. A separate discussion is necessary of that case. However, if the points where $\dot{\epsilon}$ is discrete, one can extend the validity of the differential equation by continuity.

\item The derivation of equation \eqref{covariantequationofmotion} and \eqref{equationofmotion} are independent of the Lorentz-Dirac equation. This is a point in faubour of the consistency of the approach, that is not relegated to problematic equations.

    \end{itemize}
\end{comentario}
If we assume that the maximal acceleration is infinite, then one has that $\dot\epsilon=0$ identically. In this case the relation \eqref{equationofmotionbeforerenormalization} reduces to
\begin{align}
m_b\,\ddot{x}^{\mu} =\,e\,F^{\mu}\,_{\nu}\,-\beta_2 \ddot{x}^{\mu}\,-\beta_3\,\dddot{x}^{\mu} +\,\beta_3\,a^2(\tau)\dot{x}^\mu
\dot{x}^{\mu}.
\end{align}
This is compatible with Larmor's law if
\begin{align*}
\beta_3=-\frac{2}{3}\,e^2.
\end{align*}
 In this case, the equation of motion should be the Lorentz-Dirac equation. This is the result from {\it Proposition} \ref{characterizationofLorentzDirac}. This argument makes explicit that both ingredients, maximal acceleration and generalized higher order fields are necessary to obtain a second order differential equation for point charged particles.

\subsection{General properties of the equation \eqref{equationofmotion}}

Let us consider a normal coordinate system for $\eta$.
Let us multiply equation $(\ref{equationofmotion})$ by itself and contract with the metric $g$. Using the kinetic relations of Proposition \ref{proposiciononkineticconstraints} one obtains
\begin{align*}
\,m^2\,a^2(1-\epsilon) &=\,e^2\,F^{\mu}\,_{\rho}\,\dot{x}^{\rho}F^{\nu}\,_{\lambda}\dot{x}^{\lambda}\,(1-\epsilon)\eta_{\mu\nu}
+(\frac{2}{3}\,e^2)^2\,(a^2)^2\,\dot{x}^{\mu}\dot{x}^\nu\,g_{\mu \nu}\\
 & -\,2e\frac{2}{3}e^2\,F^\mu\,_{\rho}\dot{x}^\rho\dot{x}^\nu(1-\epsilon)\eta_{\mu\nu}\\
& =\,(1-\epsilon)\big( F^2_L-\,\frac{1}{1-\epsilon}(\frac{2}{3}\,e^2)^2\,(a^2)^2\big).
\end{align*}
with the magnitude of the Lorentz force $F_L$ given by
\begin{align*}
F^2_L=\,\,e^2\,F_{\nu}\,^{\mu}F_{\mu}\,^{\rho}\,\dot{x}^{\nu}\,\dot{x}^{\lambda}\eta_{\lambda\rho}.
\end{align*}
\begin{proposicion}
For any curve solution of equation \eqref{equationofmotion} one has the following consequences,
\begin{enumerate}
\item The Lorentz force is always spacelike or zero,
\begin{align*}
F^2_L\,\geq 0.
\end{align*}
\item In the case the Lorentz force is zero, the magnitude of the acceleration is zero,
\begin{align*}
F^2_L=0\,\quad \Leftrightarrow \,a^2=0.
\end{align*}
\item If there is an external electromagnetic field, the acceleration is bounded by the strength of the corresponding Lorentz force.
\end{enumerate}
\label{propiedadesequationofmotion}
\end{proposicion}
{\bf Proof}.
 In the limit $\epsilon<<1$, one can re-write the expression
\begin{align*}
\,m^2\,a^2(1-\epsilon) =\,(1-\epsilon)\big( F^2_L-\,\frac{1}{1-\epsilon}(\frac{2}{3}\,e^2)^2\,(a^2)^2\big)
  \end{align*}
   as the following
\begin{align}
F^2_L=\,\frac{1}{1-\epsilon}(\frac{2}{3}\,e^2)^2\,(a^2)^2+\,m^2 a^2,
\label{valueofFL}
\end{align}
from which follows the three consequences.\hfill$\Box$

Solving the quadratic equation for $u=a^2$, one obtains
\begin{align}
u= \, 2m^{-2}C^{-1}\big(-1+\,\sqrt{1+F^2_L C}\big),\quad C=\frac{(\frac{4}{3} e^2)^2}{m^4}.
\label{boundoftheacceleration1}
\end{align}
In the limit $F^2_L \,C<<1$, one has
\begin{align*}
u=\,2\,m^{-2}C^{-1}\big(-1+\,\sqrt{1+F^2_L C}\big)\simeq \,2\,m^{-2}C^{-1}\big(-1+1+\frac{1}{2}F^2_L \,C\big)=\,m^{-2}\,F^2_L,
\end{align*}
in concordance with Newton second law of dynamics for the Lorentz force.

For the perturbative regime $\epsilon<< 1$ and also in the complementary limit when $F^2_L \,C>> 1$, one has that
\begin{align*}
u=\,2\,m^{-2}C^{-1}\big(-1+\,\sqrt{1+F^2_L C}\big)\simeq \,2\,m^{-2}C^{-1}F^2_L \,C\simeq \,2\,m^{-2}\,F^2_L,
\end{align*}
which is bigger by a factor $\sqrt{2}$ of the expected magnitude for the acceleration following the Lorentz law. It is not difficult to show that this asymptotic value is maximal,
\begin{teorema}
Given the Lorentz force of magnitude $F_L$, the maximal value attainable by particle following the equation of motion \eqref{equationofmotion} is
\begin{align}
a=\,\sqrt{2}\,m^{-1}\,F_L.
\label{boundinacceleration}
\end{align}
\end{teorema}
{\bf Proof}. First, note that in the perturbative limit and for forces of electromagnetic origin, $u\in\,]m^{-1}\,F_L,2\,m^{-1}\,F_L].$ This is because the function
\begin{align*}
S(F_L)=\,2\,m^{-2}C^{-1}\big(-1+\,\sqrt{1+F^2_L C}
 \end{align*}
 is monotonically increasing with $F_L$. Thus since we prove the existence of an asymptotic upper bound in the limit $F^2_L C\to \infty$, in the perturbative regime there is the bound for the acceleration \eqref{boundinacceleration}.\hfill$\Box$

One way to read this result is that for point charged particles, the perturbative regime is characterized by the bound $u\in\,]m^{-1}\,F_L,2\,m^{-1}\,F_L].$  If the bound is violated, the perturbative regime is not valid and as a consequence, equation \eqref{equationofmotion} is in principle, not valid.

\subsection*{Absence of run away solutions}

One can prove the following version of Dirac's asymptotic condition,
\begin{teorema}
For solutions of the equation \eqref{equationofmotion} it holds the following asymptotic condition,
\begin{align}
\lim_{\tau\to \infty} F_L(\tau)=0 \,\Rightarrow\,\lim_{\tau\to \infty}\,a^2=0.
\label{asymptoticondition}
\end{align}
\label{asymptoticonditiontheorem}
\end{teorema}
{\bf Proof}. From the equation \eqref{valueofFL} one has that the only non-negative solution for $u$ is just $u=0$. Then from Proposition \ref{propiedadesequationofmotion} follows the result.\hfill$\Box$

Run away solutions are solutions have the following peculiar behavior: even if the external forces have a compact domain in the spacetime, the charged particles follows accelerating without end.  Theorem \ref{asymptoticonditiontheorem} implies that equation \eqref{equationofmotion} is free of such pathological solutions,
\begin{corolario}
Equation $(\ref{equationofmotion})$ does not have {\it run-away} solutions.
\end{corolario}

An alternative way to see this is the following. Let us assume that the external field is zero for some $\tau>\tau_0$. Then the expression for the maximal acceleration is obtained again from equation (\ref{equationofmotion}),
\begin{displaymath}
-m^2\,\big(1-\frac{a^2}{A^2_{max}}\big)a^2=\,(\frac{2}{3}\,e^2)^2\,(a^2)^2.
\end{displaymath}
This condition can be read as
\begin{align*}
-m^2(1-\epsilon)\epsilon=\,A^2_{max}(\frac{2}{3})^2\epsilon^2.
\end{align*}
This implies the following consequence,
\begin{proposicion}
The equation \eqref{equationofmotion} does not have run away solutions iff
\begin{align}
A_{max}=\frac{3\, m}{2\,e^2}\,c^4.
\label{valueofthemaximalacceleration}
\end{align}
\end{proposicion}
Thus, since we have proved already that \eqref{equationofmotion} does not have run away solutions, we have
\begin{corolario}
The maximal acceleration of a point charged particle whose world-line of the equation \eqref{equationofmotion} is determined by the formula \eqref{valueofthemaximalacceleration}.
\end{corolario}
This acceleration is of the same order of magnitude than the maximal acceleration discovered by Caldirola \cite{Caldirola2}.

\subsection*{Absence of pre-acceleration for the equation \eqref{equationofmotion}}

In order to investigate the existence of pre-accelerated solutions of the equation $(\ref{equationofmotion})$, let us consider the example of a pulsed electric field \cite{Dirac}. This example conveys the discovery of the pre-accelerated solutions in the Lorentz-Dirac equation. For a electric pulse
\begin{align}
\vec{E}=(\kappa\,\delta(\tau),0,0),
\label{Diracpulse}
\end{align}
 the equation $(\ref{equationofmotion})$ in the non-relativistic limit reduces to
\begin{align}
a\ddot{x}^0=\kappa\,\delta(\tau),\quad a=\frac{3m}{2e^2}.
\end{align}
The solution of this equation is the {\it Heaviside function},
\begin{align}
a\dot{x} = & \kappa,\,\tau\geq 0,\\
& 0,\,\tau< 0.
\end{align}
This is not a pre-accelerated solution.
Therefore,
\begin{proposicion}
In the non-relativistic limit, equation $(\ref{equationofmotion})$ does not have pre-accelerated solutions of Dirac's type.
 \end{proposicion}
 Since the theory is Lorentz covariant, equation $(\ref{equationofmotion})$ does not have pre-accelerated solutions of Dirac's type in any other coordinate system. Still,
it is open the question if $(\ref{equationofmotion})$ is free of any other type of pre-accelerated solutions.
Also, note that the Dirac's pulse should be considered as an approximation, since we are considering smooth electromagnetic fields.

\subsection*{Compatibility with the covariant power radiation law}

The mechanical power performed by point charged particle whose world-line curve is a  solution of the differential equation $(\ref{equationofmotion})$ is obtained by the contraction of $m\ddot{x}^{\mu}$ with the four dimensional velocity vector $\dot{x}$,
\begin{align}
m\ddot{x}^{\mu}\dot{x}^\nu\,g_{\mu\nu}& =\,F^\mu\,_\rho\,\dot{x}^\rho\dot{x}^\nu\,g_{\mu\nu}\,-\frac{2}{3}e^2 a^2\, g(\dot{x},\dot{x})=\,\frac{2}{3}e^2 a^2.
\label{akinematicalconstraint}
\end{align}
From this expression it follows the compatibility of the equation \eqref{equationofmotion} with the covariant Larmor's law. The point charged particle takes such energy from the higher order terms of the electromagnetic field, since $\beta_2\,\dot\epsilon=-\frac{2}{3}e^2 a^2.$

Also interesting is to compute the non-relativistic limit.
There is a coordinate system where the acceleration is $a=(0,\vec{a})$ and the velocity vector is $\dot{x}=(v^0,\vec{v})$.
Two kinematical contractions $\ddot{x}^{\mu}\dot{x}^{\nu}g_{\mu\nu}=\,{1-\,\epsilon}\,\vec{a}\cdot\vec{v}$ and the contraction $g(\dot{x},\dot{x})=\,-(1-\epsilon)$. Then one obtains the rule,
\begin{align*}
m\,\vec{a}\cdot\vec{v}=\,-\frac{1}{1-\,\epsilon}\,\frac{2}{3}\,e^2 \,(\vec{a}\cdot\vec{a})\simeq \,-\frac{2}{3}\,e^2 \,(\vec{a}\cdot\vec{a})+\epsilon\,\frac{2}{3}\,e^2 \,(\vec{a}\cdot\vec{a}).
\end{align*}
Therefore,
\begin{align}
m\,\vec{a}\cdot\vec{v}=\,-\frac{2}{3}\,e^2 \,(\vec{a}\cdot\vec{a})\,+\mathcal{O}(\epsilon^2 _0).
\label{nonrelativisticlarmor}
\end{align}
 Equation $(\ref{nonrelativisticlarmor})$ is the mechanical power obtained from the forces acting on the particle by equation $(\ref{equationofmotion})$. In this sense, the non-relativistic limit of $(\ref{equationofmotion})$ is compatible the rate of energy loss by radiation. The same is true for the covariant equation of motion (\ref{covariantequationofmotion}).

\subsection*{Compatibility with kinematic constraints}

Using the kinetic relations from $g$, the relation \eqref{akinematicalconstraint} reduces to
\begin{align*}
m\,g(\ddot{x},\dot{x})=\,-\frac{2}{3}e^2 a^2(\tau)+\,\mathcal{O}(\epsilon^2).
\end{align*}
Defining the {\it characteristic time}
\begin{align*}
\tau_0:=\,\frac{2}{3m}e^2,
\end{align*}
the relation can be re-written as
\begin{align}
g(\ddot{x},\dot{x})\simeq -\tau_0\,a^2(\tau).
\label{conditionongandacceleration}
\end{align}
For accelerations taking place in a time much more larger than $\tau_0$, the condition \eqref{conditionongandacceleration} is a natural substitute to the orthogonal condition $\eta(\ddot{x},\dot{x})=0$. This interpretation also relates the maximal acceleration and the parameter $\tau_0$,
\begin{align*}
A_{max}\simeq \frac{c}{\tau_0}=\,\frac{3}{2}\frac{m\,c^4}{e^2},
\end{align*}
where $c$ is the speed of light in vacuum. Thus, this value is consistent with the exact value $A_{max}$ obtained by consistency with the requirement of absence of run away solutions.
\subsection*{Estimate of the maximal acceleration for a point charged particle}

The covariant Larmor's formula implies an estimate for the  maximal acceleration of a point charged particle. First, it is reasonable to assume that the maximal work that a point particle can realize is of the order of its rest mass $mc^2$. The characteristic time that this happens is of order of $\tau_0=r_0/c$, with $r_0$ the classical radius of the point particle. Therefore, the maximal power emitted in form of radiation,
\begin{align*}
 \big|{P}^{0}_{rad}(t)\big|=\, \frac{2}{3c^3}\,e^2 (\ddot{x}^{\rho}\,\ddot{x}^{\sigma}\eta_{\rho\sigma})(t)\,\dot{x}^{0}(t)<\,  \frac{2}{3c^3}\,e^2\,A^2_{max}=\frac{m \,c^2}{\tau_0}.
 \end{align*}
 Taking the value $(\tau_0)^{-1}= \frac{2m}{3c^3}\,e^2$, one obtains
 \begin{align*}
 A_{max}=\,\frac{3}{2}\frac{m\,c^4}{e^2}.
 \end{align*}
This estimated value coincides with the exact value for $A_{max}$ obtained by consistency in \eqref{valueofthemaximalacceleration}.
 From this value and \eqref{boundinacceleration} it follows the following bound for the value of the Lorentz force,
 \begin{align}
 F_L\leq \, \frac{1}{3}\,e^2\,m^2_e\,c^4.
 \end{align}

 Finally, note that for such maximal acceleration, the maximal power is the critical power $P_c$ introduced in {\it section 3},
 \begin{align}
 P(A_{max})=\,\frac{2}{3c^3}\,e^2 A^2_{max}=\,\frac{2}{3c^3}\,e^2\,\frac{3}{2}\frac{m\,c^4}{e^2}\, A^2_{max}=\, P_c.
 \end{align}
 Thus, since this value of the maximal acceleration is fixed by other arguments as do not have run-away solutions, we have that as consequence of {\it Proposition} \ref{conditionsforsuperluminicalmotion},
 \begin{proposicion}
There cannot be superluminical solutions of the differential equation \eqref{equationofmotion}.
 \end{proposicion}
 This is consistent with the argument in faubour of maximal acceleration based on the existence of a maximal speed of propagation for interactions and a minimal length. We conclude, that both speed and acceleration are bounded in our model of point charged particle.

\section{An effective spacetime electrodynamic theory for generalized higher order fields and point charged particles}

In this {\it section} we will introduce a generalization of Maxwell's equations for the generalized higher order fields $\bar{F}$, $\bar{G}$ and the generalized current $\bar{J}$. The metric $\eta$ is assumed flat and $M$ is four dimensional.
The theory uses an special type of cohomology elements of $\Lambda^*(M,\mathcal{F}(J^3_{0b}(M)))$, namely $\Lambda^*(M,\mathcal{F}(\tilde{J}^3_{0b}(M)))$. This is by the restriction of the general class of elements $\Lambda^*(M,\mathcal{F}(J^3_{0b}(M)))$ to be compatible with the solutions of the equation of motion \eqref{equationofmotion} (or its general covariant version \eqref{covariantequationofmotion}).
\subsection{Minimal extension of the generalized electromagnetic fields}

The {\it minimal extension} of the Faraday and excitation tensor fields  are of the form $(\ref{electromagneticfield})$ and $(\ref{excitationtensor})$. Except for very few terms, higher order terms do not contribute in the minimal extension theory. We assume that the constitutive relation $\bar{G}=\,\star \bar{F}$ holds. Since $\bar{G}$ and $\bar{F}$ are $2$-forms, the $\star$ operators for $\eta$ and $g=\,\lambda \eta$ when acting on sections of $\Lambda^{2}(M,\mathcal{F}({J}^3_0(M)))$ coincide.

Let us consider the $1$-forms
\begin{align*}
\,\{\widetilde{x^{(i)}}=\widetilde{x^{(i)}}_{\mu} \, d_4x^{\mu},\,\quad i=1,2,3\}.
\end{align*}
Then the electromagnetic and excitation fields
 $\bar{F}$, $\bar{G}\in\,\Gamma \Lambda^2(M,\mathcal{F}(J^3_0(M)))$
 can be expressed locally in the form
 \begin{align}
\bar{F}(x,\dot{x},\ddot{x},\dddot{x})=\big(\varphi(F)(x)+ \,\beta_2\,\tilde{\ddot{x}}\wedge\,\tilde{\dot{x}}\big),
\label{localexpressionforbarF}
\end{align}
\begin{align}
\bar{G}(x,\dot{x},\ddot{x},\dddot{x})=\star \,\big(\varphi(F)(x)+ \,\beta_2\,\tilde{\ddot{x}}\wedge\,\tilde{\dot{x}}\big),
\label{localexpressionforbarG}
\end{align}
where the $\star$ operator is associated with $\eta$.
The value of the coefficients in a normal coordinate system of $\eta$ is
\begin{align*}
\beta_2=\,\frac{2}{3}\,e^2\,A^2_{\max}\,\frac{(x^{(2)\rho}\,x^{(2)\lambda}\eta_{\rho\lambda})^{\frac{1}{2}}}
{2\,x^{(3)\rho}x^{(2)\lambda}\eta_{\rho\lambda}},
\end{align*}
Since $d_4$ is a skew-derivation, the following relation holds:
\begin{align}
d_4(\beta_i\,\widetilde{x^{(i)}})=\beta_i\,d_4(\epsilon)\wedge \,\widetilde{x^{(i)}}_{\mu} \, d_4x^{\mu},\,\quad i=1,2,3.
\label{equationofcoefficients}
\end{align}
The expression for $\Upsilon$ in a normal coordinate system is
\begin{align}
\Upsilon(\,^kx)=\,\frac{2}{3}\,e^2\,A^2_{\max}\,\frac{(x^{(2)\rho}\,x^{(2)\lambda}\eta_{\rho\lambda})^{\frac{1}{2}}}
{2\,x^{(3)\rho}x^{(2)\lambda}\eta_{\rho\lambda}}\widetilde{x^{(2)}}\wedge\,\widetilde{x^{(1)}}.
\label{equationforU}
\end{align}
\subsection*{Closeness of the form $\Upsilon$}
\begin{proposicion}
The form $\Upsilon\in\,\Gamma\,\Lambda^2(M,\mathcal{F}(J^3_0(M)))$ is closed,
\begin{align}
d_4\Upsilon =0,
\label{closerelationforupsilon}
\end{align}
\label{Propositiononcloserelationforupsilon}
 \end{proposicion}
 {\bf Proof}. It follows from the following calculation,
\begin{align*}
d_4\Upsilon (\,^kx)  & =\,d_4\Big(\frac{2}{3}\,e^2\,A^2_{\max}\,\frac{(x^{(2)\rho}\,x^{(2)\lambda}\eta_{\rho\lambda})^{\frac{1}{2}}}
{2\,x^{(3)\rho}x^{(2)\lambda}\eta_{\rho\lambda}}\widetilde{x^{(2)}}\wedge\,\widetilde{x^{(1)}}\Big)\\
& =d_4\Big(\frac{2}{3}\,e^2\,A^2_{\max}\,\frac{(x^{(2)\rho}\,x^{(2)\lambda}\eta_{\rho\lambda})^{\frac{1}{2}}}
{2\,x^{(3)\rho}x^{(2)\lambda}\eta_{\rho\lambda}}\Big)\,\widetilde{x^{(2)}}\wedge\,\widetilde{x^{(1)}}\\
& +\,\Big(\frac{2}{3}\,e^2\,A^2_{\max}\,\frac{(x^{(2)\rho}\,x^{(2)\lambda}\eta_{\rho\lambda})^{\frac{1}{2}}}
{2\,x^{(3)\rho}x^{(2)\lambda}\eta_{\rho\lambda}}\Big)d_4\Big(\widetilde{x^{(2)}}\wedge\,\widetilde{x^{(1)}}\Big).
\end{align*}
The first term is zero, since the operator $d_4$ acting on functions of higher order components is zero,
\begin{align*}
d_4\Big(\frac{2}{3}\,e^2\,A^2_{\max}\,\frac{(x^{(2)\rho}\,x^{(2)\lambda}\eta_{\rho\lambda})^{\frac{1}{2}}}
{2\,x^{(3)\rho}x^{(2)\lambda}\eta_{\rho\lambda}}\Big)\,\widetilde{x^{(2)}}\wedge\,\widetilde{x^{(1)}}=0.
\end{align*}
The second contribution is zero, because $d_4$ is nil-potent,
\begin{align*}
d_4\Big(\widetilde{x^{(2)}}\wedge\,\widetilde{x^{(1)}}\Big) & =\,d_4\widetilde{x^{(2)}}\wedge\,\widetilde{x^{(1)}}-\,\widetilde{x^{(2)}}\wedge\,d_4
\Big(\widetilde{x^{(2)}}\wedge\,\widetilde{x^{(1)}}\Big)\widetilde{x^{(1)}}\\
& = \,d_4({\epsilon})\wedge \,\widetilde{x^{(2)}}\wedge\widetilde{x^{(1)}}\,-\widetilde{x^{(2)}}\wedge\,d_4({\epsilon})\wedge\widetilde{x^{(1)}}\\
& =\,0.
\end{align*}\hfill$\Box$

The following are direct consequences of the algebraic structure of the generalized higher order fields $\bar{F}$ and $\bar{G}$.
\begin{corolario}
The generalized Faraday form $\bar{F}$ and the excitation form $\bar{G}$ are elements of $\Lambda^2_{cv}(M,\mathcal{F}(\tilde{J}^3_{0b}(M)))$. In particular, $[\bar{F}]\in \,H^2_{cv}(\tilde{J}^3_{0b}(M))$. Similarly, $[{\Upsilon}]\in \,H^3_{cv}(J^3_{0b}(M))$.
\end{corolario}
{\bf Proof}. It follows from the fact that for particles following the differential
equation \eqref{equationofmotion} have $n$-acceleration and the parameter
$\dot{\epsilon}^{-1}$ are bounded (if $\dot\epsilon\neq 0$). For $\epsilon=0$
it follows that $\bar{F}=0$ and $\bar{G}=0$.\hfill$\Box$
\begin{corolario}
Under the same assumptions, the fibers of $J^3_{0b}$ are such that
\begin{align*}
j^3_{0b}:=\{(x,\dot{x},\ddot{x},\dddot{x})\in j^3_{0}(x),\,s.t:
\end{align*}
\begin{itemize}
\item The non-compact condition $0<\eta(\dot{x},\dot{x})<1$ holds,

\item The non-compact condition $0<\eta(\ddot{x},\ddot{x})<A^2_{max}$ holds,

\item The compact condition $0<|\eta(\dddot{x},\dddot{x})|<c_3$ holds\}.
\end{itemize}
\label{corolarioonboundsofthefiber}
\end{corolario}
{\bf Proof}. It is clear that the domain for the first constraint is not compact.
 The domain for the second constraint , by the same reasons, the constraints $2$
 and $3$ do not have compact domain. \hfill$\Box$

This result shows how maximal acceleration implies constraints in the higher
 contributions to the generalized fields. Finally, it is natural to reconsider
 the isomorphism \ref{corolariosobreisomorfismodecohomologias}, since
  partially explains  the existence of an effective theory where the
  electromagnetic field is described by differential $2$-forms living
  on $M$, obtained by considering the averaged fields $\langle \bar{F}\rangle$ and $\langle \bar{G}\rangle$.

\subsection*{Electromagnetic vacuum}

The electromagnetic vacuum is characterized by the absence of electromagnetic field,
 $\bar{F}=0$. This is a {\it strong notion} electromagnet of vacuum. In this context, we have the following result,
\begin{proposicion}
If $\bar{F}=0$, then $\Upsilon=0$ and $\langle\bar{F}\rangle=0$.
\label{propertyofthegeneralized vacuum}
\end{proposicion}
{\bf Proof}. If $\bar{F}=0$, by identifying the corresponding components in equation
\eqref{localexpressionforbarF} with equation \eqref{estructuradebarF}, it follows
that $\varphi\langle \bar{F}\rangle=0$ and $\Upsilon=0$. \hfill$\Box$

The strong electromagnetic vacuum condition implies the weaker condition
\begin{align*}
\bar{F}=0\,\Rightarrow \,\langle \bar{F}\rangle=0.
\end{align*}
On the other hand, the weak vacuum condition does not imply the strong vacuum condition. However, it implies
\begin{align*}
\langle \bar{F}\rangle\,\Rightarrow\,\langle\Upsilon\rangle=0.
\end{align*}
The difference between the week and strong vacuum condition is  that in the weaker condition one can have $\Upsilon\neq 0$,
\begin{align*}
\frac{2}{3}\,e^2\,A^2_{\max}\,\frac{(x^{(2)\rho}\,x^{(2)\lambda}\eta_{\rho\lambda})^{\frac{1}{2}}}
{2\,x^{(3)\rho}x^{(2)\lambda}\eta_{\rho\lambda}}\widetilde{x^{(2)}}\wedge\,\widetilde{x^{(1)}}\,
\neq 0\quad \Rightarrow \quad x^{(2)\rho}\,x^{(2)\lambda}\eta_{\rho\lambda}\neq 0.
\end{align*}
Therefore, under the assumption that at infinity the electromagnetic state is the vacuum,
the weak vacuum state is not compatible with the validity of the equation of motion
\eqref{equationofmotion}. Thus if equation \eqref{equationofmotion} is valid for
$r\to \infty$, the weak condition is not enough and should be subjected to stronger constraints.

A similar argument shows that the strong vacuum condition implies
\begin{align*}
\bar{F}=0\quad \Rightarrow \frac{2}{3}\,e^2\,A^2_{\max}\,\frac{(x^{(2)\rho}\,x^{(2)\lambda}\eta_{\rho\lambda})^{\frac{1}{2}}}
{2\,x^{(3)\rho}x^{(2)\lambda}\eta_{\rho\lambda}}\widetilde{x^{(2)}}\wedge\,\widetilde{x^{(1)}}\,
= 0\quad \Rightarrow \quad x^{(2)\rho}\,x^{(2)\lambda}\eta_{\rho\lambda}=0.
\end{align*}
Thus the strong electromagnetic vacuum condition is compatible with physically
acceptable asymptotic properties of the equation \eqref{equationofmotion}.

\subsection{Maxwell's equations for higher order electromagnetic fields and currents}

One of the motivations to develop the mathematical machinery of generalized forms
in {\it section 2} and {\it section 4} is that it allows to generalize in an
straightforward way Maxwell's equations when written in covariant differential form. Thus,
we propose the following Maxwell's equations for $\bar{F}$ and $G=\star \bar{F}$. The homogeneous equations are
\begin{align}
d_4\bar{F}=0
\label{homogeneousequation}
\end{align}
Because equation $(\ref{equationofcoefficients})$ and (\ref{electromagneticfield}), the
homogeneous equation $(\ref{homogeneousequation})$ is equivalent to the standard homogeneous Maxwell's equations,
\begin{align}
d F=0.
\label{equationforF}
\end{align}
This is in accordance with the isomorphism in {\it Proposition}  \eqref{proposicionsobretildej},
\begin{align*}
H^*_{cv}(M,\mathcal{F}(\tilde{J}^3_{0b}(M)))\simeq\, H^{*-3n}_{dR}(M).
\end{align*}
This brings out the connection between the cohomology theory developed in {\it section} 2 and the theory of generalized higher order fields.

The non-homogeneous equations are
\begin{align}
d_4\star \,\bar{F}=\bar{J}:=\,\varphi({J})+d_4\,\star \Upsilon.
\label{nonhomogeneous equation}
\end{align}
For the non-homogeneous equation we note that
\begin{align*}
d_4\star\, \bar{F}=\,d_4 \star\,\varphi(F)\,+ d_4 \,\star\,\Upsilon=\,\varphi({J})\,+d_4\,\star \Upsilon.
\end{align*}
It is natural to define
\begin{align}
\bar{J}=\,\varphi(J)+d_4\,\star \Upsilon.
\label{JbarJ}
\end{align}
Therefore,  the current $J$ is such that
\begin{align*}
\varphi(J)=\,\bar{J}-\,d_4\star\,\big(\,\frac{2}{3}\,e^2\,A^2_{\max}\,\frac{(x^{(2)\rho}\,x^{(2)\lambda}\eta_{\rho\lambda})^{\frac{1}{2}}}
{2\,x^{(3)\rho}x^{(2)\lambda}\eta_{\rho\lambda}}\widetilde{x^{(2)}}\wedge\,\widetilde{x^{(1)}}\big).
\end{align*}
For such current, the non-homogeneous equations are equivalent to
\begin{align*}
d_4\star \varphi(F)=\,\varphi (J).
\end{align*}
Then one obtains an effective equation of the form
\begin{align}
d\star\, F=\,J,
\label{equationfor*F}
\end{align}
which are the standard non-homogeneous Maxwell equation.

If the current density $J(x)\in \Lambda^3M$ must be associated with physical systems, it must hold
\begin{align}
d{J}=0.
\label{equationfor J}
\end{align}
\begin{proposicion}
If (\ref{JbarJ}) holds, then $d_4\bar{J}=0\,\Leftrightarrow \,dJ=0$.
\label{relationbetweentheconservationsofcharges}
\end{proposicion}
{\bf Proof}. It is a consequence of the commutation relation $d_4(\varphi(\alpha))=\varphi(d\,\alpha)$,
because of the definition \ref{definiciondeformaintegral}.
Integrating on domains $\partial U\subset M$ one obtains the same fluxes and total charge for $\bar{J}$ and $J$.\hfill$\Box$
\subsection*{Boundary conditions for the generalized electromagnetic field}

 Let us pay attention on the boundary conditions for the generalized electromagnetic
 field $\bar{F}$. As in the discussion of the vacuum state before, the problem can be
  handle if one relates  boundary conditions of the field $\bar{F}$ with the boundary
  conditions for the field $F$. For example, if the boundary of the domain $D\subset$
  is $\partial D$ and the value on the boundary of the field is $F_0(x_0)$,
  then it is clear the following result,
\begin{proposicion}
Given an admissible boundary condition $F_0(x_0)$ for $F(x_0)$, there is a
corresponding admissible boundary condition for $\bar{F}_0$ given by
\begin{align*}
\bar{F}_0(\,^kx_0)=\,F_0(x_0)+\,\Upsilon(\,^kx_0).
\end{align*}
\label{boundaryconditions}
\end{proposicion}
In particular, one can apply this result to the initial valued problem, where the
Cauchy hypersurface $\Sigma \hookrightarrow M$ of the Maxwell's problem is lifted in $J^k_0(M)$.
However, for a prescribed field $F_0(x_0)$, the boundary condition is not necessarily
 unique, as the example of the asymptotic vacuum conditions showed.

\begin{teorema}
Let \eqref{localexpressionforbarF} and \eqref{localexpressionforbarG} be the generalized Faraday and excitation tensor fields.
If (\ref{equationforU}) holds, then the theory described by the system of equations
\eqref{homogeneousequation}, \eqref{nonhomogeneous equation}, \ref{conservationofcharge}
 and \eqref{equationofmotion} is equivalent to the theory described by the system of
 equations \eqref{equationforF}, \eqref{equationfor*F}, \eqref{equationfor J} and \eqref{equationofmotion}.
\label{teoremafinal}
\end{teorema}
{\bf Proof}. It is a direct consequence of the equivalence {\it Propositions}
\ref{closerelationforupsilon} and \ref{relationbetweentheconservationsofcharges}
and the equivalence between the possible boundary conditions for the corresponding set of equations.\hfill $\Box$

Therefore, the dynamics of generalized electromagnetic fields can be reduced to
 the standard Maxwell dynamics. The differential equation for the dynamics of
 point charged particles is described by the equation \eqref{equationofmotion}
  or the covariant version \eqref{covariantequationofmotion}.
\subsection*{Potential $1$-form for generalized higher order electromagnetic fields and gauge symmetry}

Let us introduce local gauge potentials for \eqref{localexpressionforbarF}.
From equation $(\ref{homogeneousequation})$, one can write
\begin{align*}
\bar{F}=\,d_4 \varphi(A)\,+\Upsilon,\quad A\in \,\Gamma\Lambda^1 M,
\end{align*}
 by the standard Poincar$\acute{e}$'s lemma. Then equation.  $(\ref{nonhomogeneous equation})$ can be expressed as
\begin{align*}
d_4\star(d_4 \varphi(A) +\,\Upsilon)=\,\varphi(J)\,\,+d_4\star\Upsilon
\end{align*}
with the potential $A$ satisfying the partial differential equation
\begin{align}
d_4\star d_4\,\varphi(A)=\varphi(J).
\label{equationforthepotential}
\end{align}
The $\star$ operator coincides for both metrics $\eta$ or $g$, since it is
applied acts  on two forms and $\eta$ and $g$ are related by a scalar factor.
Thus equation \eqref{equationforthepotential} is equivalent to the standard
equation as in Maxwell's theory for the electromagnetic potential $A$.
Therefore, one expects the existence and uniqueness of solutions of
equation \eqref{equationforthepotential}.

 As in standard electrodynamics, the potential is not fixed by the equations.
 Thus, it is necessary to consider an additional gauge fixing condition. For instance, the {\it Lorentz gauge condition} is
 \begin{align}
 \star_{\eta}\, d_4\star_{\eta}\, \varphi(A):=\,\delta_\eta\, \varphi(A)=0
 \label{lorentzgaugefixing}
  \end{align}
  is a valid gauge fixing condition (indeed, it corresponds to the usual
  Lorentz's condition for $\varphi(A)\in \,\lambda^1 M$).
Therefore,
\begin{corolario}
If the Lorentz gauge condition holds, $\varphi(A)$ is a solution of the wave equation
\begin{align}
\Box \varphi(A)=\,\star_{\eta}\,J.
\label{waveequation}
\end{align}
\end{corolario}
{\bf Proof}. This is direct from the definition of the wave operator
$\Box=\delta_\eta\, d+\,d\delta_\eta$ and the Lorentz gauge fixing.
\begin{comentario}
Note that we have chosen the star operator $\star_\eta$ in order to simplify the
corresponding wave equation and to have an admissible gauge fixing.
\end{comentario}
It is easy to prove that $(\ref{homogeneousequation})$ and $(\ref{nonhomogeneous equation})$
are invariant under the gauge transformation
 \begin{align*}
 A\to A\,+d\phi,
 \end{align*}
with $\phi\in \,\mathcal{F}((M))$. This gauge invariance corresponds with
the standard gauge invariance of the equations
(\ref{equationforF}) and (\ref{equationfor*F}).

\section{Discussion}

The theory presented in this work describes the motion in the vacuum  of
point charged particles interacting with the total electromagnetic field.
 It  intensively uses generalized higher order fields and the differential
 geometry theory of them. These fields are functionals acting on ordinary
 vector fields over the spacetime manifold $M$ but where the co-domain are
 spaces of functions defined on $J^3_0(M)$. The generalized higher order
 fields are used to define an extension of the notion field from standard
  fields to fields with dependence on higher derivatives of the curve
  describing the particle that  probes the field. Such extension provides
  additional degrees of freedom, allowing the fields to adapt the changes
  produced by the point charged particle used in the measurement of the field.
With such mathematical tool on hand, we were able to find an effective
electrodynamic theory (described by equations \eqref{homogeneousequation},
\eqref{nonhomogeneous equation}, \eqref{equationfor J} and  (\ref{equationofmotion}))
consistent with radiation and free from the pre-accelerated and run away pathologies
that plague the standard electrodynamic of classical point charged
particles\footnote{One still has to prove that there are not pre-accelerated
solutions of any type, not only of Dirac's type. In addition, the theory still
 uses of a mass renormalization procedure.}. Moreover, there is a general
 covariant version of the theory, where the differential equation of the point
 charged particle is the covariant version \eqref{covariantequationofmotion}.

The effective theory for the usual Faraday field $F$ and excitation
field $G$ living on $M$  is described by the equations $(\ref{covariantequationofmotion})$,
 $(\ref{equationforF})$, $(\ref{equationfor*F})$ and $(\ref{equationfor J})$.
Such effective theory is based on the relations \eqref{closerelationforupsilon},
 \eqref{localexpressionforbarF} and \eqref{equationfor J}.

\subsection{Discussion of the assumptions of the theory}

We have used of several assumptions in our considerations that because their novelty, we
should discuss and motivate them further. Highlighted below are the most relevant hypothesis considered:
\begin{itemize}

\item {\it Generalized higher order fields as fields with values in higher order jet bundles}.
This is the main new contribution of the theory. The notion is an hybrid between Mo-Papas' theory,
where dispersive forces were introduced and Wheeler-Feynman's theory. As we said in the introduction,
we do not abandon the concept of field as intermediate agent between charged particles. Also,
we think that use the language of forms is useful to describe the field, in which case the
electromagnetic force must be non-dissipative.

    Even if we persists on the physical reality of fields, we should admit that
    they have a quite {\it ghostly character} in our theory: they depend of the state
    of motion of the probe particle. In some sense, there is some similarity between
    a1scribing a physical reality to such objects is like ascribing physical reality
    to a quantum particle previous the measurement of a quantum observable. In
    the quantum case, we say that the quantum particle exists, but that its reality
    is not independent of the act of observation. Thus in a similar way, we say that
    the field exists, but its real value depend on how it is measured.

\item {\it Generalized higher order fields are horizontal.} This is one of the main
differences with other higher order field theories proposed in the literature
(see for instance \cite{Asanov, BucataruMiron, deLeonRodrigues, PfeiferWohlfarth, Vacaru, Voicu, VoicuSiparov}).
However, there are several reasons for this choice:
\begin{itemize}
\item If the generalized higher order fields are horizontal as defined in {\it section} 2,
there is a natural geometric notion of {\it local macroscopic measurement device}. After
an observer has been chosen, the measurement of a generalized higher order field  can be
magnified by using the effect on a {\it congruence of probe particles} instead of an
individual probe particle. Such {\it interaction} is mathematically described by the
flux of the generalized higher order field through the elementary $2$-dimensional element
of $TM$ determined by each pair of tangent vectors. This is an useful construction,
since the flux integral will increase the effect of the detection using a single
probe particle. On the other hand, the isomorphism \eqref{isomorphismcohomology}
implies that the flux is uniquely defined by the horizontal fields and viceversa.

\item It allows a causal description of the radiation reaction phenomena, in the sense
that changes in the fields in a region $U\subset M$ propagates with the speed of light
(thanks to equation \eqref{waveequation}). This property does not necessarily holds for
arbitrary (not necessary horizontal) fields. This is again clearly possible using
horizontal fields, by direct consequence of equation \eqref{waveequation}.

\item For horizontal fields it is easier to generalize Maxwell's equations, since one
 does not need to introduce additional dynamical equations for the non-horizontal pieces,
  which is at the present difficult to understand\footnote{In the Finslerian setting, a
  full framework of {\it Finslerian field theories} was developed in \cite{Asanov, VacaruStavrinos}.
  However, the uniqueness of such system of equations is far from being established and it needs of
   additional hypothesis and assumptions.}. This is again a consequence of the isomorphism \eqref{isomorphismcohomology}.
\end{itemize}

\item {\it Maximal $n$-acceleration and geometries of maximal acceleration.}
Maximal $n$-acceleration was introduced through the metric of maximal acceleration $g$.
It was assumed that the proper time is measured using $g$ and not the Minkowski metric $\eta$.
This is an important departing point from Special Relativity and General Relativity,
and implies that in our theory, the clock hypothesis does not hold. This is because
the metric of maximal acceleration depends on the acceleration of the particle. Thus,
this is implicit already in Caianiello's work, although seemingly not explicitly stated.

The introduction of geometries of maximal acceleration was useful to eliminate
run away solutions and to define a {\it perturbation scheme} in terms of the
parameter $\epsilon_0$. Since the asymptotic condition \eqref{asymptoticondition} holds,
 the book keeping parameter $\epsilon_0$ is related with the function
\begin{align}
\epsilon(\tau)=\,\frac{a^2(\tau)}{A^2_{max}}.
\label{perturbationparameter}
\end{align}
In the limit when the maximal acceleration $A^2_{max}$ is infinite, the asymptotic
expansions introduced in {\it section 5} are trivially zero. In such case, the procedure
used to obtain the  second order differential equation \eqref{equationofmotion} in
{\it section 6} fails. Thus, it is fundamental for our considerations to have
 a finite bound for the $n$-accelerations that a charged point particle can reach using electromagnetic fields.

We considered effective theories, where higher order terms greater than order one
in $\epsilon_0$ were disregarded. Thus the effective electromagnetic theory developed
in {\it section 7} is independent of the nature of the maximal acceleration, except
for the fact that it must be constant for each charged particle world-line $^kx:I\to M$.
However, we were able to provide a concrete value for the maximal acceleration, given by
 the formula \eqref{maximalacceleration}. The values of such acceleration for the electron,
 are of other $10^{32} m/s^2$, far away of any acceleration measured in the laboratory or
 observed in astrophysical objects. Thus the effective theory should be a
 good candidate for the current regime of experience.

Maximal $n$-acceleration has two fundamental effects in electrodynamics: from one
side, it provides a key ingredient to make the description of the motion of point
charged particles with radiation reaction consistent, since avoid run away solutions.
On the other hand, it implies the kinematics of maximal acceleration geometry, where
superluminical motion of accelerated massive particles is possible.
Moreover, the hypothesis of  maximal acceleration implies that is natural to consider
fields as sections of $T^{(p,q)}(M,\mathcal{F}(J^k_0(M)))$ with $k\geq 2$,
since generalized higher order field theory is linked with the existence of
a maximal $n$-acceleration through a generalized metric.

There are effective models of the electron where the book-keeping parameter is the
classical radios of the electron.  It could be the case that both approximations,
the existence of a universal bound in the acceleration and the existence of a finite
radius of the electron are related, since extended charged models requires a limit
for the acceleration, in order to preserve causality (see for instance \cite{Spohn2}.
 However, the maximal proper acceleration $A_{max}$ must be very large to any possible
 acceleration at the scales where the theory is applied, as it is the case of our theory.

\item {\it The equation of motion of a point charged particle must be of second order.}
This is a natural assumption, since it allows to maintain the {\it Principle of Inertia}.
Note that in the limit $A_{max}\to \infty$ the perturbative method used in {\it section 7}
 does not work. Indeed, in the framework of standard fields over $M$ and for Maxwell's
 electrodynamics, the equation of motion of a point charged particle must be Lorentz-Dirac
  equation, which is a third order differential equation. The Landau-Lifshitz equation,
  although second order, is not obtained from fundamental principles.

The equation of motion must be consistent with the rate of energy-momentum lost by the
emission of radiation. As a result, we find that (\ref{covariantequationofmotion})
(and therefore (\ref{equationofmotion})) is free of the pathologies of the Lorentz-Dirac equation.

\item {\it Point charged particle as ideal model of probe charged particle.}
We have assumed that the point charged particle is described by a one dimensional submanifold $x:I\to M$.
 However, one can should consider other classes of submanifolds in the more
 realistic description of the motion of the charge distributions. In this context,
 the use of distributional sources can be relevant \cite{Tucker}. It can happens that
 extended models could provide a solution to the still infinite Coulomb singularity.
\end{itemize}
\subsection{Phenomenological consequences of the theory}

We have found several phenomenological consequences headlined below:
\begin{itemize}
\item Consequences to zero order in the book-keeping parameter $\epsilon_0$,
\begin{itemize}
\item {\it An implicit second order differential equation for point charged particles}.
The constraint \eqref{covariantequationofmotion} is form invariant under local coordinate
transformations and it is free of run away and pre-accelerated solutions of Dirac's type.
The explicit solutions of that implicit equation are still to be fully explored, but in
principle the equation of motion provides a falsifiable model for the dynamics of point charged particles.

\item {\it Bound of the maximal covariant $n$-acceleration in the perturbative equation of motion}.
The $n$-acceleration $a^2=\eta(\ddot{x},\ddot{x})$ was predicted to have the
maximal value $A^2_{max}=\,(\frac{3 m c^4}{2 \, e^2})^2$. This maximal acceleration, for the electron, takes the value,
\begin{align}
A_{max}(e^-)=\,\frac{3}{2}\,\frac{m_e\,c^4}{e^2}=\,3.126\times 10^{32}\,m/s^2.
\label{maximalaccelerationfortheelectron}
\end{align}
This is an extraordinarily large value, which seems difficult to test in current
facilities \cite{FriedmanResin}\footnote{Several attempts to test in the laboratory
the principle of maximal acceleration has been tried by Y. Freedman and co-workers.
However, the values that the experiment try to check is of order $10^{19} m/s^2$.
Thus, the acceleration \eqref{maximalaccelerationfortheelectron} is not reachable in such attempts.}.

Independently of the value of the maximal acceleration $A_{max}$, when the equation
of motion \eqref{equationofmotion} is applicable, the $n$-acceleration must be
bounded by the value $\sqrt{2}m^{-1}\,F_L$.  This result is not in contradiction
with Newton's second law, since we know that radiation reaction can potentially
change the dynamical equation from the usual Newton's law. That this is a genuine
prediction can be seen by the fact that Lorentz-Dirac equation violates such a
bound (for instance, in run-away solutions).

\item {\it A bound on the maximal value attainable by the Lorentz force}. This is
obtained as combination of the approximate bound of the Lorentz force $F_L\leq \, \frac{1}{3}\,e^2\,m^2_e\,c^4$.

    A direct consequence is that the electric field must be also bounded,
    \begin{align}
    E_{max}=\frac{1}{3}\,e\,(m\,c^2)^2.
    \end{align}
    Notably, this bound is obtained in a linear theory, compared with the
    (undetermined) bound found in Born-Infeld electrodynamics \cite{BornInfeld}.

\item {\it Consequences for the speed of light}. In our theory, the speed of light in
vacuum is constant and it does not depend on the coordinate system.
This is a fundamental assumption. As a consequence, the theory is Lorentz invariant,
in the case that the base manifold $M$ is a four dimensional flat spacetime. Thus,
in the sense that any possible deviation from Lorentz invariance (in flat space) is not allowed.
\end{itemize}

\end{itemize}

The first three predictions does not depend on $A_{max}$ and are in principle falsifiable.
One can explicitly solved \eqref{equationofmotion} in several situations of interest
(Penning trap for example) and in principle the results can be contrasted with experiment.
For the third consequence-postulate, it is possible (in some weaker forms) to falsified or
test it (see for instance the relevant test of the speed of light in \cite{Mattingly}).

There are no predictions to higher order contributions in the parameter $\epsilon_0$,
since by construction all higher orders in $\epsilon$ were disregarded.

\subsection{Relation of the generalized higher order electrodynamics with other higher order field theories}

The use of generalized higher order fields in electrodynamics and other field theories has been previously investigated in the literature.
Without been exhaustive, some earlier work in the framework of Finsler or Finsler-Lagrange geometry
can be found in \cite{Asanov, BucataruMiron, Miron2004, Miron2006,  Vacaru, Voicu, VoicuSiparov}.
However,
distinctive from our approach is that
we have followed a minimal generalization,  from a particular motivation to introduce
generalized higher order fields, that was to find a solution to the radiation reaction
problem within a  classical theory of  electrodynamics. We have been able to prove the Theorem
\ref{teoremafinal} and show that such dynamics is free of most of the pathologies of the original theory.
Moreover, our theory does not have a direct relation with the Bopp and
Podolsky's theory, which is a higher order theory in the sense that contain partial
 differential equations with higher order derivatives than two, but with fields living on $M$.

\subsection{Future developments}

Some of the ideas presented in this paper have not been fully developed yet.
We would like to mention some of them that are relevant for the framework
presented in this work and that need further consideration. This includes
a development of the cohomology theory of generalized forms and an extended
investigation of the geometry of maximal acceleration. On the physical side,
the experimental predictions of the theory were not developed in fully detail.
Ranging in such spectrum, we headline below what could be further research
of the framework of generalized higher order fields:
\begin{itemize}
\item From the physical point of view, we find interesting to explore the following points.
\begin{itemize}
\item Experimental predictions of the theory presented in this paper.
Currently, this is a challenging point,
 since radiation reaction effects are very small in current experiments.
 However, the Penning-trap is a standard physical system where the equation
 of motion \eqref{equationofmotion} can be predictive. In particular, equation (\ref{equationofmotion})
can be solved for the Penning trap and compare it with the non-relativistic QED prediction \cite{BrownGravielse, Spohn1}.

\item We have shown that the geometry of maximal acceleration allows for
superluminical motion of massive point charged particles. However, that
this can happen for generalized metrics is not entirely surprising. For
instance, there are possible scenarios in Finsler spacetime models
(see for instance \cite{Kosteleky, PfeiferWohlfarth12}) where the spacetime metric allows that to happens.

\item Very strong constraints
for superluminical motion in neutrinos systems have been found by analyzing the
produced Cherenkov's radiation mechanism \cite{CohenGlashow11}. However, it is
also known that those constrains depend on the back-ground geometry (see for instance \cite{ChangWang}).
Thus, a natural question arises of which are the effects of Cherenkov's radiation in geometries
of maximal acceleration.

\item {\it Functional path integrals for the generalized electrodynamic theory}.
Note that a main difficulty formulating an action functional for a given field
theory is to define domain of the functional space where the fields are defined
and to provide it with a reasonable measure.  We can follow the example of functional integrals in the Finsler category, which
are integrals on the whole $TM$ of a lagrangian $L:TM\to R$ (see for instance  \cite{Voicu, PfeiferWohlfarth}).
 In a similar way, for
generalized higher order fields, one expects to have an action formulation of the theory as
\begin{align}
S[\Phi(\,^3x)]=\int_{J^3_{0b} (M)}\,dvol_4(x)\,\wedge dvol_V(y)\,\mathcal{L}(\Phi(\,^3x)).
\label{actionforphi}
\end{align}
Each fiber $j^3_{0b}$ is defined as in {\it Corollary} \ref{corolarioonboundsofthefiber}.
Thus, apart from  the constraint $\eta(\dot{x},\dot{x})$, there are additional
constraints with non-compact domains. This is an additional problem for the definition of the integrals.

 There are several natural questions related with the properties of the functional \eqref{actionforphi},
\begin{itemize}
\item  What are the Euler-Lagrange equations associated with \eqref{actionforphi}?
What is the general action such that its Euler-lagrange equations are the
(\ref{nonhomogeneous equation}) and (\ref{homogeneousequation}) equations?

\item It must be clarified the relation with the corresponding action functional of Maxwell's theory,
\begin{align}
S[A(x)]=\int_{M} \Big(-\frac{1}{2}\,F\wedge\,\star \, F\,+J\wedge\,A\Big).
\end{align}
Is it possible that $S[\Phi(\,^3x)]-S[\varphi(A)]$ is a total derivative in $J^3_0(M)$?
In the affirmative case,  the Euler-Lagrange equations will be equivalent.

\item To show the consistency of the quantization by functional integration, with generating functional given by
\begin{align}
Z[A]=\,\int \mathcal{D}\bar{A}\,exp\Big\{ \i(S[\Phi(\,^3x)(A)]+\,jA) \Big\}.
\end{align}
We believe that at least at the formal level, this tentative formulation can bring some light on the quantization process.
\end{itemize}

\item We have assumed that the electromagnetic medium is the standard vacuum. However, it seems natural to
extend the theory to more general electromagnetic media compatible with generalized
higher order fields and study the extension of the framework to models of electromagnetic media.

\item We have assumed that the probe particles are point charged particles and
described them by $1$-dimensional embedded manifolds on $M$ (world-line curves).
This is however, one of the possible source configuration. One should be able to
consider extended configurations, like sphere or charged planes. Although technically
 more complicated, it should be investigated in the view of potential applications in plasma physics.

\item The extension of the framework to Yang-Mills field models is a natural problem. Indeed, from a unification perspective,
this must be the case. In the Finslerian case such investigations has been performed for instance in \cite{VacaruStavrinos}.

\item
It is interesting to see if the generalized higher order fields correspond to a mathematical description {\it in between}
the usual local representation of gauge variables by the potential $A$ and the holonomy variables representation,
which is the natural variable for quantum gauge theory \cite{ChanTsou}.

\item We did not discuss gravitational effects. We think that a proper discussion of them should not require
additional methods. One direct issue is the following: we defined generalized
electromagnetic fields as sections of $\Lambda^2_{cv}(M,\mathcal{F}(J^3_0(M))$,
while the metric of maximal acceleration is a section of $T^{(0,2)}(M,\mathcal{F}(J^2(M))$.
Thus, if the matter fields like the generalized electromagnetic field $\bar{F}$ is related with the generalized metric $g$, the embedding
\begin{align}
T^{(0,2)}(M,\mathcal{F}(J^3(M))\hookrightarrow \,T^{(0,3)}(M,\mathcal{F}(J^2(M))
\end{align}
should play a relevant role.

\item Probably related with the above point is to provide a mechanism for maximal
acceleration. We suspect that the mechanism is of universal nature and that is of
such magnitude that accelerations related with the classical radius of the electron,
for instance, are small compared with such universal maximal acceleration. We suggest
that the mechanism is related with the quantum nature of the spacetime. Thus,
the resolution of this question will bring new insights on the problem of infinite self-energies.

\end{itemize}

\item From a mathematical point of view, we suggest the following research directions:
\begin{itemize}
\item A curvature theory for the non-linear connections on the spaces of jet bundles is
another direction that could be of interest to develop. This will include the generalization
of the fundamental results from Riemannian to generalized metrics.

\item In particular, a curvature theory for geometries of maximal acceleration is missing.
This can be introduced through the geodesic variation equation. More generally, it is of
relevance for both the theory and practice, to develop the variational theory of metrics of maximal acceleration.

\item We did not discuss in this work a natural formalism to describe  generalized higher order fields as it is the theory of sheaves
and sheaf cohomology \cite{Botttu, Warner}. We think that  the use of the sheaf cohomology
language can be useful to understand the structure of the theory, in
particular in the clarification of the relation of the several cohomologies introduced in {\it section 2}.

\item The relation between the vertical geometry and the theory of calibrations
\cite{HarveyLawson} has been already mention. This has direct relevance for our
research, since calibrations will be related with the way we perform the average of the higher order quantities.

\item There are many other problems open to mathematical investigation in the
geometry of generalized forms. We highlight here the generalization of Hodge
harmonic theory and the generalization of de Rham current theory to generalized forms.

\item A connection theory for metrics of maximal acceleration.

\item The development of a representation theory for the isometry groups of metrics of maximal acceleration.

\end{itemize}

\end{itemize}
\subsection*{Conclusion}

We have seen that in the framework of generalized higher order fields, it is possible to have a consistent theory of radiation-reaction of point charged particles at the classical level. It is also needed that the spacetime geometry is of maximal acceleration. Moreover, the theory can be extended to non-abelian theories and gravity. In a second step, one can consider the quantization of such fields theories.

\subsection*{Acknowledgements} We acknowledge to M. F. Dahl, J. Gratus, V. Perlick and J. P. Robinson
for valuable comments on previous versions of this work. This work is dedicated to Tun-I Hu for her continuous support.
Work financially supported by FAPESP, process 2010/11934-6 and the Riemann Center for Geometry and Physics, Leibniz University Hanover.
\small{

\end{document}
\begin{thebibliography}{22}

\bibitem{Asanov} G. S. Asanov, {\it Finsler
Geometry, Relativity and Gauge Theories}, D. Reidel, Dordrecht
(1985).

\bibitem{Bao} D. Bao, {\it On two curvature-driven problems
in Riemann–Finsler geometry}, Advanced Studies in Pure Mathematics 48,
Finsler Geometry, Sapporo 2005 - In Memory of Makoto Matsumoto
pp. 19–71 (2007).

\bibitem{BaoChernShen} D. Bao, S. S. Chern and Z. Shen, {\it An Introduction to Riemann-Finsler
Geometry}, Graduate Texts in Mathematics 200, Springer-Verlag (2000).

\bibitem{Beem1} J. K. Beem, {\it Indefinite Finsler Spaces and Timelike Spaces},
Canad. J. Math. {\bf 22}, 1035 (1970).

\bibitem{BeemEhrlichEasly} J. K. Beem, P. E. Ehrlich, K. L. Easly,
{\it Global Lorentzian Geometry}, Second Edition, CRC Press (1996).

\bibitem{Bonnor} W. B. Bonnor, {\it A new equation of motion for a radiating charged particle}, Proc. R. Soc. Lond. A. {\bf 337}, 591-598 (1974).
\bibitem{Bopp} F. Bopp, {\it Eine lineare theorie des Elektrons}, Ann. Physik {\bf 38}, 345-384 (1940).
\bibitem{Botttu} R. Bott and W. Tu, {\it Differential Forms in Algebraic Topology}, Springer-Verlag (1983).



\bibitem{BornInfeld} M. Born, L. Infeld, {\it Foundations of a new field theory}, Proc. R. Soc. Lond. A 144, 425 (1934).
\bibitem{BowickGiddins} M. J. Bowick and S. B. Giddins, {\it 
High-Temperature Strings}, Nucl. Phys. B {\bf 325}, 631 (1989).
\bibitem{Brandt1989} H. Brandt, {\it Maximal proprer acceleration and the structure of the space-time}, Found. Phys. Lett. 2, 39 (1989).
\bibitem{BrownGravielse} L.S. Brown and G. Gabrielse, Rev. Mod. Phys. 58, 233 (1986).

\bibitem{BucataruConstantinescuDahl} I. Bucataru, O. Constantinescu and M. F. Dahl,
{\it A geometric setting for systems of ordinary differential equations},
International Journal of Geometric Methods in Modern Physics, Vol. 8, Issue 6, pp. 1291-1327, (2011).

\bibitem{BucataruMiron} I. Bucataru and R. Miron, {\it Finsler-Lagrange geometry}, Editura Academiei Romane (2007).


\bibitem{Caianiello} E. R. Caianiello, {\it Is there a Maximal Acceleration}, Lett. Nuovo Cimento, 32, 65 (1981);
{\it Quantum and Other Physics as Systems Theory}, La Rivista 
del Nuovo Cimento, Vol. {\bf 15}, Nr 4 (1992).

\bibitem{Caldirola1956} P. Caldirola, {\it A new model of classical electron}, Suplemento al Nuovo Cimento, Vol. III, Serie X, 297-343 (1956).
\bibitem{Caldirola2} P. Caldirola, {\it On the existence of a maximal acceleration in the relativistic theory of the electron}, Lett. Nuovo. Cim. vol. 32, N. 9, 264 (1981).

\bibitem{ChanTsou} Chan Hong-Mo and Tsou Sheung Tsun, {\it Some Elementary Gauge Theory
Concepts}, World Scientific (1993).

\bibitem{ChangWang} Z. Chang, Xin Li, Sai Wang, {\it Neutrino superluminality without Cherenkov-like processes in Finslerian Special Relativity},
 Physics Letters B 710, 430-434 (2012).

\bibitem{CohenGlashow11}  A. G. Cohen and S. L. Glashow, {\it New Constraints on Neutrino Velocities},
 Phys. Rev. Lett. 107, 181803 (2011).

\bibitem{Crampin} M. Crampin, {\it On horizontal distributions on the tangent bundle of a
differentiable manifold}, J. London. Math. Soc. (2), 3, 178 (1971).

\bibitem{Crampin00} M. Crampin, {\it Connections of Berwals type}, Publ. Math. (Debrecen) 57 (2000) 455-473.
\bibitem{DengHou} Shaoqiang Deng, Zixin Hou {\it The group of isometries of a Finsler space}, Paci. Jour. Math. (207), No. 1 (2002).

\bibitem{deLeonRodrigues} M. de Le\'on, P. R. Rodrigues,
{\it Generalized Classical Mechanics and Field Theory}, North-Holland Mathematical Studies, 112 (1986).

\bibitem{Dirac} P. A. M. Dirac, {\it Classical Theory of Radiating Electrons},
Proc. R. Soc. Lond. A {\bf 137}, 148-169 (1938).

\bibitem{Einstein1922} A. Einstein, {\it The meaning of relativity}, Princenton University Press (1923).

\bibitem{Ehlers} J. Ehlers, {\it General Relativity and Kinetic Theory},
Proceedings of the International Summer School of Physics Enrico
Fermi, pg 1-70 (1971).

\bibitem{FeynmanWheeler} R. F. Feynman and J. A. Wheeler, {\it Interaction with absorvers as the mechanism of radiation}, Rev. Mod. Phys, {\bf 17}, 157 (1945).
\bibitem{FriedmanGofman} Y. Friedman and Y. Gofman, {\it A new relativistic kinematics of accelerated systems}, Physica Scripta 82, 015004 (2010).
\bibitem{FriedmanResin} Y. Friedman and E. Resin, {\it Dynamics of hydrogen-like atom bounded by maximal acceleration}, Physica Scripta, 86,  015002 (2012).

\bibitem{GallegoTorrome} R. Gallego Torrom\'e, {\it On a covariant version of Caianiello's Model}, Gen.Rel.Grav. { 39}, 1833-1845 (2007).

\bibitem{Ricardo09} R. Gallego Torrom\'e, {\it Fluid Models from Kinetic Theory
using Geometric Averaging}, Jour. Geom. Phys., 61 (2011) 829-846.

\bibitem{Ricardo012b} R. Gallego Torrom\'e, {\it A second order differential equation for point charged particles}, arXiv:1207.3627 [math-ph].

\bibitem{GirelliLiberatiSindoni} F. Girelli, S. Liberati, L. Sindoni,
{\it Planck-scale modified dispersion relations and Finsler geometry}, Phys. Rev.D 75:064015 (2007).

\bibitem{GoeckelerSchucker} M. Goeckeler, T. Schuecker,
{\it Differential geometry, gauge theories and gravity}, Cambridge University Press (1995).

\bibitem{GrallaHarteWald} S. E. Gralla, A. I. Harte, R. M. Wald,
{\it A Rigorous Derivation of Electromagnetic Self-force}, Phys. Rev. D 80: 024031 (2009).
\bibitem{HarveyLawson} R. Harvey, H. B. Lawson, {\it Calibrated geometries}, Acta Mathematica 148, 47–157 (1982).
\bibitem{HawEllis} S. W. Hawking and G. F. R. Ellis,
{\it The Large Scale Structure of the spacetime},
Cambridge Monographs on Mathematical Physics (1973).

\bibitem{Herrera1990a} L. Herrera, {\it The equation of motion for a radiating charged particle without
self interaction term}, Physics Letters A 145, 14-18, (1990).

\bibitem{Jackson} J. D. Jackson, {\it Classical Electrodynamics}, Third ed.
Wiley (1998).

\bibitem{KN} S. Kobayashi and K. Nomizu, {\it Foundations of Differential Geometry, Vol 
I}, Wiley Intersciencie, New York (1969).

\bibitem{KolarMichorSlovak} I. Kolar, P. W. Michor, J. Slovak, {\it Natural operators in differential geometry}, Springer-Verlag (1993).

\bibitem{Kosteleky} A. Kosteleky, {\it Riemann-Finsler geometry and Lorentz-violating kinematics}, Phys. Lett. B701, 137-143 (2011).

\bibitem{LandauLifshitz} L. D. Landau, E. M. Lifshitz, {\it The Classical Theory of
Fields}, Pergamon, Oxford (1962).




\bibitem{Mattingly} D. Mattingly, {\it Modern Tests of Lorentz Invariance},
Living Rev. Relativity 8 (2005), 5 http://www.livingreviews.org/lrr-2005-5.

\bibitem{Michor1980} P. W. Michor, {\it Manifolds of differentiable maps}, Shiva Publishing (1980).

\bibitem{Miron2004} R. Miron, {\it The geometry of Ingarden spaces}, Rep. on Math. Phys. 54(2), 131-147 (2004).

\bibitem{Miron2006} R. Miron, {\it Finsler-Lagrange Spaces with $(\alpha, \beta)$-Metrics and Ingarden Spaces},
Reports on Mathematical Physics, Vol. {\bf 58}, 417 (2006).

\bibitem{MoPapas} T. C. Mo and C. H. Papas, {\it New equation of motion for classical chraged particles}, Phys. Rev. D, Vol. {\bf 4}, 3566-3571 (1971).

\bibitem{MonizSharp} E. J. Moniz and D.H. Sharp, {\it Absence of runaways and divergent solutions in nonrelativistic Quantum Mechanics}, Phys. Rev. D 10, 1133 (1974); {\it Radiation reaction in nonrelativistic
quantum electrodynamics}, Phys. Rev. {\bf (15)}, 2850 (1977).

\bibitem{ParentaniPotting} R. Parentani and R. Potting, {\it Accelerating Observer and the Hagedorn 
Temperature}, Phys. Rev. Lett.{\bf 63}, 945 (1989).
\bibitem{PfeiferWohlfarth} C. Pfeifer and M. N. R. Wohlfarth, {\it Causal structure and electrodynamics on Finsler spacetimes},
Phys. Rev. D {\bf 84}:044039 (2011).

\bibitem{PfeiferWohlfarth12} C. Pfeifer, M. N. R. Wohlarth, {\it Beyond the speed of light on Finsler spacetimes}, Physics Letters B 712, 284-288 (2012).

\bibitem{Podolsky} B. Podolski, {\it A Generalized Electrodynamics I}, Nom-quantum, Phys. Rev. 62, 68-71 (1942).

\bibitem{Rindler} W. Rindler, {\it Relativity: Special, General and Cosmological}, Second Ed. Oxford Univ. Press (2006).

\bibitem{Rohrlich} F. Rohrlich, {\it The dynamics of a charge sphere of the electron}, Am. J. Phys. Vol. {\bf 65}, 1051-1056 (1997).


\bibitem{RovelliVidotto} C. Rovelli, F. Vidotto, {\it Maximal acceleration in covariant loop gravity and singularity resolution}, arXiv:1307.3228.
\bibitem{Saunders} D. J. Saunders, {\it The Geometry of Jet Bundles}, London Mathematical Society Lecture Note Series 142 (1989).



 \bibitem{Spohn1} H. Spohn, {\it The critical manifold of the Lorentz-Dirac equation}, Europhys. Lett. {\bf 50}, 287 (2000).

\bibitem{Spohn2} H. Spohn, {\it Dynamics of Charged Particles and Their Radiation Field}, Cambridge University Press (2004).

\bibitem{Toller} M. Toller, {\it Geometries of Maximal Acceleration}, hep-th/0312016; {\it La-
grangian and Presymplectic Particle Dynamics with Maximal Accelera-
tion}, hep-th/0409317.

\bibitem{Tucker} R. W. Tucker, {\it Differential form valued forms and
distributional electromagentic sources}, J. Math. Phys. {\bf 50},
033506 (2009).

\bibitem{Vacaru} S. I. Vacaru, {\it Nonholonomic distributions and gauge models of Einstein gravity}, Int. Journ. Geom. Meth. Mod. Phys. 7:215-246 (2010).

\bibitem{VacaruStavrinos} S. Vacaru, P. Stavrinos, E. Gaburon and D. Gonta, {\it Clifford and 
Riemannian-Finsler Structures in Geometric Mechanics and Gravity}, Geometry Balkan Press, 
2005, {Arxiv: gr-qc/0508023}.

\bibitem{Voicu} N. Voicu, {\it On electromagnetism and generalized energy-momentum tensor of the electromagnetic fields
in spaces with Finsler geometry}, arXiv:1012.2100 [math-ph].

\bibitem{VoicuSiparov} N. Voicu, S. Siparov, {\it A new approach to electromagnetism in anisotropic
spaces}, BSG Proc. 17, pp. 250-260 (2010).

\bibitem{Warner} F. Warner, {\it Foundations of Differentiable Manifolds and Lie Groups}, Scott, Foresman and Company (1971).

\end{thebibliography}
